\documentclass[11pt]{amsart}
\usepackage{amsmath,amsfonts,latexsym,amssymb,amscd, graphicx,pictexwd,enumerate}

\newcommand{\bpf}[1][Proof]{{\noindent {\sc #1: }}}
\newcommand{\epf}{{{\hfill $\Box$ \smallskip}}}
\newcommand{\ONE}{{\mathbf{1}}}
\newcommand{\N}{{\mathbb N}}
\newcommand{\Q}{{\mathbb Q}}
\newcommand{\Fc}{\mathcal{F}}
\newcommand{\Pp}{\mathsf{P}}
\newcommand{\Z}{{\mathbb Z}}

\newcommand{\RZ}{\R\times\Z}
\newcommand{\ZR}{\Z\times\R}

\newcommand{\E}{\mathsf{E}}

\newcommand{\R}{{\mathbb R}}

\newcommand{\HH}{{\mathbb H}}
\newcommand{\GG}{{\mathbb G}}

\newcommand{\Co}{{\rm Co}}

\newcommand{\calC}{{\mathcal C}}
\newcommand{\Ac}{{\mathcal A}}
\newcommand{\Bc}{{\mathcal B}}
\newcommand{\Cc}{{\mathcal C}}

\newcommand{\eps}{{\varepsilon }}

\newtheorem{theorem}{Theorem}[section]
\newtheorem{remark}{Remark}[section]
\newtheorem{lemma}{Lemma}[section]

\numberwithin{equation}{section}

\begin{document}
\title[Burgers equation]{Inviscid Burgers equation with random kick forcing in noncompact setting.}
\author{Yuri Bakhtin}
\maketitle

\begin{abstract}
We develop ergodic theory of the inviscid Burgers equation with random kick forcing in noncompact setting.
The results  are parallel to those in our recent work on the Burgers equation with
Poissonian forcing. However, the analysis based on the study of one-sided minimizers of the 
relevant action is different.
In contrast with previous work, finite time coalescence of the minimizers does not hold, and hyperbolicity
(exponential convergence of minimizers in reverse time) is not known.
In order to establish a One Force --- One Solution
principle on each ergodic component, we use an extremely soft method to 
prove a weakened hyperbolicity property and to construct
Busemann functions along appropriate subsequences.
\end{abstract}

\section{Introduction}
The Burgers equation was introduced by Burgers as a fluid dynamics model describing 
the evolving velocity profile of a 1D continuum of particles
in the absence of pressure. Although it did not prove to be a good model for most interesting
phenomena related to turbulence, it has re-emerged in various contexts including 
interface growth, traffic modeling,
the large-scale structure of the Universe, etc., see~\cite{Bec-Khanin:MR2318582} for a survey of the 
mathematics, physics, and numerics concerning the Burgers turbulence.

The inviscid Burgers equation is 
\begin{equation}
\label{eq:Burgers}
\partial_t u(t,x)+\partial_x\left(\frac{u^2(t,x)}{2}\right)=\partial_x F(t,x).
\end{equation}
In fluid dynamics terms,  $u(t,x)$ is the velocity of the particle located
at point $x\in\R$ at time $t\in\R$, and 
$f(t,x)=\partial_x F(t,x)$  is the external forcing term describing the acceleration of the particle at time $t$ at point $x$. 

Even if the forcing is absent or smooth,
solutions of this equation tend to develop discontinuities (shocks) in finite time. These shocks
correspond to collisions between particles and provide archetypal examples of shock waves.

Although the classical solutions are well-defined only locally due to
the inevitable singularity formation, 
one can work with generalized solutions. A natural global existence-uniqueness class of solutions is formed by so
called entropy solutions, or
viscosity solutions. They are globally well-defined and unique for a broad class of initial velocity profiles and forcing
terms. Moreover, they have a natural interpretation in terms of physics and admit a variational characterization.

\medskip

We are interested in the situation where the forcing is random, specifically, in
the ergodic properties responsible for the long-term statistics
 of solutions.  It is natural to expect
the existence of invariant distributions balancing the influx of the energy pumped in by the external forcing and
the dissipation of the energy at the shocks. This problem was considered first in compact settings, and a complete description
of invariant distributions on a circle or, in the multi-dimensional version of the problem, on a torus, was
obtained in~\cite{ekms:MR1779561},~\cite{Iturriaga:MR1952472},~\cite{Gomes-Iturriaga-Khanin:MR2241814}. The existence
and uniqueness of an invariant distribution for the Burgers dynamics with random boundary conditions were proved
in~\cite{yb:MR2299503}.

In these cases, the existence and uniqueness of an invariant distribution on the set
of velocity profiles with a given average
 followed from
the One Force --- One Solution Principle (1F1S) that asserts that for any $v\in\R$ and for almost every realization of the forcing
in the past, there is a unique velocity profile at the present that averages to $v$ and is compatible with the history of the forcing.
For any $v$, the collection of those velocity profiles indexed by all times $t\in\R$ forms a global solution that can be understood as
a random one-point pullback attractor.
Moreover, the distribution of this global solution at any fixed time is then a unique invariant distribution for the Markov semigroup generated
by the Burgers equation on velocity profiles averaging to $v$. 

In turn, to establish 1F1S, it is natural to employ a variational characterization of solutions
called the Lax--Oleinik variational principle. To introduce it, we first represent 
the velocity profile
as $u(t,x)=\partial_xU(t,x)$, where the potential $U(t,x)$
is a solution of the Hamilton--Jacobi--Bellman (HJB) equation
\begin{equation}
\label{eq:Hamilton-Jacobi-forced}
\partial_t U(t,x)+\frac{(\partial_x U(t,x))^2}{2}=F(t,x).
\end{equation}
 The viscosity solution of the Cauchy problem for this equation
with initial data $U(t_0,\cdot)=U_0(\cdot)$ can be written as
\begin{equation}
\label{eq:LO}
U(t,x)=\inf_{\gamma:
[t_0,t]\to\R}\left\{U_0(\gamma(t_0))+\frac{1}{2}\int_{t_0}^t\dot\gamma^2(s)ds+\int_{t_0}^t F(s,\gamma(s)) ds\right\},
\end{equation}
where the infimum is taken over all absolutely continuous curves $\gamma$
satisfying $\gamma(t)=x$. Then the solution $u$ of the Burgers equation can be found either by $u(t,x)=\partial_xU(t,x)$
or by using the slope of $\gamma^*$, the path on which the minimum in~\eqref{eq:LO} is attained:
$u(t,x)=\dot\gamma^*(t)$. The latter is
related to the fact that in the HJB equation~\eqref{eq:Hamilton-Jacobi-forced} the Hamiltonian is quadratic in $\partial_x U$. The sum of the last two terms on the right-hand side of~\eqref{eq:LO} is often called the (Lagrangian) action of $\gamma$. The minimizing path $\gamma^*$
can be identified with the trajectory of the particle that arrives to point $x$ at time $t$. For most points $(t,x)$ the minimizer is unique. However, there are exceptional locations corresponding to shocks where uniqueness does not hold.

This variational approach allows for an efficient analysis of the long term properties of the system via studying
the behavior of minimizers over long time intervals. The Burgers equation preserves mean velocity and is invariant with respect to
Galilean shear space-time transformations, so without loss of generality
let us confine ourselves to zero average velocity. 

To establish 1F1S, one has to prove that minimizers over increasing time 
intervals in the past converge to limiting trajectories, so called
one-sided infinite
minimizers. These are paths with infinite history such that every finite restriction of such a path minimizes Lagrangian action over paths
with the same endpoints. The entire space-time is foliated by these trajectories, and one can use their slopes to construct a global solution and prove that any global solution has to agree with this field of one-sided minimizers.

For this program to go through one needs an additional property of one-sided minimizers called hyperbolicity. It means that those paths approach
each other in the past sufficiently fast. The reason to consider this property 
is that in order to use~\eqref{eq:LO} in the proof that the velocity profile
obtained from slopes of the minimizers is, in fact, a global solution, one has to keep track not only of the velocity, but also of the velocity potential. To find the increment of the velocity potential, we can consider two minimizers approaching each other in reverse time. Although the action corresponding to each of them is infinite, one still can make sense of the difference in action between these paths since there will be a diminishingly small contribution from times in a distant past. This analysis leads to an analogy between the global
solution of the HJB equation~\eqref{eq:Hamilton-Jacobi-forced} and Busemann functions in last passage percolation theory.

For 1D periodic (or circle) setting hyperbolicity was first established in~\cite{ekms:MR1779561} and recently a simpler proof exploiting the
rigidity of 1D geometry was constructed in~\cite{Boritchev-Khanin}.

The first attempts of extending this program to noncompact settings, i.e., evolution
of velocity profiles on the entire real line with no periodicity or other compactness assumption were~\cite{Khanin-Hoang:MR1975784},
\cite{Suidan:MR2141893}, and~\cite{Bakhtin-quasicompact}. 
In~\cite{Bakhtin-quasicompact}, the program was carried out including the characterization of the global solution as a one-point random
attractor and a description of the domain of attraction. However, the random forcing still was mostly concentrated in a compact set, so 
one can describe this kind of setting as quasi-compact.

In~\cite{BCK:MR3110798}, the entire program was carried out for a fully noncompact space-time stationary setting on the real line. 
The forcing in \cite{BCK:MR3110798} was assumed to be concentrated in Poissonian points in space-time. The similarity of this problem
to last passage or first passage percolation and to Hammersley's process was exploited to adapt the methods 
of~\cite{Kesten:MR1221154}, \cite{HoNe3}, \cite{HoNe2}, \cite{HoNe},   \cite{Wu}, \cite{CaPi}, \cite{CaPi-ptrf} to this setting, although
several technical difficulties had to be overcome.

For every value $v\in\R$, the paper~\cite{BCK:MR3110798} constructs a unique family one-sided infinite minimizers with 
asymptotic slope $v$. The following strengthening of the hyperbolicity property was instrumental in constructing a global solution, proving its uniqueness
and attraction property: with probability $1$, any two of these one-sided minimizers coalesce, i.e., they meet at one of the Poissonian points
and coincide from that point on. This strengthening of the hyperbolicity property is certainly an artefact of the model where the forcing is concentrated in a discrete set of space-time points.
\smallskip

The goal of this paper is to go through the same program in a noncompact setting where the forcing is still space-time stationary, but it is applied at discrete times, so that it is smooth in space and delta-type in time. From the fluid dynamics perspective it means that in our model,
at every time $n\in\Z$
the velocity profile experiences a kick, i.e., the velocity of the particle at site $x\in\R$ is altered
by a random amount $F_\omega(n,x)$, where $F_\omega(n,x)$ is a space-time
stationary process. Between these kicks the velocity profile undergoes unforced Burgers evolution, i.e., each particle
travels with a constant velocity at least before it collides with other particles. For simplicity of the analysis
we will adopt a concrete ``shot-noise'' model for spatially smooth process $F_\omega(n,\cdot)$ such that these kicks at different times
are i.i.d. We note that kick forcing models have been studied for Burgers turbulence in compact setting in~\cite{Iturriaga:MR1952472},~\cite{Boritchev-Khanin} and for 2D Navier--Stokes system in \cite{Kuksin:MR1927233}, \cite{Kuksin:MR1845328}, \cite{Kuksin-book:MR2225710}. 

Although some steps in the program for the kicked forcing model are mere adaptations of the methods of the previous paper~\cite{BCK:MR3110798}, 
we have to tackle some new obstacles here. For example, some arguments that explicitly used the discrete nature of the Poisson point ensemble 
have to be enhanced or replaced with new ones.

One crucial difference between the old and new settings is that we can no longer rely on coalescence of one-sided minimizers. It is clear that 
distinct minimizers do not coalesce (all minimizers are solutions to the second-order Euler--Lagrange difference equation, so if 
two minimizers coincide with each other at two consecutive times, they are identical), 
but it is still not known if they are asymptotic to each other, i.e., 
if hyperbolicity holds true. This poses a serious difficulty. However, in this paper, we are able to replace hyperbolicity by a much weaker property. We prove that for any two 
one-sided minimizers $\gamma^1,\gamma^2$ with the same asymptotic slope, 
\begin{equation}
\label{eq:minimizers-approach-at-times}
\liminf_{m\to-\infty}\frac{|\gamma^1_m-\gamma^2_m|}{|m|^{-1}}=0.
\end{equation}
In other words, we prove that one-sided minimizers approach each other rather closely {\it along a sequence} $(m')$ of {\it pairing times}. It turns out that this weak hyperbolicity is sufficient to prove existence and uniqueness of a global solution of the Burgers equation in our model, and to study its domain of attraction. Of course, some limit transitions  $m\to-\infty$ are replaced by limits along the pairing sequence of times $m'$. In particular, our definition of Busemann function is based on partial limits.

The new argument we use to prove~\eqref{eq:minimizers-approach-at-times} (we actually prove and use a slightly stronger statement) is quite soft and can be
applied to other last passage percolation type models and random Lagrangian systems. However, we work only with the Burgers equation, since 
for this model system we can also complete all other steps of the program.

In the next two sections we explain the setting, main results, and the layout of the rest of the paper.

{\bf Acknowledgements.} I am grateful to Konstantin Khanin and Yakov Sinai who introduced me to the ergodic theory of the stochastic Burgers equation around the year of~2000. I am thankful to Eric Cator for educating me in last passage percolation. I thank Kostya and Eric for discussions at very early stages of this work. I would also like to thank the Banff International Research Station for Mathematical Innovation and Discovery where these discussions took place in July of 2012.

\section{The setting}
We will consider the dynamics under which at time $n\in\Z$ the velocity profile receives an instantaneous
random kick and then for time 1 solves the unforced Burgers equation before it receives the next kick, and so on. 
The forcing potential we want to consider can be informally written as
\[
\sum_{n\in\Z} F_\omega(n,x)\delta(t-n)
\]
where $\delta(\cdot)$ is Dirac's delta function and $F_\omega(n,\cdot)$, $n\in\Z$ is an i.i.d.\ sequence of stationary processes indexed by $x\in\R$ with finite range of dependence defined on a probability space~$(\Omega,\Fc,\Pp)$.

For simplicity, we work with a more concrete model where the potential $(F_\omega(n,x))_{(n,x)\in\Z\times\R}$ is given by a
shot-noise random field
\begin{equation*}
%\label{eq:shot-noise-definition}
F_\omega(n,x)=\sum_{i}\xi_{n,i} \phi\left(\frac{x-\eta_{n,i}}{\kappa_{n,i}}\right), \quad n\in\Z,\ x\in\R, 
\end{equation*}
where for each $n\in\Z$, $\{\eta_{n,i}\}_{i\in\N}$ is a Poisson point field on $\{n\}\times\R$ driven by the Lebesgue measure, 
$\phi:\R\to\R$ is a measurable
function with bounded support, and the amplitudes
 $(\xi_{n,i})_{i\in\N}$ and scaling factors $(\kappa_{n,i})_{i\in\N}$
are two bounded i.i.d.\ sequences, jointly independent of each other and the point configuration $\{\eta_i\}$. 
There is no canonical enumeration of Poissonian points, so let us now 
give a more precise and enumeration-independent definition of this model of point influences using the approach of
marked Poisson processes.

Let $\phi:\R\to\R$ be a differentiable even function such that the set~$\{x\in\R:\phi(x)\ne 0\}$ is non-empty and 
contained in $(-R_x,R_x)$ for some constant $R_x>0$.
Let $\{(\tau,\eta,\xi,\kappa)\}$ be a Poisson process on $\Z\times\R\times\R\times\R$ with driving measure
\begin{equation}
\label{eq:driving-measure}
\mu(dn\times dx\times dv\times du)=P_\tau(dn)\times P_\eta(dx) \times P_\xi(dv) \times P_\kappa(du),
\end{equation}
where $P_\tau$ is the counting measure on $\Z$, $P_\eta$ is the Lebesgue measure on~$\R$, 
$P_\xi$~and $P_\kappa$ are Borel probability  measures on $\R$ concentrated, respectively, on $[-R_\xi,R_\xi]$ and
on $(0, R_\kappa]$ for some positive constants $R_\xi,R_\kappa$.  To simplify the reasoning, we will assume that $R_x=R_\xi=R_\kappa=1$, 
although extending all our results to the arbitrary values of these constants is straightforward. Some arguments will be easier
if we assume that 
\begin{equation}
\label{eq:xi-has-values-if-both-signs}
P_\xi((0,R_\xi])>0,\quad \text{and}\quad  P_\xi([-R_\xi,0))>0,
\end{equation} 
although it is not necessary to make this assumption.

Let us recall the definition of a Poisson point field. 
It is convenient to assume that the probability space $\Omega_0$ is the space of locally finite 
point configurations $\omega$ on $\Z\times\R\times\R\times\R=\Z\times\R^3$. This space is equipped with $\sigma$-algebra $\Fc_0$
generated by maps $\omega\mapsto \omega(B)$, for bounded Borel sets $B\subset \Z\times\R^3$, where $\omega(B)$ denotes the number
of configuration points of $\omega$ in $B$. Then the probability measure $\Pp_0$ on $\Fc_0$ is defined via the following two properties:
(i) for any mutually disjoint bounded Borel sets $B_1,\ldots,B_k$, $\omega(B_1),\ldots,\omega(B_k)$ are independent random variables; (ii)
for a bounded Borel set $B$, $\omega(B)$ has Poisson distribution with parameter $\mu(B)$, where $\mu$ is defined in~\eqref{eq:driving-measure}.

The space-time projections $\{(\tau,\eta)\}$ of Poissonian $\{(\tau,\eta,\xi,\kappa)\}$ points form a Poisson point field in $\ZR$ driven
by $P_\tau(dn)\times P_\eta(dx)$. We will
often refer to these points as space-time footprints of the original configuration points.

Now we can define the random kick forcing potential via
\begin{equation}
\label{eq:shot-noise-definition}
F_\omega(n,x)=\sum_{\eta,\xi,\kappa: (n,\eta,\xi,\kappa)\in \omega}\xi \phi\left(\frac{x-\eta}{\kappa}\right),\quad (n,x)\in\Z\times\R.
\end{equation}
This sum is well-defined and differentiable for all $(n,x)$ because nonzero contributions come only from finitely many Poissonian points with
spatial component $\eta$ satisfying $|\eta-x|\le 1$. We will often omit the argument $\omega$ of $F(\cdot)=F_\omega(\cdot)$ for brevity.

It is immediate to see that 
$F(\cdot)$ is a space-time stationary process, i.e.,  the finite-dimensional distributions of $(F(n,x))_{(n,x)\in\Z\times\R}$
are the same as those of $(F(n+k,x+y))_{(n,x)\in\Z\times\R}$ for any choice of $(k,y)\in\Z\times\R$. Moreover,
$(F_\omega(n,\cdot))_{n\in\Z}$, is an i.i.d.\ sequence of stationary processes. Each of these stationary processes has radius of
dependence bounded by 2, i.e., $\sigma$-algebras generated by $(F(n,x))_{x\le r}$ and $(F(n,x))_{x\ge r+2}$ are independent for any choice of $r\in\R$, since the former depends on Poissonian points with spatial component $\eta$ satisfying $\eta\le r+1$, the latter depends
on Poissonian points with $\eta$ satisfying $\eta\ge r+1$, and the probability to have $\eta=r+1$ is zero. Let us also note that all
moments of 
\begin{equation}
F_\omega^*(n,x)=\max\{|F_\omega^*(n,y)|:\ y\in[x,x+1]\},
\label{eq:F-star}
\end{equation}
are finite and can be easily estimated by the corresponding moments of the Poisson distribution. Moreover,
\begin{equation}
\label{eq:tail-requirement2}
\varphi(\lambda)=\ln\E e^{\lambda |F_\omega^*(n,x)|}\in(0,\infty),\quad \lambda\in\R.
\end{equation}

It will be convenient in this paper to work on a modified probability space
\[
\Omega=\left\{\omega\in\Omega_0: \lim_{|x|\to\infty}\frac{F_\omega(n,x)}{|x|}=0,\quad n\in\Z\right\}\in\Fc_0.
\]
instead of $\Omega_0$. The reason is that $\Pp_0(\Omega)=1$ due to stationarity and moment assumptions
on $F_\omega$, and, as we will see later, 
on this set 
the Burgers dynamical system possesses some nice properties. Moreover, $\Omega$ is invariant
under (i) Galilean  space-time shear transformations $L^{a,v}$ moving 
each Poissonian point $(\tau,\eta,\xi,\kappa)$ to $(\tau,\eta+a+v\tau,\xi,\kappa)$, and 
 (ii) space-time translations $\theta^{n,x}$
moving each Poissonian point $(\tau,\eta,\xi,\kappa)$
to $(\tau-n,\eta-x,\xi,\kappa)$. Here $n\in\Z,$ $x\in\R,$ $v\in\R$.

We denote the restrictions of $\Fc_0$ and $\Pp_0$ onto $\Omega$ by $\Fc$ and $\Pp$.
From
now on for convenience we remove from $\Omega_0$ the zero measure complement to $\Omega$ and work with the probability
space $(\Omega,\Fc,\Pp)$. Under this transition, all the distributional properties are preserved.

\medskip

Let us now define the Burgers dynamics associated with this random forcing potential and the corresponding dynamics for HJB equation.
The latter may be called Hamilton -- Jacobi -- Bellman -- Hopf -- Lax -- Oleinik (HJBHLO) dynamics.

The space of velocity potentials that we will consider will be~$\HH$, the space of all locally
Lipschitz functions $W:\R\to\R$
 satisfying
\begin{align*}
 \liminf_{x\to\pm\infty}\frac{W(x)}{|x|}&>-\infty.
% \limsup_{x\to-\infty}\frac{W(x)}{x}&<+\infty.\\
\end{align*}

For a function $W\in\HH$, a forcing realization $F_\omega$ determined by $\omega\in\Omega$, times $n_0,n_1\in\Z$ satisfying $n_0<n_1$,
any a sequence of points $\gamma=(\gamma_{n_0},\ldots,\gamma_{n_1})$,
 we define the following
action:
\begin{equation}
\label{eq:action}
A_\omega^{n_0n_1}(W,\gamma)=W(\gamma_{n_0})+S^{n_0,n_1}(\gamma)+F_\omega^{n_0,n_1}(\gamma),
\end{equation}
where
\[
S^{n_0,n_1}(\gamma)=\frac{1}{2}\sum_{k=n_0}^{n_1-1}(\gamma_{k+1}-\gamma_k)^2
\]
is the kinetic action associated to $\gamma$, and the potential action
\begin{equation}
\label{eq:potential-action}
F_\omega^{n_0,n_1}(\gamma)=\sum_{k=n_0}^{n_1-1}F_\omega(k,\gamma_k)
\end{equation}
is responsible for the interaction with the external forcing
potential~$F_\omega$. Sometimes, we will identify the sequence of points $\gamma$ with the planar broken line consisting of
segments connecting $(n_0,\gamma_{n_0})$ to  $(n_1,\gamma_{n_1})$, $(n_1,\gamma_{n_1})$ to $(n_2,\gamma_{n_2})$, etc. Very
often we will omit the argument $\omega$ of  $A_\omega^{n_0n_1}(W,\gamma)$, $F_\omega^{n_0,n_1}(\gamma)$
and other random variables.

Let us now consider the following minimization problem:
\begin{align}
 \label{eq:optimization_problem}
 A_\omega^{n_0,n_1}(W,\gamma)&\to \inf,\\
\gamma_{n_1}&=x,\notag
\end{align}
where the infimum is taken over all sequences $\gamma:\{n_0,\ldots,n_1\}\to\R$.
Let us denote the
infimum value in~\eqref{eq:optimization_problem} by $\Phi^{n_0,n_1}_\omega W(x)$.

If $n_1=n_0+1$, then one can use 
the optimization problem~\eqref{eq:optimization_problem} to find the viscosity
solution of the unforced Burgers equation or the associated HJB equation with initial condition $W(\cdot)+F_\omega(n_0,\cdot)$, 
see, e.g.,~\cite[Section~3.4]{Evans-PDE-book:MR1625845}. Namely, $\Phi^{n_0,n_0+1}_\omega W(x)$ is well-defined 
and equal to the solution of the HJB at point $(n_0+1,x)$; and $\partial_x\Phi^{n_0,n_0+1}_\omega W(x)$ is the solution of the Burgers
equation at point $x$. The latter can be also represented as $x-\gamma_{n_0}$, where $(\gamma_{n_0},\gamma_{n_0+1})=(\gamma_{n_0},x)$
is the optimal two-point path solving~\eqref{eq:optimization_problem}.

 Note that in~\cite{Evans-PDE-book:MR1625845}, the
assumption of global Lipschitzness of the initial data is used. Although that assumption does not hold for 
$W(\cdot)+F_\omega(n_0,\cdot)$, it can be easily replaced by local Lipschitzness and at most linear growth to $-\infty$.
These properties hold for all $\omega\in\Omega$.

The general optimization problem~\eqref{eq:optimization_problem} with an arbitrary gap between $n_0$ and $n_1$ corresponds to
iterative applications of this one-step procedure, so the result matches the informal description of the dynamics given
in the introduction: at time $n_0$ we start with a velocity profile $W$,  alter it by $F(n_0,\cdot)$, then we solve the unforced
Burgers--HJB equation for time~$1$, and then, at time $n_0+1$ we alter the solution by $F(n_0+1,\cdot)$ and then solve unforced Burgers equation
for time~$1$, etc.

The family of random nonlinear
operators
$(\Phi^{n_0,n_1}_\omega)_{n_0\le n_1}$ is the main object in this paper. Our main goal is to understand the
asymptotics of $\Phi^{n_0,n_1}_\omega$ as $n_1-n_0\to\infty$. 
Let us now formulate the most important  preliminary  facts about the operators $\Phi^{n_0,n_1}_\omega$. We do not give proofs since they are just minor modifications of the textbook material in~\cite[Section~3.4]{Evans-PDE-book:MR1625845}.

\begin{lemma}\label{lem:properties-of-evolution-operator} For any $\omega\in\Omega,$ $W\in\HH$ and $n_0,n_1\in\Z$ with $n_0<n_1$,
the following holds true:
\begin{enumerate}
\item For any $x\in\R$ there is a path $\gamma(x)$ that realizes
the minimum in~\eqref{eq:optimization_problem}. In particular, the operator $\Phi_\omega^{n_0,n_1}$ is well-defined on $\HH$. 
\item \label{it:continuity_of_HJ_solution} The function $x\mapsto\Phi_\omega^{n_0,n_1}W(x)$ is locally Lipschitz.
\item The set $O$ of points $x$ such that the variational problem~\eqref{eq:optimization_problem} admits a unique solution is open
and dense in~$\R$. The complement to $O$ is at most countable.
\item If $x_0\in O$,  then $\Phi_\omega^{n_0,n_1}W(x)$
is differentiable at $x_0$ with respect to~$x$ and
\[
\partial_x\Phi_\omega^{n_0,n_1}W(x)\bigr|_{x=x_0}=\gamma_{n_1}(x_0)-\gamma_{n_1-1}(x_0).
\]
If $\gamma:[n_0,n_1]\to\R$ is the continuous curve linearly interpolating between points of the sequence $\gamma(x_0)$, then
 the right-hand side of this identity can be interpreted as $\dot\gamma(n_1)$, where the dot denotes (left) time derivative. 
 
\item If $x_0\in\R\setminus O$, then the right and left derivatives of
$\Phi_\omega^{n_0,n_1}W(x)$ w.r.t.~$x$
are well-defined at $x_0$. They are equal to the slope of, respectively, the leftmost and rightmost minimizers realizing
$\Phi_\omega^{n_0,n_1}W(x_0)$.
\item Let $\gamma(x)$ denote a minimizer with endpoint $(n_1,x)$. Then for every pair of points $x,y\in\R$ satisfying $x<y$, and every 
$n\in\{n_0,\ldots,n_1\}$, $\gamma_n(x)<\gamma_n(y)$. Moreover, $\lim_{x\to\pm\infty}\gamma_n(x)=\pm\infty$.
\end{enumerate}
\end{lemma}

Optimal sequences $\gamma(x)$ and their continuous interpolations can be viewed as particle trajectories.

The following statement is the cocycle property for the operator family~$(\Phi^{n_0,n_1}_\omega)$. It is a direct
consequence of Bellman's principle of dynamic programming.

\begin{lemma}\label{lem:cocycle} If $\omega\in\Omega_1$, then for any $W\in\HH$, any $n_0,n_1,n_2$ satisfying
$n_0<n_1<n_2$,
$\Phi_\omega^{n_1,n_2}\Phi_\omega^{n_0,n_1}W$ is well-defined and equals $\Phi_\omega^{n_0,n_2}W$. For any $x\in\R,$ if $\gamma$ is an optimal path
realizing $\Phi_\omega^{n_0,n_2}W(x)$, then the restrictions of $\gamma$ onto $\{n_0,\ldots, n_1\}$ and $\{n_1,\ldots,n_2\}$
are optimal paths realizing $\Phi_\omega^{n_0,n_1}W(\gamma_{n_1})$ and   $\Phi_\omega^{n_1,n_2}(\Phi_\omega^{n_0,n_1}W)(x)$.
\end{lemma}

Introducing $\Phi^n_{\omega}=\Phi^{0,n}_{\omega}$ we can rewrite the cocycle property as
\[
 \Phi^{n_1+n_2}_{\omega}W=\Phi^{n_2}_{\theta^{n_1}\omega} \Phi^{n_1}_\omega W,\quad n_1,n_2\in\N,\quad\omega\in\Omega,
\]
where $\theta^{n}=\theta^{n,0}$ denotes the time-shift on point configurations by $n$ time units.

Let us denote 
\[
\HH(v_-,v_+)=\left\{W\in\HH:\ \lim_{x\to \pm\infty} \frac{W(x)}{x}=v_\pm \right\},\quad v_-,v_+\in\R.
\] 
The following result shows that these spaces are invariant under HJBHLO dynamics.
Along with Lemma~\ref{lem:cocycle} it allows to treat the dynamics as
a random dynamical system with perfect cocycle property (see, e.g.,~\cite[Section~1.1]{Arnold:MR1723992}).
We give a proof of this lemma in Section~\ref{sec:aux}.

\begin{lemma} \label{lem:invariant_spaces} For any $\omega\in\Omega$,
for any $n_0,n_1$ with $n_0<n_1$,
\begin{enumerate}
 \item If $W\in\HH$, then $\Phi^{n_0,n_1}_\omega W\in\HH$.
 \item If $W\in\HH(v_-,v_+)$ for some  $v_-,v_+$, then $\Phi^{n_0,n_1}_\omega W\in\HH(v_-,v_+)$.
\end{enumerate}
\end{lemma}

In nonrandom setting the family of operators $(\Phi^{n_0,n_1})$ constructed via a variational
problem of type~\eqref{eq:optimization_problem}  is called a HJBHLO evolution
semigroup, see~\cite[Definition 7.33]{Villani},
but in our setting it would be more precise to
call it a HJBHLO cocycle.

\medskip

Potentials are naturally defined up to an additive constant. It is thus convenient to
work with $\hat \HH$, the space of equivalence classes of potentials from $\HH$. 
Also, we can introduce spaces $\hat\HH(v_-,v_+)$ 
as classes of potentials in $\HH(v_-,v_+)$ coinciding up to an additive constant.
The cocycle $\Phi$ can be
projected on $\hat\HH$ in a natural way. We denote the resulting cocycle on $\hat\HH$ by~$\hat\Phi$.

Since a velocity field determines its potential uniquely up to an additive constant, we can also introduce
dynamics on velocity fields. We can introduce the Burgers dynamics on the space $\HH'$ of functions $w$ (actually, classes of equivalence of
functions since we do not distinguish two functions coinciding almost everywhere) such that
for some function $W\in\HH$ and almost every $x$,  $w(x)=W'(x)=\partial_x W(x)$.
For all $v_-,v-+\in\R$, we can also introduce $\HH'(v_-,v_+)$, the space consisting of functions $w$ such that
the potential~$W$ defined by $W(x)=\int_0^x w(y)dy$ belongs to $\HH(v_-,v_+)$. One can interpret this space
as the space of velocity profile with well-defined one-sided averages $v_-$ and $v_+$.

 We will
write $ w_2=\Psi^{n_0,n_1}_\omega w_1$ if $w_1=W'_1$, $w_2=W'_2$, and $W_2=\Phi^{n_0,n_1}_\omega
W_1$  for some $W_1,W_2\in\HH$. Of course, the maps belonging to the Burgers cocycle $(\Phi^{n_0,n_1})$ map $\HH'$ into itself
and the spaces $\HH'(v_-,v_+)$  are also invariant. However, these maps have additional regularity
that we are about to exploit.

It follows from Lemma~\ref{lem:properties-of-evolution-operator}
that for all $\omega\in\Omega$ and every $w\in\HH'$ there is a {\it cadlag} version of $\Psi^{n_0,n_1}_\omega w$ (i.e., it is right-continuous
and has left limits at all points) such that
the function 
\[
M(x)= x-\Psi^{n_0,n_1}_\omega w (x),\quad x\in\R,
\] is strictly increasing and satisfies $\lim_{x\to\pm\infty}M(x)=\pm\infty$. 
In particular, $M$ has at most countably many discontinuities, and so does $\Psi^{n_0,n_1}_\omega w$.
In the fluid dynamics interpretation, the  particle that arrives 
at $x$ at time $n_1$ with velocity $\Psi^{n_0,n_1}_\omega w (x)$, is
located at $M(x)$ at time $n_1-1$. Monotonicity means that the paths of those particles do not cross before $n_1$.
Each discontinuity of $M$ corresponds to a pair of particles arriving at the same point at time $n_1$ with different velocities and creating a shock.
Let us also note that one can recover $\Psi^{n_0,n_1}_\omega w$ from $M$:
\[
\Psi^{n_0,n_1}_\omega w (x)= x-M(x),\quad x\in\R.
\]

This description of solutions of the Burgers equation in terms of monotone maps 
allows one to study the Burgers dynamics as evolution in spaces
$\GG$ or $\GG(v_-,v_+)$ consisting of cadlag functions $w\in\HH'$ (or, respectively, $w\in\HH'(v_-,v_+)$) such that 
$M_w:x\mapsto x-w(x)$ is an increasing function satisfying $\lim_{x\to\pm\infty}M(x)=\pm\infty$. Of course,
$M_w$ retains all the information about $w$ since $w(x)=x-M_w(x)$.

There are some advantages of work with monotone functions. For example,  we will utilize
the following mode of convergence: a sequence of $\GG$-functions $(w_n)_{n\in\N}$ converges to $w\in\GG$
iff  $\lim_{n\to\infty}M_{w_n}(x)=M_w(x)$ for all $x$ in $\Cc(w)$, the set of continuity points $x$ of $w$. This is equivalent to
$\lim_{n\to\infty}w_n(x)=w(x)$ for all $x\in \Cc(w)$. It is easy to define a metric on $\GG$ compatible  with this mode of convergence. 
We discuss one such metric in Section~\ref{sec:metric}.

\section{Main results}\label{sec:main_results}
 We say that $u(n,x)=u_\omega(n,x)$, $(n,x)\in\Z\times\R$ is a global solution for the
cocycle $\Psi$ if there is a set $\Omega'\in\Fc$ with $\Pp(\Omega')=1$ such that for all
$\omega\in\Omega'$, all $m$ and $n$ with $m<n$, we have $\Psi^{m,n}_\omega u_\omega(m,\cdot)= u_\omega(n,\cdot)$.
We can also introduce the global solution as a skew-invariant function: $u_\omega(x)$, $x\in\R$ is called skew-invariant
if there is a set $\Omega'\in\Fc$ with $\Pp(\Omega')=1$ such that for any $n\in\Z$, $\theta^n\Omega'=\Omega'$, and for any
$n\in\N$ and $\omega\in\Omega'$,
$\Psi^n_\omega u_\omega =u_{\theta^n\omega}$. Here and further on $\theta^n=\theta^{n,0}$ is the time shift on $\Omega$ (we recall that
space-time shifts $\theta^{n,x}$ were introduced in the previous section).

If $u_\omega(x)$ is a skew-invariant function, then
$u_\omega(n,x)=u_{\theta^n\omega}(x)$ is a global solution. One can naturally view the potentials of $u_\omega(x)$ and
$u_\omega(n,x)$ as a  skew-invariant
function and global solution for the cocycle $\hat\Phi$.

To state our first result, a description of global solutions, we need more notation.  For a subset $A$ of $\Z\times\R$, we denote by $\Fc_A$ the $\sigma$-sub-algebra
of $\Fc$  generated by $\omega|_A$, the restriction of Poisson point configuration $\omega$ to $A\times\R\times\R$, i.e., by random variables $\omega(B)$ where $B$ runs through Borel subsets of $A\times \R\times\R$. In other words, $\Fc_A$ is generated by Poissonian points with
space-time footprint in $A$.

\begin{theorem}\label{thm:global_solutions}
For every $v\in\R$
there is a unique (up to zero-measure modifications) skew-invariant function
$u_v:\Omega\to\HH'$ such that for almost every $\omega\in\Omega$, $u_{v,\omega}\in
\HH'(v,v)$.

The potential $U_{v,\omega}$ defined by $U_{v,\omega}(x)=\int^x u_{v,\omega}(y)dy$ is a unique skew-invariant potential in $\hat\HH(v,v)$. The skew-invariant functions
$U_{v,\omega}$ and $u_{v,\omega}$ are measurable  w.r.t.\ $\Fc|_{(-\N)\times\R}$, i.e., they depend only on the
history of the forcing.   The spatial random process $(u_{v,\omega}(x))_{x\in\R}$ is stationary
and ergodic with respect to space shifts.
\end{theorem}
%
%
%
%\marginpar{check if the stronger version with ``With prob.1 for all $v$'' holds}
%
%
%

%\begin{remark}\label{rem:stationary-solution}\rm 

Notice that this theorem can be interpreted as a 1F1S Principle: for
any velocity value $v$, the
solution at time $0$ with mean velocity $v$ is uniquely determined by
the history of the forcing: $u_{v,\omega}\stackrel{\rm a.s.}{=}\chi_v(\omega|_{(-\N)\times\R\times\R\times\R })$ for some
deterministic functional $\chi_v$ of the point configurations in the past, i.e., in $(-\N)\times\R\times\R\times\R$ (we actually describe
$\chi_v$ in the proof, it is constructed via one-sided action minimizers). Since the forcing is stationary in time, we obtain that $u_{v,\theta^n\omega}$ is a stationary
process in $n$, and the distribution of $u_{v,\omega}$ is an invariant distribution for the corresponding Markov
semi-group, concentrated on $\HH'(v,v)$.
%\end{remark}
\medskip

The next result shows that each of the global solutions constructed in Theorem~\ref{thm:global_solutions} plays the
role of a one-point pullback attractor.
 To describe the domains of attraction  we need to introduce several assumptions on initial potentials
$W\in\HH$. Namely, we will assume that there is $v\in\R$ such that $W$ and $v$ satisfy one of the following sets of
conditions:
\begin{align}
v&=0,\notag\\
\liminf_{x\to+\infty} \frac{W(x)}{x}&\ge 0,  \label{eq:no_flux_from_infinity}\\
\limsup_{x\to-\infty} \frac{W(x)}{x}&\le 0,\notag
\end{align}
or
\begin{align}
v&> 0,\notag \\
\lim_{x\to-\infty} \frac{W(x)}{x}&= v,\label{eq:flux_from_the_left_wins}\\
\liminf_{x\to+\infty} \frac{W(x)}{x}&> -v,\notag
\end{align}
or
\begin{align}
v&< 0,\notag\\
\lim_{x\to+\infty} \frac{W(x)}{x}&= v,\label{eq:flux_from_the_right_wins}\\
\limsup_{x\to-\infty} \frac{W(x)}{x}&< -v.\notag
\end{align}

Condition~\eqref{eq:no_flux_from_infinity} means that there is no macroscopic flux of particles from infinity toward
the origin for the initial velocity profile $W'$. In particular, any $W\in\HH(0,0)$ or any $W\in\HH(v_-,v_+)$ with
$v_-\le 0$ and $v_+\ge 0$ satisfies~\eqref{eq:no_flux_from_infinity}. It is natural to call the arising phenomenon
a rarefaction fan.
We will see that in this case the long-term behavior is described
by the global solution $u_0$ with mean velocity $v=0$.

Condition~\eqref{eq:flux_from_the_left_wins} means that the initial velocity profile $W'$ creates the influx of
particles from $-\infty$ with effective velocity $v\ge 0$, and the influence of the particles at $+\infty$ is not as
strong.
In particular, any $W\in\HH(v,v_+)$ with
$v\ge 0$ and $v_+> -v$ (e.g., $v_+=v$) satisfies \eqref{eq:flux_from_the_left_wins}. We will see that in this case the
long-term
behavior is described by the global solution $u_v$.

Condition~\eqref{eq:flux_from_the_right_wins} describes a situation symmetric to~\eqref{eq:flux_from_the_left_wins},
where
in the long run the system is dominated by the flux  of particles from $+\infty$.

The following precise statement supplements Theorem~\ref{thm:global_solutions} and describes the basins of attraction of
the global solutions $u_v$ in terms of
conditions~\eqref{eq:no_flux_from_infinity}--\eqref{eq:flux_from_the_right_wins}.

\begin{theorem}\label{thm:pullback_attraction} For every $v\in\R$,
there is a set $\tilde \Omega\in\Fc$ with $\Pp(\tilde \Omega)=1$ such that if
$\omega\in\tilde \Omega$, $W\in \HH$, and one of
conditions~\eqref{eq:no_flux_from_infinity},\eqref{eq:flux_from_the_left_wins},\eqref{eq:flux_from_the_right_wins}
holds,
then $w=W'$ belongs to the domain of pullback attraction of $u_v$ in the sense of space $\GG$:
for any $n\in\R$ and any  $x\in \Cc(u_v(n,\cdot))$,
\[
\lim_{m\to-\infty} \Psi^{m,n}_\omega w(x) = u_{v,\omega}(n,x).
\]
\end{theorem}

The last statement of the theorem implies that for every $v\in\R$, the invariant
measure on $\HH'(v,v)$ described after Theorem~\ref{thm:global_solutions} is unique and for any initial condition
$w=W'\in\HH'$ satisfying one of
conditions~\eqref{eq:no_flux_from_infinity},\eqref{eq:flux_from_the_left_wins}, and~\eqref{eq:flux_from_the_right_wins},
the distribution of the random velocity profile at time $n$ weakly converges to the unique stationary
distribution on $\HH'(v,v)$ as $n\to\infty$, in the topology of the space $\GG$. However, our approach does not produce any convergence rate estimates.

\medskip

The proofs of Theorems~\ref{thm:global_solutions} and~\ref{thm:pullback_attraction} are
given
in Sections~\ref{sec:global_solutions} and~\ref{sec:attractor}, but most of the preparatory work is carried out in
Sections~\ref{sec:subadd}--\ref{sec:weak-hyperbolicity}.

The long-term behavior of the cocycles $\Phi$ and $\Psi$ defined
through the optimization
problem~\eqref{eq:optimization_problem} depends on the asymptotic behavior of the action minimizers over long time
intervals.
The natural notion that plays a crucial role in this paper is the notion of backward one-sided infinite minimizers or
geodesics. A path $\gamma:\{\ldots, n-1,n\}\to\R$ with $\gamma(n)=x$ is called a backward one-sided minimizer if its restriction onto
any time interval $\{m,\ldots,n\}$ provides the minimum to the action $A_\omega^{m,n}(W,\cdot)$ defined in~\eqref{eq:action}
among paths connecting $\gamma(m)$ to $x$.

It can be shown (see Lemma~\ref{lem:geodesic direction}) that any backward minimizer $\gamma$ has an
asymptotic slope $v=\lim_{n\to
-\infty}(\gamma(n)/n)$. On the other hand, for every space-time point $(n,x)$ and every $v\in \R$ there is a backward
minimizer with slope $v$ and endpoint $(n,x)$. The following theorem describes the most important properties of
backward minimizers associated with the random potential $F_\omega$.

\begin{theorem}
 For every $v\in\R$ there is a set of full measure $\Omega'$ such that for all $\omega\in\Omega'$ and 
all $(n,x)\in\Z\times\R$ except countably many there is a unique
backward minimizer with asymptotic slope $v$. For $\omega\in\Omega'$, any two one-sided minimizers $\gamma^1,\gamma^2$ 
with asymptotic slope $v$, satisfy
\begin{equation}
\label{eq:weak-hyperb}
\liminf_{m\to-\infty}\frac{|\gamma^1_m-\gamma^2_m|}{|m|^{-1}}=0.
\end{equation}

\end{theorem}

The proof of this core statement of this paper is spread over
Sections~\ref{sec:subadd} through~\ref{sec:weak-hyperbolicity}.
Sections~\ref{sec:subadd} through~\ref{sec:one-sided_minimizers} are technical extensions of the corresponding
results for the Poissonian forcing from~\cite{BCK:MR3110798}, although there are some new difficulties and some simplifications.
In Section~\ref{sec:subadd} we apply the sub-additive ergodic theorem to derive the linear growth of action. In
Section~\ref{sec:concentration} we prove quantitative estimates on deviations from the linear growth. We use these
results in Section~\ref{sec:straightness} to analyze deviations of optimal paths from straight lines.
In  Section~\ref{sec:one-sided_minimizers},  we prove
the existence of infinite one-sided optimal paths and their
properties.
In Section~\ref{sec:weak-hyperbolicity} we prove a weak hyperbolicity property that is actually slightly stronger 
than~\eqref{eq:weak-hyperb}. This is where new ideas are introduced and 
the exposition differs significantly from~\cite{BCK:MR3110798}, where every two minimizers are
proved to coalesce in finite time. It is this weak hyperbolicity that is used in Sections~\ref{sec:global_solutions} and~\ref{sec:attractor}
to construct global solutions and study their properties. 
In Section~\ref{sec:metric}, we discuss the convergence in space $\GG$.
Section~\ref{sec:aux} contains some auxiliary lemmas.

\section{Optimal action asymptotics and the shape function}\label{sec:subadd}

Let us start with a note that the definition of the action has a slight time-reversal asymmetry: in the potential
action, one of the endpoints of $\{n_0,\ldots,n_1\}$ is included in summation in~\eqref{eq:potential-action},
and the other is not, whereas sometimes it is useful or convenient to work with the time-reversed version of the
potential action and define it as
 \begin{equation}
\label{eq:potential-action2}
F_\omega^{n_0,n_1}(\gamma)=\sum_{k=n_0+1}^{n_1}F_\omega(k,\gamma_k).
\end{equation}
Several estimates in this and forthcoming sections hold true for both definitions of the action, and to avoid
separate discussion of both cases we will introduce an additional parameter $p\in[0,1]$ and use the following
generalized definition of the potential action:
\begin{align}
\label{eq:gen-potential-action}
F_\omega^{n_0,n_1}(\gamma)&=pF_\omega(n_0,\gamma_{n_0})+\sum_{k=n_0+1}^{n_1-1}F_\omega(k,\gamma_k)+(1-p)F_\omega(n_1,\gamma_{n_1})
\\&=\sum_{k=n_0}^{n_1-1}\bigl(pF_\omega(k,\gamma_k)+(1-p)F_\omega(k+1,\gamma_{k+1})\bigr).\notag
\end{align}
If $p=1$ or $0$, then we recover definition~\eqref{eq:potential-action} or, respectively,~\eqref{eq:potential-action2}.

\bigskip

In this section we study the asymptotic behavior of the optimal action between space-time points $(n,x)$ and $(m,y)$
denoted by
\begin{align}\label{eq:optimal_action_between_two_points}
  A^{n,m}(x,y)=A^{n,m}_\omega(x,y)&=\min_{\gamma:\gamma_n=x,\gamma_m=y}A^{n,m}_\omega(0,\gamma)
\\ &=
\min_{\gamma:\gamma_n=x,\gamma_m=y}\left(S^{n,m}(\gamma) + F_\omega^{n,m}(\gamma)\right).
\notag
\end{align}

We will use the generalized action that involves generalized potential action $F_\omega^{n_0,n_1}(\gamma)$
defined in~\eqref{eq:gen-potential-action} for an arbitrary value of parameter $p\in[0,1]$.

Although to construct stationary solutions for the Burgers equation, we
will need the asymptotic behavior as $n\to{-\infty}$, it is more convenient and equally useful to work
with positive times and studying the limit $m\to -\infty$.

We begin with some simple observations on Galilean shear transformations of the point field.
\begin{lemma}\label{lem:shear} Let $a,v\in\R$ and let $L=L^{a,v}$ be a transformation of space-time defined by
$L(n,x)=(n,x+a+vn)$.
\begin{enumerate}
 \item
Suppose that  $\gamma$ is a path defined on a time interval
$\{n_0,\ldots,n_1\}$ and let $\bar \gamma$ be defined by $(n,\bar\gamma_n)=L(n,\gamma_n)$.
Then
\[
 S^{n_0,n_1}(\bar\gamma)=S^{n_0,n_1}(\gamma)+(\gamma_{n_1} -\gamma_{n_0})v+\frac{(n_1-n_0)v^2}{2}.
\]

\item  For any $\omega\in\Omega$, let us define $L\omega$ by pointwise application of $L$ to the space-time footprint
of configuration points:   $(\tau,\eta,\xi,\kappa)\in\omega$ iff $(\tau,\eta+a+v\tau,\xi,\kappa)\in L\omega$.
Then for any $\omega\in\Omega$, any time interval $\{n_0,\ldots, n_1\}$, and  and any points $x_0,x_1,\bar x_0, \bar x_1\in\R$ satisfying
$L(n_0,x_0)=(n_0,\bar x_0)$ and $L(n_1, x_1)=(n_1,\bar x_1)$,
\[
A_{L\omega}^{n_0,n_1}(\bar x_0, \bar x_1)= A_{\omega}^{n_0,n_1}(x_0, x_1)+(x_1-x_0)v+\frac{(n_1-n_0)v^2}{2},
\]
and $L$ maps minimizers realizing $A_{\omega}^{n_0,n_1}(x_0, x_1)$ onto minimizers
realizing $A_{L\omega}^{n_0,n_1}(\bar x_0, \bar x_1)$.
\item 
The measure~$\Pp$ is invariant under~$L$.

\item For any points $x_0,x_1,\bar x_0,\bar x_1$ and any time interval $\{n_0,\ldots, n_1\}$,
\[
  A^{n_0,n_1}(\bar x_0, \bar x_1)\stackrel{distr}{=} A^{n_0,n_1}(x_0,
x_1)+(x_1-x_0)v+\frac{(n_1-n_0)v^2}{2},
\]
where
\[
 v=\frac{(\bar x_1-x_1)-(\bar x_0-x_0)}{n_1-n_0}.
\]

\end{enumerate}
\end{lemma}
\bpf The first part of the Lemma is a simple computation:
\begin{align*}
 S^{n_0,n_1}(\bar\gamma)& 
 =
 \frac{1}{2}\sum_{n=n_0}^{n_1-1}(\gamma_{n+1}+v(n+1)-\gamma_n-vn)^2
 =
 \frac{1}{2}\sum_{n=n_0}^{n_1-1}(\gamma_{n+1}-\gamma_n+v)^2\\
 &=\frac{1}{2}\sum_{n=n_0}^{n_1-1}(\gamma_{n+1}-\gamma_n)^2+ v\sum_{n=n_0}^{n_1-1}(\gamma_{n+1}-\gamma_n)+
 \frac{1}{2}\sum_{n=n_0}^{n_1-1}v^2.
\end{align*}
The second part follows from the first one. 
To prove the third part, it is sufficient to notice that the measure $\mu$ driving the Poisson process is invariant under $L$.
The last part is a consequence of the previous two parts, once one finds an
appropriate Galilean transformation sending $(n_0,x_0)$ to $(n_0, \bar x_0)$ and $(n_1,x_1)$ to $(n_1, \bar x_1)$.
\epf

\medskip

The next useful property is the sub-additivity of action along any direction: for any velocity $v\in
\R$, and any $n,m\geq 0$, we have
\[ A^{0,n+m}(0,v(n+m))\leq A^{0,n}(0,vn) + A^{n,n+m}(vn,v(n+m)).\]
This means that we can apply Kingman's sub-additive ergodic theorem to the function $n\mapsto A^{0,n}(0,vn)$ if we can
show that $-\E A^{0,n}(0,vn)$ grows at most linearly in $n$. We claim this linear bound in
the following result:
\begin{lemma}\label{lem:linbound}
Let $v\in \R$. There exist constants $C=C(v)>0$ and $n_0>0$ such that for all $n\geq n_0$
\[ \E |A^{0,n}(0,vn)| \leq Cn.\]
\end{lemma}
\bpf This statement and its proof are adapted from~\cite{BCK:MR3110798}.
Lemma~\ref{lem:shear} implies that it is enough to prove this for $v=0$. So in this proof we work with
$A^n=A^n(0,0)$.

Let $(\gamma_0,\ldots\gamma_n)$ be a  path realizing $A^n$.  Let us denote 
\begin{equation}
\Sigma(\gamma)=\sum_{j=0}^{n-1}k_j, 
\label{eq:Sigma}
\end{equation}
where $k_j=|[\gamma_{j+1}]-[\gamma_j]|+1$, and set
%One can use Lemma~\ref{lem:properties-of-evolution-operator}  to derive that for every two points $(n_0,x_0)$
%and $(n_1,x_1)$ with $n_0<n_1$, there are minimizers $\gamma^*$ and $\gamma_*$ realizing $A^{n_0,n_1}(x_0,x_1)$ such that
%all other minimizers $\gamma$ realizing $A^{n_0,n_1}(x_0,x_1)$ satisfy $\gamma_*\le \gamma_n\le \gamma^*_n$. We call $\gamma^*$
%and $\gamma_*$ the rightmost minimizer and the leftmost minimizer. This notion makes the following definition possible:
\[ E_{n,m} =\{ \Sigma(\gamma) = m\ \textrm{for some}\ \gamma\ \textrm{realizing}\ A^n\}.\]

\begin{lemma}\label{lem:distribution-of-animal-size} There
are constants $C_1>0$ and $R\in\N$ such that if $m\geq R n$, then
\[ \Pp(E_{n,m}) \leq \exp(-C_1 m^2/n).\]
\end{lemma}

\bpf If a path $\gamma$ realizes $E_{n,m}$, then
\begin{equation}
A^n\ge \frac12\sum_{j=0}^{n-1}\left(k_j-2\right)_+^2 - F_\omega^*([\gamma_0],\ldots,[\gamma_n]),
\label{eq:lower-estimate-on-An}
\end{equation}
where $a_+=0\vee a$, and 
\[
F_\omega^*(i_0,\ldots,i_n)=\sum_{j=0}^{n}F_\omega^*(j,i_j),\quad i_0,\ldots,i_n\in\Z.
\]
We recall that $F_\omega^*(j,k)$ was introduced in~\eqref{eq:F-star}.
Since $A^n$ is optimal,
\begin{equation}
A^n\le A_\omega^{0,n}(0,0,\ldots,0)=pF_\omega(0,0)+\sum_{j=1}^{n-1}F_\omega(j,0)+ (1-p)F_\omega(n,0).
\label{eq:upper-estimate-on-An}
\end{equation}
%
%Since for any $\lambda>0$ we can choose $r$ so that $-\lambda r+\phi(\lambda)<0$, we conclude that 
Since $a\mapsto (a-2)_+^2$ is convex, we can use Jensen's inequality to see that
\[ \frac12\sum_{j=0}^{n-1}\left(k_j-2\right)_+^2 \ge \frac12 n 
\left(\frac{m}{n}-2\right)_+^2.\]
Combining this with~\eqref{eq:lower-estimate-on-An} and~\eqref{eq:upper-estimate-on-An}, we obtain
\begin{equation*}
 pF_\omega(0,0)+\sum_{j=1}^{n-1}F_\omega(j,0)+ (1-p)F_\omega(n,0)
 \ge \frac12 n \left(\frac{m}{n}-2\right)_+^2 - F_\omega^*([\gamma_0],\ldots,[\gamma_n]).
%\label{eq:comparing_zero_path_to_best}
\end{equation*}
We conclude that on $E_{n,m}$,
\begin{equation}
2 F^*_{\omega, n,m} 
 \ge \frac12 n \left(\frac{m}{n} -2\right)_+^2,
\label{eq:upper_bound_on_animal_weight}
\end{equation}
where
\[
F^*_{\omega, n,m}=\max \Biggl \{F_\omega^*(i_0,\ldots,i_n):\ \sum_{j=0}^{n-1} (|i_{j+1}-i_j|+1)\le m\Biggr\}.
\]
So we need a tail estimate for $F^*_{\omega,n,m}$. Since the distribution of $F_\omega^*(i_0,\ldots,i_n)$ does not depend
on the choice of $(i_0,\ldots,i_n)$, we obtain that for any $r>0$,
\begin{equation}
\Pp\{F_{\omega,n,m}^*>r\}\le N_{n,m}\Pp\left\{\sum_{j=0}^{n}F_\omega^*(0,j)>r\right\},
\label{eq:animal-estimate-1}
\end{equation}
where $N_{n,m}$ is the size of the set $\{(i_0,\ldots,i_n): \sum_{j=0}^{n-1} (|i_{j+1}-i_j|+1)\le m\}$, $n\le m$. Let us estimate
$N_{n,m}$ first.
The number of ways to represent $m$ as a sum of $n$ ordered nonnegative terms is ${\binom{m+n-1}{n-1}}\le 2^{m+n-1}$.
Since we also may choose the sign of $i_{j+1}-i_j$, we obtain an additional factor of $2^n$, so we obtain a crude estimate
\begin{equation}
N_{n,m}\le e^{\rho m},\quad m\ge n,
\label{eq:animal-estimate-2}
\end{equation}
for some $\rho>0$.
We also have,
\[
\Pp\left\{\sum_{j=0}^{n}F_\omega^*(0,j)>r\right\}\le e^{-\lambda r}\E e^{\lambda \sum_{j=0}^{n}F_\omega^*(0,j)}=e^{-\lambda r+\varphi(\lambda)(n+1)},\quad r,\lambda>0.
\]
and, combining this with~\eqref{eq:animal-estimate-1} and~\eqref{eq:animal-estimate-2}, we obtain
\[
\Pp\{F^*_{\omega,n,m}>my\}\le e^{\rho m-\lambda my +\varphi(\lambda)(n+1)}\le e^{m(\rho-\lambda y +2\varphi(\lambda))},\quad m\ge n.
\]
So choosing first any  $\lambda>0$, then any $y_0>0$ such that 
\[\rho-\lambda y_0 +2\varphi(\lambda) < -\lambda y_0/2,\] we obtain
\begin{equation}
\Pp\{F^*_{\omega,n,m}>my\}\le e^{-Kmy},\quad y\ge y_0,\ m\ge n,
\label{eq:large-dev-of-omega-star}
\end{equation}
where $K=\lambda/2$.
If $R>2$ and $m\ge Rn$, then using~\eqref{eq:upper_bound_on_animal_weight}, denoting 
\[
y=\frac{1}{4}\frac{n}{m}\left(\frac{m}{n}-2\right)^2=\frac{1}{4}\frac{(m-2n)^2}{nm}\ge \frac{m(1-2R^{-1})^2}{4n},
\] 
noticing that the right-hand side is bounded below by 
$R(1-2R^{-1})^2/4$ which exceeds $y_0$ for sufficiently large $R$, we obtain 
from~\eqref{eq:large-dev-of-omega-star}:
\[
\Pp(E_{m,n})\le e^{-K m y}\le e^{-Km\frac{m(1-2R^{-1})^2}{4n}},
\]
and the lemma follows with $C_1=K(1-2R^{-1})^2/4$. \epf

\begin{lemma}\label{eq:moments-of-omega-star} For any $k\ge1$, there is $c_k>0$ such that for all $n,m\in\N$ with $m\ge n$,
\[
\E F_{\omega,n,m}^{*k} \leq c_k m^k.
\]
\end{lemma}
\bpf
Clearly,
\[ \E F_{\omega,n,m}^{*k} = \E F_{\omega,n,m}^{*k}\ONE_{\{F_{\omega,n,m}^*\le y_0m\}}+\E F_{\omega,n,m}^{*k}\ONE_{\{F_{\omega,n,m}^*> y_0m\}}.
\]
We can bound the first term simply by $(y_0 m)^k$.
For the second term we can use \eqref{eq:large-dev-of-omega-star}:
\begin{align*}
\E F_{\omega,n,m}^{*k}\ONE_{\{F_{\omega,n,m}^*> y_0m\}}
&\le \sum_{i=0}^\infty ((y_0+i+1)m)^k\Pp\{F_{\omega,n,m}^*\in((y_0+i)m, (y_0+i+1)m]\}\\
&\le \sum_{i=0}^\infty (y_0+i+1)^km^k\Pp\{F_{\omega,n,m}^*>(y_0+i)m\}\\
&\le m^k  \sum_{i=0}^\infty (y_0+i+1)^ke^{-Km(y_0+i)}.
\end{align*}
The series on the right-hand side is uniformly convergent in $m$, so the proof is completed. 
\epf

We can now prove Lemma~\ref{lem:linbound}.
From~\eqref{eq:lower-estimate-on-An} and~\eqref{eq:upper-estimate-on-An} we know
that $|A^n|\le F_{\omega,n,m}^*$ on $E_{n,m}$. So, using Lemmas~\ref{lem:distribution-of-animal-size}
and~\ref{eq:moments-of-omega-star}, and the fact that $F_{\omega,n,m}^*$ is nondecreasing in $m$, we obtain
\begin{align*}
\E |A^n| & =  \E  \sum_{n\le m\le Rn} |A^n|\ONE_{E_{n,m}} + \sum_{m> Rn} \E |A^n|\ONE_{E_{n,m}}\\
& \le  \E\sum_{n\le m\le Rn}  F_{\omega,n,m}^*\ONE_{E_{n,m}} +  \sum_{m> Rn } \E F_{\omega,n,m}^*\ONE_{E_{n,m}}\\
& \le  \E F_{\omega,n,Rn}^* + \sum_{m>Rn} \sqrt{\E F_{\omega,n,m}^{*2}}\sqrt{\Pp(E_{n,m})}\\
& \le  c_1Rn + \sqrt{c_2}\sum_{m>Rn}   m\exp(-C_1m^2/(2n))\\
& \le  Cn,
\end{align*}
for $C$ big enough. \epf

In fact, we can use the last calculation to obtain the following generalization of Lemma~\ref{lem:linbound}  for
higher moments of $A^n$:

\begin{lemma} \label{lem:moments_of_At} Let $k\in\N$. Then there is a  constant $C(k)>0$ such that
\[ \E |A^n|^k\le C(k) n^k,\quad n\in\N.\]
\end{lemma}

\begin{remark}\label{rem:lemmas-hold-for-restricted-optimization} \rm In fact, 
the analysis of all the proofs in this section shows that all the results above are valid if
$A^n$ is replaced by 
\[
\tilde A^n=\min_{\gamma\in \Gamma_n}\left(S^{0,n}(\gamma) + F_\omega^{0,n}(\gamma)\right), 
\]
where $\Gamma_n$ is any set satisfying two conditions: (i) all elements of $\Gamma_n$ are paths $\gamma:\{0,\ldots,n\}\to\R$ 
such that $\gamma_0=0$, $\gamma_n=0$; (ii) $\Gamma_n$ contains the path $(0,0,\ldots,0)$. 
We will need this in the proof of Lemma~\ref{lem:difference_of_rectangle_expectation_from_true}.
\end{remark}

\medskip

Now we arrive to the main result of this section describing the {\it shape function} $\alpha$ for our model. 
\begin{lemma}\label{lem:shape-function} 
For each $v\in\R$, the number $\alpha(v)\in\R$ defined by
\begin{equation}
\label{eq:definition-of-alpha-via-subadditive-ergodic}
\alpha(v)=\inf_n \frac{\E A^{0,n}(0,vn)}{n},
\end{equation}
satisfies
\begin{equation}
\frac{A^{0,n}(0,vn)}{n}\to\alpha(v),\quad \text{\rm a.s.\  and
in\ } L^1,\quad n\to\infty,
\label{eq:def_of_shape_function}
\end{equation}
and does not depend on the choice of constant $p$ in definition~\eqref{eq:gen-potential-action}.
Moreover,
\begin{equation}
 \alpha(v)=\alpha(0)+\frac{v^2}{2},\quad v\in\R.
\label{eq:shape-function-is-quadratic}
\end{equation}
\end{lemma}
\bpf The number $\alpha(v)$ is finite due to Lemma~\ref{lem:linbound}.
 The sub-additive ergodic theorem now implies~\eqref{eq:def_of_shape_function}. The independence on~$p$ follows from
 \[
 \lim_{n\to\infty}\frac{F_\omega(n,vn)}{n}= \lim_{n\to\infty}\frac{F_\omega(0,0)}{n}=0,\quad \text{\rm a.s.\  and
in\ } L^1.
 \]

To prove~\eqref{eq:shape-function-is-quadratic} we notice that the
Galilean shear map $(n,x)\mapsto(n,x+vn)$ transforms the paths connecting $(0,0)$ to
$(n,0)$ into paths connecting $(0,0)$ to $(n,vn)$.
 Lemma~\ref{lem:shear} implies that 
under this map
the optimal action over these paths is altered by a deterministic
correction $v^2n/2$, but the measure $\Pp$ is invariant under the lift of this transformation onto $\Omega$, so the lemma follows. 
\epf

We know now from~\eqref{eq:def_of_shape_function} that  $A^{0,n}(0,vn)\sim\alpha(v) n$ as $n\to\infty$ with
probability~1. However, this is not enough for our purposes since we need quantitative estimates on deviations of
$A^{0,n}(0,vn)$ from~$\alpha(v) n$. This is the material of the next section.

\section{Concentration inequality for optimal action}\label{sec:concentration}

The goal of this section is to prove a concentration inequality for
$A^n(vn)=A^n(0,vn)=A^n_\omega(0,vn)=A^{0,n}_\omega(0,vn)$. The methods
are similar to those of~\cite{BCK:MR3110798}, except for some technical moments. In particular, we use the
Azuma--Hoeffding  inequality instead of the Kesten inequality. Throughout this section we work with
the version of action defined in~\eqref{eq:gen-potential-action}.

\begin{theorem}\label{thm:concentration_around_alphat} There are positive
constants $c_0,c_1,c_2,c_3,c_4$ such that for any $v\in\R$, all
$n>c_0$, and
all
$u\in(c_3n^{1/2}\ln^2 n, c_4n^{3/2}\ln n]$,
\[
 \Pp\{|A^n(0,vn)-\alpha(v) n |>u\}\le c_1\exp\left\{-c_2\frac{u}{n^{1/2}\ln n}\right\}.
\]
\end{theorem}
\begin{remark}\rm In our setting one can also prove a similar bound for small values of $n$, but we will
mostly need the theorem as it is stated.
\end{remark}

Due to the invariance under shear transformations (Lemmas~\ref{lem:shear} and~\ref{lem:shape-function}), it is
sufficient to prove this theorem for $v=0$.
We will first derive a similar inequality with $\alpha(0)n$ replaced by $\E A^n$, and then we will have to
estimate the corresponding approximation error.

\begin{lemma}\label{lem:concentration_around_mean} There are positive constants $b_0,b_1,b_2,b_3,b_4$ such that for all
$n>b_0$ and all
$u\in(b_3n^{1/2}\ln n, b_4n^{3/2}\ln n]$,
\[
 \Pp\{|A^n-\E A^n|>u\}\le b_1\exp\left\{-b_2\frac{u}{n^{1/2}\ln n}\right\}.
\]
\end{lemma}

We will need an estimate on probabilities of the following events:
\[
B_n(u)=\left\{ \max_{0\le k \le n} |\gamma_k|>u\ \text{\rm for some}\ \gamma\ \text{\rm realizing}\ A^n\right\}.
\]
\begin{lemma}\label{lem:path_in_wide_rectangle_whp} There is a constant $C_2>0$ such that if $n\in\N$ and $u\ge Rn$, then
\[
\Pp(B_n(u))\le C_2 \exp(-C_1 u^2/n),
\]
where constants $C_1$ and $R$ were introduced in Lemma~\ref{lem:distribution-of-animal-size}.
\end{lemma}
\bpf
If $\max_{0\le k \le n}|\gamma_k|>u$, then $\Sigma(\gamma)\ge u$. Lemma~\ref{lem:distribution-of-animal-size}
implies
\begin{align*}
\Pp(B_n(u))& \le \sum_{m\ge u} \exp(-C_1m^2/n) \\ &\le \exp(-C_1u^2/n) \sum_{m\ge u} \exp(-C_1(m^2-u^2)/n)
\\ & \le \exp(-C_1u^2/n) \sum_{m\ge u} \exp(-C_1 R (m-u))
 \\&\le C_2\exp(-C_1 u^2/n),
\end{align*}
where $C_2=\sum_{m=0}^\infty e^{-C_1Rm}$.
\epf

\medskip

Having Lemma~\ref{lem:path_in_wide_rectangle_whp} in mind, we define $\tilde A^n$ to be the optimal action over all
paths connecting $(0,0)$ to $(0,n)$ and staying within $[-Rn,Rn]$.

\begin{lemma}\label{lem:probability_that_Kesten_action_worse} Let constants $R,C_1,C_2$ be defined in
Lemmas~\ref{lem:distribution-of-animal-size} and~\ref{lem:path_in_wide_rectangle_whp}. For any $n\in\N$,
\[
 \Pp\{A^n\ne \tilde A^n\}\le  C_2\exp(-R^2C_1 n).
\]
\end{lemma}
\bpf It is sufficient to notice that
\[
 \Pp\{A^n\ne \tilde A^n\}\le\Pp(B_n(Rn))
\]
and apply Lemma~\ref{lem:path_in_wide_rectangle_whp}.
\epf

\begin{lemma}\label{lem:difference_of_rectangle_expectation_from_true} There is a constant $D_1$ such that for all
$n\in\N$,
\[
 0\le \E \tilde A^n-\E A^n\le\E  (| A^n|+|\tilde A^n|)\ONE_{B_n(Rn)} \le D_1.
\]
\end{lemma}
\bpf
The first inequality is obvious, since $A^n\le \tilde A^n$.
We also have
\begin{align}
\E \tilde A^n-\E A^n \notag
&\le \E  (\tilde A^n-A^n)\ONE_{B_n(Rn)}
\\ &\le \E  | A^n|\ONE_{B_n(Rn)}+\E  |\tilde A^n|\ONE_{B_n(Rn)}.
\label{eq:auxiliary-inequality-in-lemma-difference_of_rectangle_expectation_from_true}
\end{align}
%
%HHHHHHHHHHHHHHHHHHHH
%
%To estimate the first term we apply Lemmas~\ref{lem:distribution-of-animal-size} and~\ref{eq:moments-of-omega-star}:
%\begin{align*}
% \E  | A^n|\ONE_{\left\{\max_{0\le k \le n}|\gamma_k|>Rn\right\}}
%&\le \sum_{m>Rn} \E( F_{\omega,n,m}^*\ONE_{E_{n,m}})
%\\&\le \sum_{m>Rn} \sqrt{\E F_{\omega,n,m}^{*2} }\sqrt{\Pp(E_{n,m})}
%\\&\le  \sum_{m>Rn} \sqrt{c_2}n\exp(-C_1m^2/(2n)),
%\end{align*}
%and the series in the right-hand side is 
%convergent uniformly in $n$.
%It is also easy to see that Lemma~\ref{eq:moments-of-omega-star} holds true if $A^n$ is replaced by $\tilde A^n$,
%so we can apply the same reasoning to the second term of~\eqref{eq:auxiliary-inequality-in-lemma-difference_of_rectangle_expectation_from_true}.
%HHHHHHHHHHHHHHHHHHHH
%
Lemmas~\ref{lem:moments_of_At} and~\ref{lem:path_in_wide_rectangle_whp} give:
\begin{align*}
\E | A^n|\ONE_{B_n(Rn)}&\le \sqrt{\E |A^n|^2}\sqrt{\Pp(B_n(Rn))}
\\ &\le \sqrt{C(2)} n \sqrt{C_2\exp(-C_1 R^2n)}.
\end{align*}
Remark~\ref{rem:lemmas-hold-for-restricted-optimization} shows that the same estimate applies to $\tilde A^n$:
\[
\E | \tilde A^n|\ONE_{B_n(Rn)}\le \sqrt{C(2)} n \sqrt{C_2\exp(-C_1 R^2n)},
\]
so the the right-hand side of~\eqref{eq:auxiliary-inequality-in-lemma-difference_of_rectangle_expectation_from_true} is uniformly bounded in $n$.
\epf

As in~\cite{BCK:MR3110798}, we could use Kesten's concentration inequality to  
estimate the deviations of $\tilde A^n$ from its mean. However, in our discrete time setting, we also
can use a more basic tool, the following Azuma--Hoeffding inequality \cite{MR0144363:Hoeffding}:

\begin{lemma}\label{lem:Azuma} Let $(\Fc_k)_{0\le k\le N}$ be a filtration. Suppose $(M_k)_{0\le k\le N}$ is a martingale with respect to $(\Fc_k)_{0\le k\le N}$ such that for some constant $c>0$ the increments $\Delta_k=M_k-M_{k-1}$ satisfy
\[
 |\Delta_k|<c,\quad k=1,\ldots,N.
\] 
Then
\[
 \Pp\{M_N-M_0\ge x\}\le \exp\left(-\frac{x^2}{2Nc^2}\right).
\]
\end{lemma}

To use this lemma in our framework, we must introduce an appropriate martingale.
For a given natural $n\ge 2$, we will use $N=n-1$. To define a filtration $(\Fc_k)_{0\le k\le n-1}$,
we introduce $Q_j = \{j\}\times [-Rn,Rn)$ and $Q^+_j=\{j\}\times [-Rn-1,Rn+1)$, $j=1,\ldots, n-1$.

Using the notation  
introduced before the statement of Theorem~\ref{thm:global_solutions}, we set
$\Fc_0=\{\emptyset,\Omega\}$ and
\[
 \Fc_k=\sigma\left(\omega\bigr|_{\bigcup_{j=1}^k Q^+_{j}} \right),\quad k=1,\ldots,n-1.
\]
We  introduce a martingale $(M_k,\Fc_k)_{0\le k\le n-1}$ by
\[
 M_k=\E (\tilde A^n|\Fc_k),\quad 0\le k\le n-1.
\]
Note that $M_0=\E \tilde A^n$ and $M_{n-1}= \tilde A^n$.

Let us denote by $P_k$ the distribution of $\omega\bigr|_{Q^+_k}$ on the sample space $\Omega_k$ of finite
point configurations in $Q^+_k\times \R\times\R$. For $\omega,\sigma\in\prod_{k=1}^{n-1}\Omega_k$ we write
\[
 [\omega,\sigma]_k=(\omega_1,\ldots,\omega_k,\sigma_{k+1},\ldots,\sigma_{N})\in\prod_{k=1}^{n-1}\Omega_k.
\]
Then, for $k=1,\ldots,n-1$,
\begin{align*}
 \Delta_k(\omega_1,\ldots,\omega_k)&:=M_k-M_{k-1}\\
 &= \int \tilde A^n_{[\omega,\sigma]_k}\prod_{j=k+1}^{n-1} dP_j(\sigma_j)-\int
\tilde A^n_{[\omega,\sigma]_{k-1}}\prod_{j=k}^{n-1}
dP_j(\sigma_j)\\
&=\int \left(\tilde A^n_{[\omega,\sigma]_k}- \tilde A^n_{[\omega,\sigma]_{k-1}}\right)\prod_{j=k}^{n-1}
dP_j(\sigma_j).
\end{align*}

For any set $B\in\Z\times\R$, we will denote 
\begin{equation}
\label{eq:potential-maximum}
F^*_\omega(B)=\sup\{|F_\omega(n,x)|: (n,x)\in B\}.  
\end{equation}

\begin{lemma}
\label{lm:estimating_integrand_in_Kesten_martingale_difference}
Let  $k\in\{1,\ldots,n-1\}$. Then
\[
|\tilde A^n_{[\omega,\sigma]_k}- \tilde A^n_{[\omega,\sigma]_{k-1}}| \le
F^*_{[\omega,\sigma]_k}(Q_{k})+F^*_{[\omega,\sigma]_{k-1}}( Q_{k}).
\]
\end{lemma}

\bpf Changing $\omega_k$ to $\sigma_k$ 
we change the action of any path passing through $Q_k$ by at most
$F^*_{[\omega,\sigma]_k}(Q_k)+F^*_{[\omega,\sigma]_{k-1}}(Q_k)$, and the statement follows.
\epf

\bigskip

The next step is to introduce a truncation of configuration $\omega$.
For $j\in\{1,\ldots,n-1\}$ and $i\in\{-Rn-1,\ldots, Rn\}$, we denote $Q_{ji}=\{j\}\times [i,i+1)$.
We define $\bar\omega$ by erasing all configuration points of
$\omega$ in each
block $Q_{ji}\times\R\times\R$ with $\omega(Q_{ji}\times\R\times\R)> b\ln n$. The value $b>0$ will be chosen later.
 The restrictions of $\bar \omega$ to blocks $Q_{ji}\times\R\times\R$ are
mutually independent.
Lemma~\ref{lm:estimating_integrand_in_Kesten_martingale_difference} applies to truncated configurations as well.
Since for any segment $\{k\}\times[a,b]$ and any point configuration~$\omega$,
\[
 F^*_{\omega}(\{k\}\times[a,b])\le \omega(\{k\}\times[a-1,b+1]\times\R\times\R),  
\]
Lemma~\ref{lm:estimating_integrand_in_Kesten_martingale_difference} gives:
\[
| \tilde A^n_{[\bar\omega,\bar\sigma]_k}- \tilde A^n_{[\bar\omega,\bar\sigma]_{k-1}}| \le 6 b \ln n,\quad k=1,\ldots,n-1,
\]
where $\bar\sigma$ is the truncation of $\sigma$.
Therefore,
\[
 |\Delta_k(\bar\omega_1,\ldots,\bar\omega_k)|\le
 \int   6 b\ln n \prod_{j=k}^N
dP_j(\sigma_j)\le 6 b\ln n.
\]

Now Lemma~\ref{lem:Azuma} directly implies the following estimate:
\begin{lemma}\label{lem:outcome_of_Kestens_lemma}For all $x>0$,
\[
 \Pp\left\{|\tilde A^n(\bar\omega)-\E \tilde A_n(\bar\omega)|>x\right\}\le 2\exp\left(-\frac{x^2}{72(n-1)b^2\ln^2 n}\right).
\]
\end{lemma}

Let us now estimate the discrepancy between $\tilde A^n(\bar\omega)$ and $\tilde A^n(\omega)$:

\begin{lemma} 
\label{lem:deviation_of_truncation_from_original}
If $b>2$, then there is $n_0$ such that  for any $x>0$ and any $n>n_0$,
\[
\Pp\{|\tilde A^n(\bar \omega) -\tilde A^n(\omega)|>x\}\le 2e^{-x}.
\]
\end{lemma}
\bpf Let us define $\xi_{ji}=\omega(Q_{ji}\times\R\times\R)$. Then
\[
|\tilde A^n(\bar \omega) -\tilde A^n(\omega)|\le \sum_{j,i} \xi_{ji}\ONE_{\{\xi_{ji}>b\ln n\}}.
\]
By Markov's inequality and mutual independence of $\xi_{ji}$,
\[
 \Pp\left\{\sum_{j,i} \xi_{ji}\ONE_{\{\xi_{ji}>b \ln n\}}>x\right\}\le
e^{-x}\left[
\E e^{\xi\ONE_{\{\xi>b\ln n\}}}
\right]^{2R(n+1)(n-1)},
\]
where $\xi$ is a r.v.\ with the same distribution as any of $\xi_{ji}$.
The lemma will follow from
\begin{equation*}%\label{eq:lim_factor_in_Markov_equals_1}
\lim_{n\to\infty}\left[\E e^{\xi\ONE_{\{\xi>b \ln t\}}}\right]^{2R(n+1)(n-1)}=1,
\end{equation*}
which is implied by
\[
 \E e^{\xi\ONE_{\{\xi>b \ln n\}}}\le 1+\frac{\E  e^{2\xi}}{e^{b\ln n}}\le 1+\frac{\E
e^{2\xi}}{n^b},
\]
and $b>2$.\epf

We also need an estimate on the discrepancy between $\E \tilde A^n(\bar\omega)$ and $\E \tilde A^n(\omega)$.
It is a direct consequence of Lemma~\ref{lem:deviation_of_truncation_from_original}:
\begin{lemma}\label{lem:expectation_of_truncation-expectation_over_Kesten_rectangle}
There is a constant $D_2$ such that for all
$n\in\N$,
\[
|\E \tilde A^n(\bar\omega)-\E \tilde A^n(\omega)| <D_2.
\]
\end{lemma}

\smallskip

\bpf[Proof of Lemma~\ref{lem:concentration_around_mean}]
Lemmas~\ref{lem:difference_of_rectangle_expectation_from_true}
and~\ref{lem:expectation_of_truncation-expectation_over_Kesten_rectangle}
 imply that
 for $u>D_1+D_2$
\begin{align*}
 \Pp\{|A^n(\omega)-\E A^n(\omega)|>u\}\le& \Pp\{|A^n(\omega)-\tilde A^n(\omega)|>(u-D_1-D_2)/3 \}
\\&+\Pp\{|\tilde
A^n(\omega)-\tilde A^n(\bar\omega)|>(u-D_1-D_2)/3\}
\\&+\Pp\{|\tilde
A^n(\bar\omega)-\E\tilde A^n(\bar\omega)|>(u-D_1-D_2)/3\}.
\end{align*}
The lemma follows from the estimates of the three terms provided by
Lemmas~\ref{lem:probability_that_Kesten_action_worse},~\ref{lem:outcome_of_Kestens_lemma},
and
\ref{lem:deviation_of_truncation_from_original}.\epf

The following lemma gives an estimate on how $\E A^n$ changes under argument doubling. We
will use this lemma  estimate on $\E A^n-\alpha(0)n$ and bridge the gap between
Lemma~\ref{lem:concentration_around_mean} and Theorem~\ref{thm:concentration_around_alphat}.

\begin{lemma}\label{lem:doubling}
There is a number $b_0>0$ such that for any $n>n_0$,
\[
 0\le 2\E A^n - \E A^{2n} \le b_0n^{1/2}\ln^2n.
\]
\end{lemma}
\bpf The first inequality follows from $A^{0,2n}(0,0)\le A^{0,n}(0,0)+A^{n,2n}(0,0)$. Let us prove the second one.

Let $\gamma$ be the (rightmost, for definiteness) minimizer from $(0,0)$ to $(2n,0)$. Then
\[
 A^{2n}\ge \min_{|x|\le 2Rn} A^{0,n}(0,x)+\min_{|x|\le
2Rn}A^{n,2n}(x,0)+A^{2n}\ONE_{\left\{\max_{0\le k \le 2n}|\gamma_k| >2Rn\right\}}.
\]
Therefore, by  Lemma~\ref{lem:difference_of_rectangle_expectation_from_true},
\begin{equation}
\E A^{2n}\ge  \E \min_{|x|\le 2Rn} A^{n}(0,x)+\E\min_{|x|\le
2Rn}A^{n,2n}(x,0) - D_1.
\label{eq:doubling+constant}
\end{equation}
We will estimate the first term of the right-hand side. To that end we define 
$I_n=\{-2Rn,\ldots,2Rn-2,2Rn-1\}$ and $\bar I_n = I_n\cup \{2R n\}$. 
Let now $\gamma$ be the minimizer from $(0,0)$ to $(x,n)$, with $x\in[k,k+1]$ for some
$k\in I_n$. 
Let us introduce $\gamma^+$ and $\gamma^-$ satisfying $\gamma^+_{j}=\gamma^-_{j}=\gamma_j$ for all $j<n$ and $\gamma^+_{n}=k+1$,
$\gamma^-_{n}=k$. Comparing $\gamma$ to $\gamma_+$ if $\gamma_{n-1}\ge k+{1/2}$ and to $\gamma_-$ if $\gamma_{n-1}< k+{1/2}$,
we see that 
\[
S^{0,n}(\gamma_\pm)\le S^{0,n-1}(\gamma) +\frac{(1/2)^2}{2} =S^{0,n-1}(\gamma) +\frac{1}{8},
\]
so
\[
A^n(x)\ge \min \{A^n (k), A^n(k+1)\} -  2F^*(n,k)-\frac{1}{8},
\] 
But 
\[
\E \max_{k\in I_n}  F^*_\omega(n,k)\le 3\E \max  \bigl\{\omega(\{n\}\times[k,k+1]): {k=-2Rn-1,\ldots,k=2Rn+1}\bigr\},
\]
and so there is a constant $c>0$ such that
\[
 \E \min_{|x|\le 2Rn} A^{n}(0,x) \ge  \E \min_{k\in \bar I_n} A^{n}(0,k) - c(\ln n +1).
\]
Lemma~\ref{lem:shear} implies
$\min_{x} \E A^n(x)=\E A^n(0)$. Therefore,
denoting
\[
 X_n=\max_{k\in \bar I_n} \{(\E A^n(k)-A^n(k))_+\},
\]
we obtain
\begin{align}
 \E \min_{|x|\le 2Rn} A^{n}(0,x) &\ge \min_{k\in \bar I_n}\E A^n(k) - \E X_n - c(\ln n+1) \notag
\\&\ge \E A^n - \E X_n - c(\ln n+1), \label{eq:doubling1}
\end{align}
Similarly, we obtain for the second term in~\eqref{eq:doubling+constant}
\begin{align}
\E\min_{|x|\le 2Rn}A^{n,2n}(x,0)&\ge \E A^{n,2n}(0,0) -\E Y_n -c(\ln n +1),\label{eq:doubling2}
\end{align}
where
\[
Y_n=\max_{k\in \bar I_n} \{(\E A^{n,2n}(k,0)-A^{n,2n}(k,0))_+\}.
\]
For a constant $r$ to be determined later, we introduce the event
\[
 E=\{X_n+Y_n\le r(\ln^2 n)\sqrt{n}\}
\]
Then
\[
 X_n+Y_n\le r(\ln^2 n)\sqrt{n}\ONE_{E}+ (X_n+Y_n)\ONE_{E^c}.
\]
Therefore,
\begin{equation}
 \E (X_n+Y_n)\le r(\ln^2 n)\sqrt{n}+\sqrt{\E X_{n}^2\Pp(E^c)}+\sqrt{\E Y_{n}^2\Pp(E^c)}.
\label{eq:M_t}
\end{equation}

Let us estimate the second term in~\eqref{eq:M_t}. According to Lemma~\ref{lem:shear}, the random variables $A^n(k)-\E A^n(k)$, $k\in
I_n$ have the same distribution, so replacing the maximum in the definition of $X_n^2$ with summation we
obtain
\begin{equation}
\E X_{n}^2\le (4Rn+1) \E (A^n-\E A^n)_+^2\le (4Rn+1) \E (A^n)^2\le Cn^3,
\end{equation}
for some $C>0$ and all $n\in\N$, where we used Lemma~\ref{lem:moments_of_At} in the last inequality.

Also, Lemma~\ref{lem:concentration_around_mean} shows that
\begin{align}
 \Pp(E^c)\le &\sum_{k\in \bar I_n}\Pp\left\{A^n(k)-\E A^n(k) > \frac{r}{2}(\ln^2 n)\sqrt{n}\right\}\notag \\ 
    &+\sum_{k\in \bar I_n}\Pp\left\{A^{n,2n}(k,0)-\E A^{n,2n}(k,0) > \frac{r}{2}(\ln^2 n)\sqrt{n}\right\}\label{eq:probability_of_D_complement}
\\ \le & 2(4R n+1) b_1\exp\left(-b_2\frac{r\ln n}{2} \right).\notag
\end{align}

The same estimates apply to the third term in~\eqref{eq:M_t}, and we can now finish the proof by choosing $r$ to be large enough and combining estimates
\eqref{eq:doubling+constant}--\eqref{eq:probability_of_D_complement}.
\epf

\medskip

With this lemma at hand we can now use a discrete version of~Lemma~4.2 from~\cite{HoNe}. 
Its proof literally repeats the proof of the original lemma where the argument of $a$ and $g$
was assumed to be continuous.

\begin{lemma}Suppose the functions $a: \N\to\R$ and $g:\N\to\R_+$  satisfy
the following conditions:
$a(n)/n\to \nu\in\R$ and  $g(n)/n \to 0$ as $n\to\infty$, $a(2n)\ge 2a(n)-g(n)$
and $\psi\equiv \limsup_{n\to\infty} g(2n) /g(n)< 2$. Then, for any $c > 1/(2-\psi)$,
and for all large $n$,
\[a(n) \leq \nu n + cg(n).\]
\end{lemma}
Taking $a(n)=\E A^n$, $\nu=\alpha(0)$, $g(n)=b_0n^{1/2}\ln^2n$, $\psi=\sqrt{2}$, $c=2$,
we conclude that for $b'_0=2b_0$ and large $n$,
\[
0 \le \E A^n-\alpha(0) n\le b_0'n^{1/2}\ln^2n,
\]
and Theorem~\ref{thm:concentration_around_alphat} follows from this estimate,
Lemma~\ref{lem:concentration_around_mean}, and the shear invariance established in Lemma~\ref{lem:shear}.\epf

\section{Straightness estimates}\label{sec:straightness}
 As in~\cite{BCK:MR3110798},\cite{CaPi}, and~\cite{CaPi-ptrf} we keep following the ideas from~\cite{HoNe} and~\cite{Wu}, adapting the program to our specific situation. The step that we make in this section is to estimate deviations of minimizers from straight lines. 

We will need a curvature estimate for the shape function constructed in Section
\ref{sec:subadd}. 
Recalling that the shape function $\alpha:\R\to\R$ was introduced in Lemma~\ref{lem:shape-function}, we define
$\alpha_0=\alpha(0)$ and extend $\alpha$ to a function on $\N\times\R$:
\[
\alpha(n,x)= n\alpha\left(\frac{x}{n}\right)=n\left(\alpha_0+\frac{1}{2}\left(\frac{x}{n}\right)^2\right)=\alpha_0 n +\frac{x^2}{2n},
\quad (n,x)\in\N\times\R.
\]
We remark that, in contrast with the situation in~\cite{BCK:MR3110798}, we do not know the sign of~$\alpha_0$.

Lemma~\ref{lem:shape-function} implies that for all $n\in\N$ and $x\in\R$,
\[ \lim_{m\to \infty} \frac{A^{0,mn}(0,mx)}{m} \stackrel{\rm a.s.}= \alpha(n,x).\]
We need a convexity estimate of this function $\alpha$. For $(n,x)\in \N\times \R$, $L>0$,
we define
\[ \calC(n,x,L) := \left\{ (m,y)\in\Z\times\R:\ m\in \{n+1,\ldots, 2n\}\ \mbox{and}\ \left|y - \frac{m}{n}\,x\right| \le
L\right\}.\]
So $\calC(n,x,L)$ is a parallelogram of width $2L$ with one pair of sides parallel to the $x$-coordinate axis and the
other one parallel to $[(n,x),(2n,2x)]$  (for any two points $p,q$ on the plane $\R\times\R\supset \Z\times\R$,
$[p,q]$ denotes the straight line segment connecting these two points):
\[ \partial_S^\pm\calC(n,x,L) := \left\{ (m,y)\in\Z\times\R:\ m\in \{n+1,\ldots, 2n\}\ \mbox{and}\  y - \frac{m}{n}\,x=\pm
L\right\}.\]
 The union of $ \partial_S^+\calC(n,x,L)$ and $ \partial_S^-\calC(n,x,L)$ is
\[ \partial_S\calC(n,x,L) := \left\{ (m,y)\in\Z\times\R:\ m\in \{n+1,\ldots, 2n\}\ \mbox{and}\ \left|y - \frac{m}{n}\,x\right|=
L\right\}.\]
The following lemma will play the role of Lemma 2.1 in \cite{Wu}.
\begin{lemma}\label{lem:convest}
For all $(n,x),(m,y)\in \N\times \R$, such that $m>n$,
we have
\begin{equation}\label{eq:convest}
\alpha(m-n,y-x) + \alpha(n,x) = \alpha(m,y) + \frac{n}{2m(m-n)} \left(y - \frac{m}{n}\,x\right)^2.
\end{equation}
\end{lemma}
\bpf
\begin{eqnarray*}
\alpha(m-n,y-x) + \alpha(n,x) & = & \alpha_0 (m-n) + \frac{(y-x)^2}{2(m-n)} + \alpha_0 n + \frac{x^2}{2n}\\
& = & \alpha_0 m + \frac{(y-x)^2}{2(m-n)} + \frac{x^2}{2n}\\
& = & \alpha(m,y)-\frac{y^2}{2m} + \frac{(y-x)^2}{2(m-n)} + \frac{x^2}{2n}.
\end{eqnarray*}
Identity~\eqref{eq:convest} now follows from a straightforward comparison of the algebraic expressions involved.

One may also argue that for fixed $m$, $\alpha(m-n,y-x) + \alpha(n,x) -  \alpha(m,y)$ is a quadratic function in $y$. The minimum of this function
equals $0$, and is attained at $\tilde y$ such that $(m,\tilde y)$ is a multiple of $(n,x)$, i.e.,  $\tilde y= mx/n$. The lemma follows by computing
the coefficient in front of $y^2$.
\epf

This deterministic convexity lemma, together with the concentration bound of Section~\ref{sec:concentration}, will help
us to show that minimizers cannot deviate  from a straight line too much. 

We will need one more auxiliary estimate:
\begin{lemma} \label{lem:increment-of-shape-function-in-time} Let $c_1,c_2>0$,  $0<n_1<n_2<n_1+c_1$, $|x|<c_2 n_1$. Then
\[
\alpha(n_1,x)-\alpha(n_2,x) <|\alpha_0|c_1+c_2^2c_1. 
\]
\end{lemma}
\bpf A straightforward computation gives:
\begin{align*}
\alpha(n_1,x)-\alpha(n_2,x)&=\alpha_0 n_1 +\frac{x^2}{2n_1}-\alpha_0 n_2 +\frac{x^2}{2n_2}\\
&\le |\alpha_0| |n_1-n_2| +\frac{x^2}{2}\frac{n_2-n_1}{n_1^2},
\end{align*}
and the lemma follows.
\epf

%For $(n,x)\in \N\times\R$, we define
%\begin{equation}
%K(n,x,R):= \{n\}\times[x,x+R].
%\label{eq:K-square}
%\end{equation}
For $(n,x),(m,y)\in \Z\times\R$ satisfying $n<m$, we denote by $\gamma^{(n,x),(m,y)}$ the rightmost point-to-point minimizer
connecting $(n,x)$ to $(m,y)$
and by $A((n,x),(m,y))=A^{n,m}(x,y)$ the associated optimal action.
For
$(n,x),(m,y)\in \Z\times\R$ with $n<m$, we define the events
\begin{align*}
 G^+((n,x),(m,y)) = \left\{ \exists \tilde{0}\in [0,1],\ \exists\tilde{y}>y:\ \gamma^{(0,\tilde{0}),(m,\tilde{y})}_n\in [x,x+1]
\right\},\\
G^-((n,x),(m,y)) = \left\{ \exists \tilde{0}\in [0,1],\ \exists\tilde{y}<y:\ \gamma^{(0,\tilde{0}),(m,\tilde{y})}_n\in [x,x+1]
\right\}.
\end{align*}
%and $G((n,x),(m,y))=G^-((n,x),(m,y))\cup G^+((n,x),(m,y))$.
These events say that there is a minimizer starting close to $0$ at time $0$, ending at time $m$ to the right (respectively, left) of $y$, and
 passing close to $x$ at time $n$. To estimate
 the probability
of this event, we first have to control the action to and from points close to $(n,x)$.

\begin{lemma}\label{lem:boundaction}
Suppose $(n,x), (m,y)\in \Z\times\R$ and $n<m$. Let $\tilde x\in [x,x+1]$, $\tilde{y}\in [y,y+1]$.
Then
\[
A^{n,m}(\tilde x,\tilde{y})\ge A^{n-1,m+1}(x,y) -F^*(n-1,x)-F^*(n,x)- F^*(m,y)- F^*(m+1,y)-1.
\]
If, additionally, $n+3\le m$, then
\[ 
A^{n,m}(\tilde x,\tilde{y})\le A^{n+1,m-1}(x,y)+ F^*(n,x) + F^*(n+1,x) +F^*(m-1,y) +F^*(m,y)+1.\]
\end{lemma}
\bpf We recall that we work with action defined by~\eqref{eq:gen-potential-action}.
The first inequality is a result of a direct comparison of the action of $\gamma^{(n-1,x),(m+1,y)}$ to that of
\[
\gamma_k =
\begin{cases}  x,& k=n-1,\\ \gamma^{(n,x),(m,y)}_k,& k\in\{n,\ldots,m\},\\ y,& k=m+1.
\end{cases}
\]
The second inequality 
is a result of a direct comparison of the action of 
\[
\gamma_k =
\begin{cases}  \tilde x,& k=n,\\ \gamma^{(n+1,x),(m-1,y)}_k,& k\in\{n+1,\ldots,m-1\},\\ \tilde y,& k=m,
\end{cases}
\]
to the action of $\gamma^{(n,\tilde x),(m,\tilde y)}$. \epf

\begin{lemma}\label{lem:dstraight1}
Fix $\delta\in(0,1/4)$ and $v>0$. There exist constants $c_1, c_2, M>0$
such that  for all $(n,x)$ with $n>M$ and $(m,y)\in \partial_S^\pm \calC(n,x,n^{1-\delta})$, 
with $m>n$,
we have
\[ \Pp(G^\pm((n,x),(m,y)))\le c_1\exp\left(-c_2n^{1/2-2\delta}/\ln n\right).\]
\end{lemma}
\bpf 
Let us consider only the case of $\partial_S^+\calC(n,x,n^{1-\delta})$ and $G^+((n,x),(m,y))$. The shear invariance implies
that it is sufficient to consider $x=-1$, so $[x,x+1]=[-1,0]$.

On $G^+((n,x),(m,y))$, there are numbers
 $\tilde{0}\in [0,1]$, $\tilde{x}\in [-1,0]$ and $\tilde y>y$ such
that
\begin{equation}
\label{eq:action-decomposes-in-two}
%A((0,\tilde{0}),(m,\tilde{y})) = A((0,\tilde{0}),(n,\tilde{x})) + A((n,\tilde{x}),(m,\tilde{y})).
A^{0,m}(\tilde{0},\tilde{y})=A^{0,n}(\tilde{0},\tilde{x})+A^{n,m}(\tilde{x},\tilde{y})
\end{equation}
Let $i=[y]$.
Let us first consider the case where $i\ge 4Rn$. Then
\begin{equation}
\frac{n}{m}i > Rm.
\end{equation}
By monotonicity of dependence of point-to-point minimizers on the endpoints,  $\gamma^{(0,0),(m,i)}_n<0$.
It means that on $[0,m]$, $\gamma^{(0,0),(m,i)}$ deviates from a straight line by $ni/m$. Denoting the latter
event by $B(n,m,i)$, we can use Lemma~\ref{lem:path_in_wide_rectangle_whp} to write:
\[
\Pp(B(n,m,i))\le C_2 \exp(-C_1 (ni/m)^2/m)\le C_2\exp(-C_1 i^2n^2/m^3).
\]
For $B(n,m,4Rn+)=\bigcup_{i\ge 4Rn} B(n,m,i)$, we obtain
\[
\Pp\left( B(n,m,4Rn+)\right)\le C_2\sum_{i\ge 4Rn}\exp(-C_1 i^2n^2/m^3).
\]
Since the first term in this series is bounded by $C_2\exp(-2C_1 R^2n)$,
and the ratio of two consecutive terms is bounded by $\exp(-C_1R)$, we conclude that
\begin{equation}
\label{eq:far-exits}
\Pp\left( B(n,m, 4Rn+)\right)\le C_3 \exp(-2C_1 R^2n).
\end{equation}

Let us now consider the case where $[\tilde y]=i\in(n^{1-\delta}-1, 4Rn)$.
Lemma~\ref{lem:boundaction} implies that
\begin{align*}
 A^{0,m}(\tilde 0,\tilde{y}) &\le A^{1,m-1}(0,i)+ F^*(0,0)+ F^*(0,1)+F^*(m-1,i)+F^*(m,i)+1,\\
 A^{0,n}(\tilde 0,\tilde{x}) &\ge A^{-1,n+1}(0,0)) -F^*(-1,0)-F^*(0,0) - F^*(n,0)-F^*(n+1,0)-1,\\
 A^{n,m}(\tilde x,\tilde{y}) &\ge A^{n-1,m+1}(0,i) -F^*(n-1,0)-F^*(n,0) - F^*(m,i)-F^*(m+1,i)-1.
\end{align*}
Along with~\eqref{eq:action-decomposes-in-two}, this implies
\begin{equation}
%\notag
A^{1,m-1}(0,i)-A^{-1,n+1}(0,0)-A^{n-1,m+1}(0,i)\ge -3- X(n, m,i), 
\label{eq:subadditivity-wrong-way}
\end{equation}
where
\begin{align*}
X(n, m, i)=&F^*(0,0)+ F^*(0,1)+F^*(m-1,i)+F^*(m,i) \\ &+F^*(-1,0)+F^*(0,0) + F^*(n,0)+F^*(n+1,0)\\&+F^*(n-1,0)+F^*(n,0) + F^*(m,i)+F^*(m+1,i).
\end{align*}
Let us now approximate the left-hand side of~\eqref{eq:subadditivity-wrong-way} using the extension of shape function $\alpha$
introduced in the beginning of this section.
Lemma~\ref{lem:increment-of-shape-function-in-time} implies that 
 there is a constant $L$ such that under the constraints we have imposed on $n,m$, and $i$,
\begin{equation}
\label{eq:approx-by-alpha1}
\alpha(m-2,i)<\alpha(m+2,i)+L
\end{equation}
and
\begin{equation}
\label{eq:approx-by-alpha2}
\alpha(n+2,0)>\alpha(n,0)-L.
\end{equation}
Since $m+2\le 3n$ and $m+2-n\le 2n$, Lemma~\ref{lem:convest} implies that
\begin{equation}
\label{eq:approx-by-alpha3}
\alpha(m-n+2,i) + \alpha(n,0) \geq \alpha(m+2,i) + \frac{i^2}{12n} .
\end{equation}
Combining \eqref{eq:approx-by-alpha1},\eqref{eq:approx-by-alpha2},\eqref{eq:approx-by-alpha3}  with \eqref{eq:subadditivity-wrong-way}, we obtain
\begin{align*}
&(A^{1,m-1}(0,i)-\alpha(m-2,i))-(A^{-1,n+1}(0,0)-\alpha(n+2,0))\\ &-(A^{n-1,m+1}(0,i)-\alpha(m-n+2,i))
+X(n, m,i)\\ 
\ge  & -3 -2L +\frac{i^2}{12n}\ge \frac{i^2}{15n}
\end{align*}
if $n>M$ for sufficiently large $M$.

Let us define events 
\begin{align*}
E_1(n,m,i)&=\left\{ A^{1,m-1}(0,i)-\alpha(m-2,i) \ge \frac{i^2}{60 n }\right\},\\
E_2(n,m,i)&=\left\{ A^{-1,n+1}(0,0)-\alpha(n+2,0) \le - \frac{i^2}{60 n }\right\},\\
E_3(n,m,i)&=\left\{ A^{n-1,m+1}(0,i)-\alpha(m-n+2,i) \le -\frac{i^2}{60 n }\right\},\\
E_4(n,m,i)&=\left\{ X(n,m,i) \ge \frac{i^2}{60 n }\right\}.
\end{align*}
We would like to show that for some constants
$c'_1,c'_2$,
\begin{multline}
\label{eq:4events}
\Pp(E_1(n,m,i))+\Pp(E_2(n,m,i))+\Pp(E_3(n,m,i))+\Pp(E_4(n,m,i))\\ \le c'_1\exp\left(- c'_2 \frac{i^2}{n^{3/2}\ln n}\right).
\end{multline}
Sufficient estimates for first two terms follow from Theorem~\ref{thm:concentration_around_alphat},
and for the
last one from~\eqref{eq:tail-requirement2}. Let us estimate $\Pp(E_3(n,m,i))$. By the shear invariance,
\[
\Pp(E_3(n,m,i))=\Pp\left\{ A^{m-n+2}-\alpha_0\cdot(m-n+2) \le - \frac{i^2}{60 n }\right\}.
\]
Since $m-n+2$ may be small, Theorem~\ref{thm:concentration_around_alphat} does not apply directly. However, to estimate
the right-hand side, we notice that $A^{m-n+2}\ge A^m- A^{m-n+2,m}$ and, therefore, 
\[
\Pp(E_3(n,m,i))\le \Pp\left\{ |A^{m}-\alpha_0 m| \geq \frac{i^2}{120 n }\right\}+\Pp\left\{ |A^{n-2}-\alpha_0(n-2)| \geq \frac{i^2}{120 n }\right\},
\]
and the desired estimate follows from an application of Theorem~\ref{thm:concentration_around_alphat} to
both terms on the right-hand side. So, \eqref{eq:4events} follows, and we obtain
\begin{align*}
       & \sum_{n^{1-\delta}-1< i < 4Rn} (\Pp(E_1(n,m,i))+\Pp(E_2(n,m,i))+\Pp(E_3(n,m,i))+\Pp(E_4(n,m,i)))
\\ \le &  c'_1 \exp\left(- c'_2 \frac{n^{2-2\delta}}{4n^{3/2}\ln n}\right) 
\sum_{i> n^{1-\delta}/2}\exp\left(- c'_2 \frac{i^2-n^{2-2\delta}/4}{n^{3/2}\ln n}\right),
\end{align*}
where the series factor in the right-hand side is uniformly bounded in $n$.
This, along with~\eqref{eq:far-exits} implies the theorem.
\epf

The above Lemma can be used to show that a minimal path starting close to the origin and passing  close to $(n,x)$, with
high probability will not exit the parallelogram $\calC(n,x,n^{1-\delta})$ through the lateral sides. 

\begin{lemma}\label{lem:dstraight2}
Fix $\delta\in(0,1/4)$ and $v>0$. There exist constants $c_1, c_2, M, \kappa>0$ 
such that  for all $(n,x)$ with $n>M$,
\[ \Pp(H(n,x))\le c_1\exp\left(-c_2n^\kappa\right),\]
where $H(n,x)= H^+(n,x)\cup H^-(n,x)$ and
\[
 H^\pm(n,x)=\bigcup_{(m,y)\in \partial_S^\pm \calC(n,x,n^{1-\delta})} G^\pm((n,x),(m,y)).
\]
\end{lemma}
\bpf
There are $n$ admissible values of $m$, so we can apply Lemma~\ref{lem:dstraight1} and choose any $\kappa\in(0,1/2-2\delta)$. 
\epf

Let us now prove $\delta$-straightness of minimizers, as was introduced by Newman in \cite{Ne}. For a path $\gamma$ and $n\in \Z$, we define
\[ \gamma^{\rm out}(n) = \{ (m,\gamma_m):\ m>n\}.\]
For $v>0$ we define
\[\Co(v) = \{ (n,x)\in\N\times \R:\ |x|\leq n v\}.\]
For  $(n,x)\in \N\times \R$ and $\eta>0$ we define:
\begin{equation}
\label{eq:Co}
{\rm Co}(n,x,\eta) = \{ (m,y)\in \N\times \R\ :\ |y/m - x/n| \le\eta \}.
\end{equation}

\begin{lemma}[{\bf $\delta$-straightness}]\label{lem:delta_straightness}
For
$\delta\in(0,1/4)$ and $v>0$ we have with probability one that there exists $M=M(\omega)>0$ (depending on $v$
and $\delta$) and nonrandom $Q>0$ (depending only on $\delta$), such that for all $\tilde{0}\in [0,1]$, for all $(m,y)\in
\N\times \R^+$ and for all 
$(k,z)\in \gamma^{(0,\tilde 0),(m,y)}\cap {\rm Co}(v)$ with $k>M$, we have
\[ \gamma^{\rm out}(k) \subset {\rm Co}(k,z,Qk^{-\delta}),\]
for $\gamma=\gamma^{(0,\tilde 0),(m,y)}$.
%\marginpar{$M$ is random}
\end{lemma}

This lemma states that if a minimizer starting near $(0,0)$ passes through a remote point $(k,z)$, it has to stay in a
narrow cone around the ray $\N\cdot (k,z) $.

\bpf  Using Lemma~\ref{lem:dstraight2} and the Borel--Cantelli Lemma, we see that for any $v'>v$,  there is a (random) $M>0$
such that if $j>M$ and $(j,i)\in \Co(v')$ then  $H^c(j,i)$ holds. 

We conclude that any minimizer $\gamma$ passing through $\{0\}\times[0,1]$ and $\{j\}\times [i,i+1]\subset \Co(v')$ with $j>M$
satisfies
\begin{equation}
\left|\gamma_n - \frac{ni}{j} \right|< j^{1-\delta},\quad n\in \{j+1,\ldots,2j\}. 
\label{eq:deviation-in-upper-part-of-rectangle} 
\end{equation}
So for any such minimizer and point $(k,z)$ satisfying the conditions of the lemma and $M$ chosen as above,
we can choose numbers $i_0=[z], i_1,i_2,\ldots$ such that for each $l$, 
$\gamma_{2^l k}\in [i_l,i_l+1]$,
and
\[
\left|\gamma_n - \frac{ni_l}{2^lk} \right|< (2^lk)^{1-\delta},\quad n\in \{2^l k+1,2^l k+2,\ldots 2^{l+1} k\}. 
\]
In particular, for all $l$
\[
\left|i_{l+1} - 2i_l \right|< (2^lk)^{1-\delta}+1 \le 2(2^lk)^{1-\delta}.
\]
Therefore, for all $l$,
\[
|i_l-2^l i_0|\le \sum_{j=0}^{l-1} 2^{l-j-1}|i_{j+1}-2i_j|\le 
\sum_{j=0}^{l-1} 2^{l-j-1}2 (2^jk)^{1-\delta}\le 
 2^{l}k^{1-\delta}\sum_{i=0}^{l-1}2^{-\delta j},
\]
so
\begin{equation}
\left|\frac{i_l}{2^lk}-  \frac{i_0}{k}\right|\le c k^{-\delta}
\label{eq:displacement-for-powers-of-2} 
\end{equation}
for some constant $c$.

For $n\in [2^lk,2^{l+1}]$, $l\ge 1$, \eqref{eq:deviation-in-upper-part-of-rectangle} implies
\begin{equation}
\label{eq:displacement-for-powers-of-2a}
 \left|\gamma_n - \frac{ni_l}{2^lk} \right|< (2^lk)^{1-\delta}\le n^{1-\delta}k^{1-\delta}.
\end{equation}
Now we obtain the lemma by combining \eqref{eq:displacement-for-powers-of-2} and \eqref{eq:displacement-for-powers-of-2a}.
\epf

As a side product of the proof of Lemma~\ref{lem:delta_straightness} we obtain the following statement:
\begin{lemma}\label{lem:quantitative-straightness}
 For
$\delta\in(0,1/4)$ and $v>0$ there are nonrandom positive numbers $M,Q,C_1,C_2,\kappa$, such that if
\begin{multline*}
 G_n=\Bigl\{\exists \tilde 0\in[0,1],\ (m,y)\in\Z\times\R,\ (k,z)\in\gamma^{(0,\tilde 0),(m,y)}\cap\Co(v):  \\
k>n,\ \gamma^{\rm out}(k) \not\subset {\rm Co}(k,z,Qk^{-\delta})\Bigr\},
\end{multline*}
then for $n\ge M$,
\[
 \Pp(G_n)\le C_1e^{-C_2n^\kappa}.
\]
\end{lemma}
\bpf Arguing as in the proof of Lemma~\ref{lem:delta_straightness}, we obtain that 
if $G_n$ holds, then at least one event $H(j,i)$ is violated for $j\ge n$, and the lemma follows now from 
Lemma~\ref{lem:dstraight2}.
\epf

\section{Existence and uniqueness of semi-infinite minimizers}\label{sec:one-sided_minimizers}
\subsection{Existence.}
With $\delta$-straightness at hand, we can prove some important properties of minimizing paths.
A semi-infinite minimizer starting at  $(n,x)\in\Z\times\R$ is a path $\gamma:\{n,n+1,\ldots\}\to \R$ such that
$\gamma_n=x$ and the restriction of $\gamma$ to any finite time interval is
a point-to-point minimizer. We call  $(n,x)$ the endpoint of $\gamma$.

\begin{lemma}\label{lem:geodesic direction}
With probability one, all semi-infinite minimizers have an asymptotic slope (velocity, direction): for every
minimizer $\gamma$, there
is $v\in \R\cup\{\pm \infty\}$  depending on $\gamma$ such that
\[ \lim_{n\to\infty} \frac{\gamma_n}{n} = v.\]
\end{lemma}
\bpf Let us fix a sequence $v_n\to \infty$. Using the translation invariance of the forcing potential $F$,
with probability one, for any $(j,i)\in \Z^2$ we can choose  a
corresponding sequence of constants $M_n(j,i)>0$ such that the statement in Lemma~\ref{lem:delta_straightness} holds for
the entire sequence, for paths starting in $\{j\}\times[i,i+1]$.

Let us take some one-sided minimizer $\gamma$. If $\gamma_k/k\to+\infty$ or $-\infty$, then the desired statement
is automatically true. In the opposite case we have
\[ \liminf_{k\to\infty} \frac{|\gamma_k|}{k} < \infty.\]
This implies that there exist $n\geq 1$ and a sequence $k_m\to \infty$ such that $|\gamma_{k_m}|/k_m \leq v_n$. Let us
and choose $(j,i)\in\Z^2$ such
that $(k_1,\gamma_{k_1})\in j\times[i,i+1]$. For $m$ large enough, we will have that $k_m>M_n(j,i)$ and, therefore,
\[ \gamma^{\rm out}(k_m)\subset  (j,i) + {\rm Co}(k_m-j, \gamma_{k_m}-i,Q|k_m-j|^{-\delta}),\]
for a constant $Q>0$ and $m$ large enough. Therefore,
$\bar v=\limsup_{k\to\infty}\gamma_k/k$ and $\underline{v}=\liminf_{k\to\infty}\gamma_{k}/k$ are
well-defined and satisfy 
\[
\bar v -\underline{v}<Q|k_m-j|^{-\delta}.
\] 
Since the right-hand side converges to $0$ as $m\to\infty$, we obtain $\bar v =\underline{v}$.
\epf

\begin{lemma} Let $\gamma$ be a point-to-point minimizer between $(n,x)$ and $(m,y)$. Then for all $k\in\{n+1,\ldots,m-1\}$,
\begin{equation}
\label{eq:euler-lagrange-0}
\gamma_{k+1}-\gamma_k=\gamma_k-\gamma_{k-1}+f(k,\gamma_k),
\end{equation}
or, equivalently,
\begin{equation}
\label{eq:euler-lagrange-1}
\gamma_{k+1}=2\gamma_k-\gamma_{k-1}+f(k,\gamma_k),
\end{equation}
where $f(k,x)=\partial_x F(k,x)$ for all $(k,x)\in\Z\times\R$.
\end{lemma}
\bpf
For all $\omega\in\Omega$, $F_\omega$ is a sum of finitely many smooth functions on every interval $\{k\}\times(x-1,x+1)$. Therefore,
action is a smooth function of paths. The lemma follows from equating the partial derivative of the action with respect to $\gamma_k$ to $0$.
\epf

Equations~\ref{eq:euler-lagrange-0} and~\eqref{eq:euler-lagrange-1} are discrete time versions of Euler--Lagrange equations. Their meaning is
that path $\gamma$ instantaneously changes its velocity by $f(k,\gamma_k)=\partial_x F(k,\gamma_k)$ at time $k$.

The following is the last technical lemma we need to prove existence of one-sided minimizers.
\begin{lemma}\label{lem:one-more-straightness-lemma} Let $n_i\uparrow \infty$ and $x_i,z_i$ satisfy $|x_i|<|z_i|$ and $|x_i-z_i|>2Rn_i$ for all $i\in\N$.  
 Then with probability $1$, there is $N(\omega)$ such that if $n_i>N(\omega)$, then
 \[
|\gamma^{(0,0),(n_i+1,x_i)}_{n_i}| < z_i.
 \]
\end{lemma}
\bpf This Lemma is a direct consequence of Lemma~\ref{lem:path_in_wide_rectangle_whp} and the Borel--Cantelli Lemma.
\epf

The next lemma states existence of one-sided minimizers with a given asymptotic slope and provides a way to construct them.
\begin{lemma}\label{lem:exist-geodesic}
With probability one, for every $v\in\R$ and for every sequence $(m_i,y_i)\in \Z\times\R$ with $m_i\to \infty$ and
\[ \lim_{i\to \infty} \frac{y_i}{m_i} = v,\]
and for every $(n,x)\in \ZR$, there exists a subsequence $(i_k)$ such that the minimizing paths
$\gamma^{(n,x),(m_{i_k},y_{i_k})}$ converge pointwise to a semi-infinite minimizer starting
at $(n,x)$ and with asymptotic slope equal to $v$.
\end{lemma}
\bpf Without loss of generality, we can assume that $(n,x)\in\{0\}\times[0,1]$ since
$\ZR$ can be represented as a countable union of shifts of this set. 
We fix $\delta\in(0,1/4)$  and a sequence $v_l\to \infty$, and then
choose $Q>0$ and $M_l=M_l(\omega)\to \infty$ such that the statement of Lemma~\ref{lem:delta_straightness} holds for every
triplet $(v_l,M_l,Q)$. From these triplets we choose a triplet $(v_0,M,Q)$ such that $v_0>v+2QM^{-\delta}+2R$.

Passing to a subsequence if needed, we can make sure that $(m_i)$ is an increasing sequence satisfying
 $m_i\ge M$ for all $i$,  and $(m_j,y_j)\in {\rm
Co}(1,v,Qm_{i}^{-\delta})$ for all $j>i$. 

Consider the paths $\gamma^j = \gamma^{(0,x),(m_j, y_j)}$. 
%We claim that there is $C>0$ 
%such that for any $t_0>0$ there is $t\in\N$
%such that $t>t_0$ and infinitely many paths $\gamma^j$ satisfy
%\begin{equation}
%\label{eq:minimizers-in-cone}
%\max \{|\gamma^j_t-vt|, |\gamma^j_{t+1}-vt|\}<Ct^{1-\delta}.
%\end{equation}
We claim that there is $t\in\N$ and $C>0$ such that infinitely many paths $\gamma^j$ 
satisfy
\begin{equation}
\label{eq:minimizers-in-cone}
|\gamma^j_t|, |\gamma^j_{t+1}|<C.
\end{equation}
To see that, let us  analyse the restrictions on paths $\gamma_j$ that we have.

The situation where $\gamma^j$ visits a point $(p,z)\in {\rm Co}(v_0)\setminus {\rm Co}(1,v,2Qm_i^{-\delta})$
satisfying $p>M$  is impossible since
this path would violate the $\delta$-straightness condition: the relevant cone through
$(p,z)$ will not overlap with ${\rm Co}(1,v,Qm_i^{-\delta})$, and therefore it cannot contain $(m_j,y_j)$.

Suppose now that $\gamma^j$ does not visit ${\rm Co}(v_0)\setminus {\rm Co}(1,v,2Qm_i^{-\delta})$ but 
visits  ${\rm Co}(v_0)^c$. Then there are points $(p,z_1)\notin {\rm Co}(v_0)$ and  $(p+1,z_2)\in{\rm Co}(1,v,2Qm_i^{-\delta})$
such that $\gamma^{(0,x),(p+1,z_2)}_p=z_1$. Let us consider the case $z_1>v_0 p$ (the case  $z_1<-v_0 p$ is treated similarly).
Due to monotonicity of dependence of minimizers on endpoints, $\gamma^{(0,1),(p+1, (v+2Qm_i^{-\delta})(p+1))}_p>v_0p$.
Now Lemma~\ref{lem:one-more-straightness-lemma} implies that this can happen only for finitely many
values of $p$. So there exists $N$ such that no path $\gamma^j$ visits ${\rm Co}(v_0)^c$ after time $N$.

Combining the claims of last two paragraphs, we obtain that for each $i\ge 1$
and $j>i$, the path $\gamma^j$  lies in the cone ${\rm Co}(1,v,2Qm_i^{-\delta})$ for times larger than
$m_i\vee N$. In particular, our claim about\eqref{eq:minimizers-in-cone} holds true.

Using compactness of $[-C,C]$, we can find a subsequence of minimizers $\gamma_j$ 
converging to a limit at times $t$ and $t+1$.  Let us recall that
 minimizers solve the
Euler--Lagrange equation~\eqref{eq:euler-lagrange-1}. The values of a 
solution of~\eqref{eq:euler-lagrange-1} at $t$ and $t+1$
determine the entire solution uniquely. Also, due to continuous dependence of solutions on the initial data,
the convergence of the solution at times $t,t+1$ implies the pointwise convergence at all times. 
Due to Lemma~\ref{lem:limit-of-minimizers}, the resulting limiting infinite trajectory is a one-sided
minimizer. Since all the pre-limiting finite minimizers lie within the cones  ${\rm Co}(1,v,2Qm_i^{-\delta})$
for sufficiently large times, so does the limiting one-sided minimizer. Therefore, its asymptotic slope equals $v$.\epf

\subsection{Uniqueness.} Let us now prove uniqueness of one-sided minimizers with given asymptotic slope.
\begin{lemma}
\label{lem:minimizers-do-not-intersect}
 Let $\omega\in\Omega$ and $(n,x)\in\ZR$. If $\gamma^1$ and $\gamma^2$ are two distinct one-sided minimizers with endpoint
$(n,x)$, then they do not intersect as curves in in $\R\times\R$.
\end{lemma}
\bpf This follows directly from Lemma~\ref{lem:finite-horizon-minimizers-do-not-intersect}. \epf

\begin{lemma} \label{lem:uniqueness-of-minimizer-for-1-point}
Let $(n,x)\in\RZ$ and $v\in\R$. With probability 1, there is a unique one-sided minimizer with slope $v$ and endpoint $(n,x)$.
\end{lemma}
\bpf
We fix $(n,x)$. First, we know that with probability~1, for each $v\in\R$, there is a one sided minimizer with asymptotic slope~$v$. 
To each direction $v\in\R$ we assign an interval $I_v=(a_v,b_v)$, where $a_v=\inf \gamma_1$ and
$b_v=\sup \gamma_1$
with infimum and supremum taken over all one-sided minimizers $\gamma$ with slope $v$. 
If there is a unique one-sided minimizer with slope $v$,  we set $I_v=\emptyset$.
Lemma~\ref{lem:minimizers-do-not-intersect} implies that if $v_1<v_2$, then every minimizer $\gamma^1$  with slope $v_1$ and every minimizer $\gamma^2$  with slope $v_2$ satisfy $\gamma^1_k<\gamma^2_k$,
for all times $k>n$. In particular $I_{v_1}\cap I_{v_2}=\emptyset$. Since one can place at most countably many disjoint open intervals in $\R$,
we have that there are at most countably many values of $v$ such that $I_v\ne\emptyset$, i.e.,  $a_v\ne b_v$. However, 
$\Pp\{I_v\ne\emptyset\}=p$ does not depend on $v$ due to shear invariance. Therefore, we can take any probability density $f$ on $\R$
and write
\[
p=\int_\R \Pp\{I_v\ne\emptyset\}f(v)dv=\int_\R \E \ONE_{\{I_v\ne\emptyset\}}f(v)dv=\E \int_\R \ONE_{\{I_v\ne\emptyset\}}f(v)dv =0,
\]
since $I_v\ne\emptyset$ can be true at most for countably many $v$. So, for any $v\in\R$, $\Pp\{I_v\ne\emptyset\}=0$, and the lemma follows.
\epf

\begin{lemma}\label{lem:limits-of-minimizers-uniquely-defined} 
Under the conditions of the previous lemma, let $\gamma$ denote the (a.s.-unique) minimizer with endpoint $(n,x)$ and slope $v$. Then
there is an event of probability 1 such that on that event, for any sequence of points $(m_i,y_i)\in\ZR$
such that $m_i\to\infty$ and $y_i/m_i\to v$ and all $k>n$, we have $\gamma^{(n,x),(m_i,y_i)}_k\to\gamma_k$.
\end{lemma}
\bpf Let us check that for any sequence $(i')$
we can choose a subsequence $(i'')$ such that the corresponding minimizers $\gamma^{(n,x),(m_{i''},y_{i''})}$ converge to $\gamma$.
In fact, Lemma~\ref{lem:exist-geodesic} allows to find a  subsubsequence $(i'')$
such that the corresponding minimizers $\gamma^{(n,x),(m_{i''},y_{i''})}$ converge pointwise to a limiting
infinite one-sided minimizer. However, there is a unique one-sided minimizer~$\gamma$, and the desired convergence follows.
\epf

\begin{lemma} \label{lem:minimizers-do-not-intersect-if-one-is-unique}
Let $(n,x)\in\ZR$  and $v\in\R$. With probability $1$, the following holds true:  there is a unique minimizer $\gamma$ with endpoint  $(n,x)\in\ZR$  and slope $v\in\R$;
if $y\ne x$, then no minimizer $\tilde \gamma$ with endpoint $(n,y)$ and slope $v$ can intersect $\gamma$.
\end{lemma}

\bpf Without loss of generality let us assume that $y>x$. Minimizers $\gamma$ and $\tilde \gamma$ have at most one intersection point (the argument for this claim is similar to the proof of Lemma~\ref{lem:finite-horizon-minimizers-do-not-intersect}).
Suppose they have exactly one intersection point. Then there is 
an integer $m>n$ such that $\tilde\gamma$ and $\gamma$ intersect strictly between times $m-1$ and $m+1$.  
 
 Let us consider the sequence of minimizers $\beta^k=\gamma^{(n,x),(k,\tilde \gamma_k)}$.
Lemma~\ref{lem:limits-of-minimizers-uniquely-defined} implies that these minimizers converge pointwise to $\gamma$. In particular,
$\beta^k_{j}\to \gamma_{j}$, $j\in\{m-1,m,m+1\}$.
Therefore, for sufficiently large $k$ the point-to-point 
minimizers $\beta^k$ and $\tilde \gamma$ coincide at time $k$ and intersect
between times $m-1$ and $m+1$. This is a contradiction with Lemma~\ref{lem:finite-horizon-minimizers-do-not-intersect},
and the proof is completed.
\epf

\begin{lemma} Let $v\in\R$. With probability 1 the following holds true: for all $(n,x)\in\ZR$ there is a one-sided minimizer with
endpoint $(n,x)$ and slope~$v$; the minimizers are unique for all $(n,x)$ except countably many.
\end{lemma}
\bpf Let us fix $v$. Lemma~\ref{lem:limits-of-minimizers-uniquely-defined} implies that with probability 1, unique minimizers $\gamma(n,q)$ with slope~$v$ exist simultaneously for all points $(n,q)\in\Z\times\Q$. Let us take a point $(n,x)\in\ZR$ with $x\notin\Q$. Due to Lemma~\ref{lem:minimizers-do-not-intersect-if-one-is-unique},
none of $\gamma(n,q)$ with  $(n,q)\in\Z\times(\Q\cap [[x],[x]+1])$ intersect each other. In particular, they all are squeezed between
$\gamma(n,[x])$ and  $\gamma(n,[x]+1)$. So we can take a sequence of rational points $q_i\downarrow x$ and 
use Lemma~\ref{lem:limit-of-minimizers} to
conclude that the pointwise limit of $\gamma(n,q_i)$ is a minimizer with endpoint $(n,x)$. Since this minimizer $\gamma^+(n,x)$
lies between $\gamma(n,[x])$ and  $\gamma(n,[x]+1)$, its asymptotic slope is also equal to $v$.

We can also construct another minimizer $\gamma^-(n,x)$ for any point $(n,x)$ taking left limits $q_i\uparrow x$. The set of
all points $(n,x)$ such that these two minimizers do not coincide is at most countable set since each discrepancy defines an interval and these
intervals are disjoint. Lemma~\ref{lem:minimizers-do-not-intersect-if-one-is-unique} implies that if the limits on the right and on the left coincide, then the constructed minimizer is unique, so there is at most countable set of shock points where the minimizer is not unique.
\epf

A useful point of view at the field of minimizers is via families of monotone maps associated with them. Namely,
for each $n\in\Z$ we can consider a monotone right-continuous map defined by $x\mapsto \gamma^+_{n-1}(n,x)$. This map has
at most countably many discontinuities and, in terms of fluid dynamics, takes the particle
that arrives to point $x$ at time $n$ and outputs its position at time $n-1$. Of course, one can also consider the inverses of those
maps. They are also monotone maps that have intervals of constancy corresponding
to particles absorbed into shocks. Compositions of maps from both families correspond to particle dynamics over longer time intervals.

\section{Weak hyperbolicity}\label{sec:weak-hyperbolicity} Here we switch back to the point of view where the action is given by~\eqref{eq:action} and
one-sided minimizers are ``backward'' ones. For an endpoint $(n,x)$ a one-sided minimizer $\gamma$ is defined on $\{\ldots,n-1,n\}$. Let us recall that the guiding idea is to construct global solutions of the Burgers equation using these minimizers collecting information about the history of the forcing before time $n$.

In contrast with the existing work on the Burgers equation and last-passage type percolation models, for our model, we are not able to prove
hyperbolicity, i.e., the asymptotic closeness of one-sided minimizers in reverse time. However, in this section we prove a weakened hyperbolicity property
that is still useful for our model and hopefully will be useful for other related models. The most important new ideas of this
paper are introduced in this section and the following ones.

From the previous section we know that for a fixed $v\in\R$ , with probability~$1$, to each point $(n,x)\in\ZR$, we can assign the rightmost
minimizer $\gamma^+(n,x)$ and leftmost minimizer $\gamma^-(n,x)$ with slope $v$. Throughout this section we will be working with $v=0$
(although the results are valid for all $v\in\R$ due to the Galilean invariance of the system),
so we suppress the dependence on $v$ in the notation. For all but countably many $(n,x)$, $\gamma^+(n,x)=\gamma^-(n,x)$.

Let us define for $n\in\Z$ and $x,y\in\R$ satisfying $x\le y$,
\[
W_k(n,x,y)=\sum_{i=0}^{2}(\gamma^+_{k-i}(n,y)-\gamma^-_{k-i}(n,x)),\quad k\le n.
\]
We also define $W_k(n,x)=W_k(n,x,x)$ for $(n,x)\in\ZR$. Let us state the main result of this section.

\begin{theorem}\label{thm:weak-hyperbolicity} With probability $1$, for all $n\in\Z$ and $x,y\in\R$ satisfying $x\le y$,
\begin{equation}
\liminf_{k\to-\infty}\frac{W_k(n,x,y)}{(n-k)^{-1}}=0.
\label{eq:weak-hyperbolicity-different-endpoints}
\end{equation}

\end{theorem}

This theorem is a corollary of the following fact:
\begin{lemma}\label{lem:weak-hyperbolicity} With probability $1$, for all $(n,x)\in\ZR$,
\begin{equation}
\liminf_{k\to-\infty}\frac{W_k(n,x)}{(n-k)^{-1}}=0.
\label{eq:weak-hyperbolicity-one-endpoint}
\end{equation}
\end{lemma}
\bpf[Derivation of Theorem~\ref{thm:weak-hyperbolicity} from Lemma~\ref{lem:weak-hyperbolicity}]
Suppose that with positive probability there are $n,x,y$ such 
that~\eqref{eq:weak-hyperbolicity-different-endpoints} fails.
Our goal is to prove that then with positive probability there is $(n',x')$ such that~\eqref{eq:weak-hyperbolicity-one-endpoint}
fails for $(n,x)$ replaced with $(n',x')$.

Due to space-time stationarity we obtain that there is $r>0$ with the following property: with positive probability
there are $x,y\in[-r,r]$ such that $x<y$, \eqref{eq:weak-hyperbolicity-different-endpoints} fails, and the unique one-sided
minimizers $\gamma^x$ for $(0,x)$ and $\gamma^y$ for $(0,y)$ satisfy $-r <\gamma^x_{-1}<\gamma^y_{-1}<r$. We denote this event by $A$.

Let $B$ be the following event:  $F_\omega(0,w)=0$ for all $w\in [-r,r]$ and 
$F_\omega(0,-r-1)\vee F_\omega(0,r+1)<-3(r+1)^2$.
Then $P(B)>0$ due to~\eqref{eq:xi-has-values-if-both-signs} and the fact that Poisson random variables are unbounded. 
Also, events $A$ and $B$ are independent since $A$ depends only on the realization of $\omega$ for negative times. Therefore, $\Pp(A\cap B)>0$.

Let us prove that if $A$ and $B$ hold, then  no minimizer with endpoint $(1,z)$ for some $z\in\R$ can pass between $x$ and $y$.

Suppose that the opposite holds, namely, there is $z\in\R$ such that $p=\gamma_0(1,z)$ satisfies $x<p<y$. Then $\gamma(1,z)$ has
to pass between $\gamma^x$ and $\gamma^y$, so, denoting  $q=\gamma_{-1}(1,z)$, we obtain $-r<q<r$. We claim that,
contrary to our assumption,
the path $(\gamma_{-1}{(1,z)},\gamma_0{(1,z)},\gamma_1{(1,z)})= ((-1,q),(0,p),(1,z))$ is not 
a minimizing path between $(-1,q)$ and $(1,z)$.
Namely, we claim that one of the the paths $\gamma^+=((-1,q),(0,r+1) ,(1,z))$ and  $\gamma^-((-1,q),(0,-r-1),(1,z))$ has smaller action on time interval $\{-1,0,1\}$
than~$\gamma(1,z)$. 
Suppose $z>0$, then, noticing that the contribution of $(-1,q)$ is common for all paths under consideration and using
the definition of $B$, we obtain
\begin{align*}
A^{-1,1}(\gamma(1,z))-A^{-1,1}(\gamma^+)= & \frac{1}{2}(p-q)^2+\frac{1}{2}(z-p)^2+F_\omega(0,p) \\ &- \frac{1}{2}(r+1-q)^2-\frac{1}{2}(z-r-1)^2-F_\omega(0,r+1)
\\ >& (r+1-p)z \ge 0,
\end{align*}
and, similarly, for $z<0$,
\[
A^{-1,1}(\gamma(1,z))-A^{-1,1}(\gamma^-)> 0.
\]
In both cases we obtain a contradiction that proves our claim that on
$A\cap B$ no one-sided
minimizer with endpoint after time $0$, can pass between~$\gamma^x$ and~$\gamma^y$. 

Due to stationarity, it is impossible for all
minimizers with endpoint after time $0$ to pass on the left of $x$ at time $0$, or for all of these minimizers to pass on the right of $y$ at time $0$. Therefore, due to the monotonicity of the minimizers with respect to the endpoint, there is $z\in\R$ such that for all $z'<z$,
$\gamma_0(1,z')\le x$ and for all $z'>z$, $\gamma_0(1,z')\ge y$. In particular $(1,z)$ is a shock point, i.e., there are
two one-sided minimizers with endpoint $(1,z)$, one passing on the right of $(0,y)$ 
and another on the left of $(0,x)$. Since~\eqref{eq:weak-hyperbolicity-different-endpoints} fails
and $W_k(1,z)>W_k(0,x,y)$ for all $k<0$, we conclude that on the event $A\cap B$ of positive probability, there is a point $(1,z)$
such that~\eqref{eq:weak-hyperbolicity-one-endpoint} fails with $(n,x)$ replaced by $(1,z)$. This  contradicts
Lemma~\ref{lem:weak-hyperbolicity}, and the proof is completed.
\epf

For the proof of Lemma~\ref{lem:weak-hyperbolicity}, we need some notation and auxiliary results.

For every shock $(n,x)$ we define the area absorbed in this shock by
\[
\Lambda(n,x)=\{(m,y):\  m\le n,\  \gamma^-_m(n,x)\le y\le \gamma^+_m(n,x)\}.
\]
If $(m,y)\in \Lambda(n,x)$, we write $(n,x)\prec (m,y)$  and say that $(n,x)$ absorbs $(m,y)$ or inherits from $(m,y)$.
In that case $(n,x)$ is called a  successor of $(m,y)$, and $(m,y)$ is called a predecessor of $(n,x)$. We have chosen
the notation to ensure that every shock counts as its own predecessor.

\begin{lemma} Every shock $(n,x)$ has a unique successor at time $n+1$, i.e., there is a unique shock $(n+1,z)$ such that $(n+1,z)\prec (n,x)$.
\end{lemma}
\bpf As in the derivation of Theorem~\ref{thm:weak-hyperbolicity} from Lemma~\ref{lem:weak-hyperbolicity}, it is not possible for all minimizers corresponding to points $(n+1,z'), z\in\R$ to pass on one side of $(n,x)$.
So there is a point $z$ such that for all $z'<z$,
$\gamma_{n}(n+1,z')\le x$ and for all $z'>z$, $\gamma_n(n+1,z')\ge x$. Moreover,
$\gamma_n(n+1,z')$ cannot be continuous at $z'=z$, since otherwise, due to the Euler--Lagrange 
equation~\eqref{eq:euler-lagrange-1}, $\gamma_{n-1}(n+1,z')$ would also be continuous at $z'=z$, in contradiction 
with the presence of shock at $(n,x)$. So, $(n+1,z)$ is a shock point and $(n+1,z)\prec (n,x)$.
\epf

Two shocks $(n,x)$ and $(n,y)$ are said to {\it merge} or {\it coalesce} if they have a common successor at time $n+1$, i.e., there is $z\in\R$ such that $(n+1,z)\prec(n,x) $ and $ (n+1,z)\prec(n,y)$.

\begin{lemma}
\label{lem:no-systematic-slope}
 Let $N,M\in\N$, $\delta>0$. Then $P(B_\delta)=0$, where
\begin{multline*}
B_\delta=\Bigl\{\text{\rm there is } r_0>0:  
\text{\rm for all}\ r>r_0,\ \max_{1\le i\le M}\ \max_{-N\le k\le i} \gamma^+_k(i,r)>(1+\delta)r\Bigr\}.
\end{multline*}
\end{lemma}
\bpf The lemma follows directly from the fact that
\[
X(r)=\max_{1\le i\le M}\ \max_{-N\le k\le i} \gamma^+_k(i,r)-r,\quad r\in\R,
\]
is a stationary process.
\epf

%\begin{align*}
%\rho^+(A)&=\lim_{r\to+\infty} \frac{|A\cap [0,r]|}{r}\in[0,+\infty],\\
%\rho^-(A)&=\lim_{r\to+\infty} \frac{|A\cap [-r,0]|}{r}\in[0,+\infty],
%\end{align*}
%
% If these limits are equal to each other, we define $\rho(A)$ to be equal to their
%common value.

\bpf[Proof of Lemma~\ref{lem:weak-hyperbolicity}] Suppose that with positive probability there is a point $(n,x)$ such that 
\begin{equation}
\label{eq:assumption-on-existence-of-fat-shock}
\liminf_{k\to-\infty}\frac{W_k(n,x)}{(n-k)^{-1}}>0.
\end{equation}

Then there is a number $q>0$ such that with positive probability there is a point $(n,x)$ satisfying
\begin{equation}
\label{eq:condition-on-non-decay-of-width}
W_k(n,x) >\frac{q}{n-k},\quad k<n.
\end{equation}
Our goal is to show that the systematic presence of such shocks all over the space-time leads to a contradiction.
The strategy is to show how to construct a family of shocks with disjoint absorbed areas and such that the union
of those areas is too large to fit into the space-time.

For any set $A\subset\R$ we introduce its density
\[
\rho(A)=\lim_{r\to+\infty} \frac{|A\cap [0,r]|}{r}\in[0,+\infty],
\]
if the limit is well-defined. Here bars mean the number of elements in the set.

Let $A_n$, $n\in\N$ be the (random) set of points $(n,x)$ satisfying~\eqref{eq:condition-on-non-decay-of-width}. 
By ergodic theorem, there is a deterministic number $a>0$ such that $\rho(A_n)=a$  almost surely for all $n\in\Z$.

Since the shock areas (areas between leftmost and rightmost one-sided minimizers)
are disjoint for different points in $A_n$, we see that
\[
\sum_{i=0}^2 (\gamma_{n-1-i}^+(n,r)-\gamma_{n-1-i}^-(n,0)) \ge
\sum_{x\in A_n\cap[0,r]}W_{n-1}(n,x)\ge q|A_n\cap[0,r]|.
\]
Since $|A_n\cap[0,r]|/r\to a$, we conclude that
\[
\liminf_{r\to\infty}\frac{\sum_{i=0}^2\gamma_{n-1-i}^+(n,r)}{r}\ge qa. 
\] 
Due to Lemma~\ref{lem:no-systematic-slope}, this is impossible unless $a\le q^{-1}<\infty$.

For $i\in\N$, let $B_i$ consist of all shocks in $A_i$ that are successors of some shocks in $A_{i-1}$.
Let $C_i$ consist of those shocks in $A_i$ that are not successors of any shocks in $A_{i-1}$. Let
$D_i$ consist of those shocks in $A_{i-1}$ whose successors are not in $A_i$. Let $E_i$ consist
of all shocks $(i,x)$ in  $A_{i-1}$ such that there is $y<x$ satisfying $(i,y)\in A_{i-1}$ and both
$(i,x)$ and $(i,y)$ have the same successor belonging to~$A_i$.  Let $F_i$ be the set of all shocks $(i-1,x)$ in $A_{i-1}$ such that no other point $(i-1,y)\in A_{i-1}$
has the same successor in $A_i$ as $(i-1,x)$  and satisfies $y<x$. 
We then have
\begin{equation}
\label{eq:density-identity-1}
\rho(A_i)=\rho(B_i)+\rho(C_i),
\end{equation}
and
\begin{equation}
\label{eq:density-identity-2}
\rho(A_{i-1})=\rho(D_i)+\rho(E_i)+\rho(F_i).
\end{equation}

To see this, first notice that all of the densities involved are well-defined deterministic numbers due to ergodic theorem, since they are defined via skew-invariant (w.r.t.\ spatial shifts) functionals of the ergodic environment. 
Then identities~\eqref{eq:density-identity-1} and~\eqref{eq:density-identity-2} 
follow from the fact that $A_i=B_i\cup C_i$ and $A_{i-1}= D_i\cup E_i\cup F_i$.

We already know that  $\rho(A_i)=\rho(A_{i-1})$. Also, $\rho(B_i)=\rho(F_i)$, since the
relation ``$\prec$'' defines a one-to-one monotone map between $B_i$ and $F_i$ and it preserves the density due to 
Lemma~\ref{lem:no-systematic-slope}.
So, \eqref{eq:density-identity-1} 
and~\eqref{eq:density-identity-2} imply $\rho(C_i)=\rho(D_i)+\rho(E_i)$. Noticing that $\rho(E_i)>0$
does not depend on $i$ and denoting the common value by $\rho_E>0$, we obtain
\begin{equation}
\label{eq:density-identity-3}
\rho(C_i)-\rho(D_i)=\rho_E>0.
\end{equation}

Let us recall that we are proving Lemma~\ref{lem:weak-hyperbolicity}. We made 
an assumption
that there are shocks with absorbed areas that are not too thin, i.e., that satisfy~\eqref{eq:assumption-on-existence-of-fat-shock}. 
Identity~\eqref{eq:density-identity-3} computes the (positive) balance of densities of such shocks emerging and disappearing at time $i$.
Our goal now is to exploit this identity to see that there is a systematic addition of area due to newly emerging shocks 
satisfying~\eqref{eq:assumption-on-existence-of-fat-shock}. 
%To that end, let us first switch to a more abstract combinatorial setting
%and then fill it with concrete content.

%Now let us map the structure of shocks in sets $A_1,\ldots,A_n$ onto the structure described above.
Let $n\in\N$. Let $\Ac^n=\bigcup_{i=1}^n A_i$ and let $\Bc^n$ be the set of all shocks $(i,x)$ such that
$i\in\{1,\ldots,n\}$ and $(i,x)$ is a successor of some shock from $\Ac^n$. Since according to our definition
each shock is its own successor, we have $\Ac^n\subset \Bc^n$. We want to partition $\Bc^n$ into sequences
of successive shocks and estimate the area absorbed by each of those sequences at times $-2,-1,-0$. We must be careful to
avoid overlaps of these areas.

We say that $(i,x)\in\Bc^n$ is the {\it main predecessor} of $(i+1,z)\in\Bc^n$ if 

\begin{enumerate}
\item $(i+1,z)\prec (i,x)$;
\item there is no $y<x$ such that $(i+1,z)\prec (i,y)$ and $(i,y)\in \Ac^n$;
\item if $(i,x)\in\Bc^n\setminus \Ac^n$ then (i)~there is no  $y<x$ such that $(i+1,z)\prec (i,y)$ and $(i,y)\in \Bc^n$;
(ii)~there is no $y$ such that $(i+1,z)\prec (i,y)$ and $(i,y)\in\Ac^n$.
\end{enumerate}
In other words: if at time $i$ there are predecessors of $(i+1,z)$ among $\Ac^n$, we choose the leftmost
of them; if not, we choose the leftmost predecessor  among $\Bc^n$.
Clearly, if $(i+1,z)\in\Bc^n$ has some predecessors in $\Bc^n$ at time $i$, then exactly one of them is the main one.

For any $k,m\in\N$ satisfying $ k \le m \le n$ we define $Q_{k,m}$ as the set of sequences 
$\zeta:\{k,\ldots,m\}\to \R$ satisfying the following conditions: 
\begin{enumerate}[(i)]
\item $(i,\zeta_i)$ is a shock for all $i\in \{k,\ldots,m\}$;
\item  $(k,\zeta_k)\in A_k$;
\item for all $i\in\{1,\ldots,k-1\}$ and all $(i,x)\in A_i$, $(k,\zeta_k)\not\prec (i,x)$; 
\item for each $i\in \{k,\ldots, m-1\}$, $(i,\zeta_i)$ is the main predecessor of $(i+1,\zeta_{i+1})$;
\item if $m<n$, then $(m,\zeta_m)$ is not the main predecessor of its successor.
\end{enumerate}

Sequences in $Q=\bigcup_{1\le k\le m\le n} Q_{k,m}$ viewed as sets of space-time points
generate a partition of $\Bc^n$,  and each shock sequence in $Q$ is uniquely defined by its birth place, i.e., its first entry.
Since the birth places of two distinct sequences cannot be successors of each other, 
the areas absorbed by them are mutually disjoint.

For $k,m\in\N$ such that $k\le m \le n$, we define
\[
R_{k,m}=\bigl\{a:\{k,\ldots,m\}\to \{0,1\}:\ a_k=1\bigr\}
\]
and $a:Q_{k,m}\to R_{k,m}$ by
\[
a_i(\zeta)=\begin{cases} 
                  1,& (i,\zeta_i)\in A_i,\\
                  0,& (i,\zeta_i)\notin A_k.
\end{cases}
\]
For $1\le k\le m\le n$ and $a\in R_{k,m}$, we introduce 
\[
Q_{k,m}^a = \{\zeta\in Q_{k,m}:\ a(\zeta)=a\}
\]
and note that all sets $Q_{k,m}^a$ with all possible values of $k,m$, and $a$ are mutually disjoint, and there are
finitely many of them.

If $1\le k\le m\le n$, $a\in R_{k,m}$, and $i\in \{1,\ldots, n\}$, then  we define 
$Q_{k,m}^a(i)$ to be the section of $Q_{k,m}^a$ at level $i$, i.e., the set consisting of points $x\in\R$ such that $(i,x)$ belongs to a path
from $Q_{k,m}^a(i)$.

\begin{lemma}\label{lem:density-independent-of-the-level-of-section} 
%\begin{enumerate}
%\item For all $i\notin \{k,\ldots,m\}$, $\rho(Q_{k,m}^a(i))=0$.
%\item 
For all $i,j\in\{k,\ldots,m\}$, $\rho(Q_{k,m}^a(i))=\rho(Q_{k,m}^a(j)).$
%\end{enumerate}
\end{lemma}
\bpf This
is a corollary of Lemma~\ref{lem:no-systematic-slope}.
\epf

If $a\in R_{k,m}$ we write $k(a)=k$, $m(a)=m$.
For $a\in S_n=\bigcup_{1\le k\le m\le n} R_{k,m}$, we introduce $\bar a\in\{0,1\}^{\{0,\ldots,n\}}$ by
\[
\bar a_i=
\begin{cases} a_i,& k(a)\le i \le m(a),\\
0,& \text{\rm otherwise.}
\end{cases}
\]
and $\Delta(a)\in\{-1,0,1\}^{\{1,\ldots,n\}}$ by
\[
\Delta_k(a)=\bar a_k-\bar a_{k-1},\quad k=1,\ldots, n,
\]
so $\Delta(a)$ is an alternating sequence of $+1$'s and $-1$'s with possibly some $0$'s between and around them. 
We define also
\[
v(a)=\sum_{k=1}^{n}\frac{1}{k}\Delta_k(a),\quad a\in S_n,
\]
and
\[
w(a)=\frac{1}{k(a)},\quad a\in S_n.
\]
Clearly, $w(a)\ge v(a)>0$ for all $a\in S_n$, due to the alternating character of the sequence $\Delta(a)$.

Let us define $\rho(Q_{k,m}^a)=\rho(Q_{k,m}^a(k))$ and
$\rho(a)=\rho(Q_{k(a),m(a)}^a).$
Then
\begin{align}
\label{eq:total-width-1}
\sum_{a\in S_n}w(a)\rho(a)& \ge  \sum_{a\in S_n}v(a)\rho(a)= \sum_{a\in S_n}\sum_{k=1}^{n}\frac{1}{k}\Delta_k(a) \rho(a)
= \sum_{k=1}^{n}
\frac{c_k-d_k}{k},
%\notag
\end{align}
where for $k=1,\ldots,n$б
\begin{align*}
%\label{eq:transitions-in-and-out-1}
c_k&=\sum_{a\in S_n: \Delta_k(a)=1} \rho(a),\\ 
%\label{eq:transitions-in-and-out-2}
d_k&=\sum_{a\in S_n: \Delta_k(a)=-1} \rho(a).
\end{align*}
Since
\begin{align*}
C_i&=\bigcup_{a\in S_n: \Delta_i(a)=1} Q_{k(a),m(a)}^a(i),\\
D_i&=\bigcup_{a\in S_n: \Delta_i(a)=-1} Q_{k(a),m(a)}^a(i-1),
\end{align*}
we obtain from %~\eqref{eq:transitions-in-and-out-1},\eqref{eq:transitions-in-and-out-2}, and 
Lemma~\ref{lem:density-independent-of-the-level-of-section}:
\begin{align*}
c_i&=\sum_{a\in S_n: \Delta_i(a)=1} \rho(Q_{k(a),m(a)}^a) = \sum_{a\in S_n: \Delta_i(a)=1} \rho(Q_{k(a),m(a)}^a(i)) =\rho(C_i),\\
d_i&=\sum_{a\in S_n: \Delta_i(a)=-1} \rho(Q_{k(a),m(a)}^a)=  \sum_{a\in S_n: \Delta_i(a)=-1} \rho(Q_{k(a),m(a)}^a(i-1))=\rho(D_i).
\end{align*}
So, \eqref{eq:total-width-1} and \eqref{eq:density-identity-3} imply
\begin{equation}
\label{eq:total-w}
\sum_{a\in S_n}w(a)\rho(a)\ge\sum_{i=1}^n\frac{\rho(C_i)-\rho(D_i)}{i}\ge \rho_E  \sum_{i=1}^n \frac{1}{i}.
\end{equation}

\smallskip

By construction and by the definition of sets $A_n$,  to each $a$ and each sequence $\zeta\in Q_{k(a),m(a)}^a$ we can associate a triplet of intervals 
 $J_l(\zeta)\subset\{l\}\times\R$, $l\in\{0,-1,-2\}$ with the following properties:
\begin{equation}
\sum_{l=-2}^0 |J_l(\zeta)|\ge q w(a), \quad \bigcup_{l=-2}^0 J_l(\zeta)\subset \Lambda(k(a),\zeta_{k(a)}).
\label{eq:traces}
\end{equation}
In particular, intervals $J_l(\zeta^1)$ and $J_l(\zeta^2)$ for any $l\in\{0,-1,-2\}$ and any 
two distinct sequences $\zeta^1,\zeta^2$ do not overlap.
Let us denote by
\[
T_n(r)= \sum_{a\in S_n}\ \sum_{\gamma\in Q_{k(a),m(a)}^a:  \zeta_{k(a)}\in[0,r]}\ \sum_{l=-2}^0|J_l(\zeta)|,
\]
the total trace at times $-2,-1,0$ of non-overlapping shocks with birth locations in $\{1,\ldots,n\}\times [0,r]$.
Due to \eqref{eq:total-w} and~\eqref{eq:traces}, 
\[
\liminf_{r\to\infty} \frac{T_n(r)}{r}\ge L_n,
\]
where
\[
L_n=q \rho_E  \sum_{i=1}^n \frac{1}{i}.
\]
Since $\lim_{n\to\infty}L_n=\infty$, we can choose $n$ large enough to ensure $L_n>8$.
Then, for all sufficiently large $r$, we are guaranteed to have $T_n(r)>7r$. Therefore, for all sufficiently large $r$
we will have points in $\{1,\ldots,n\}\times [0,r]$ with associated one-sided minimizers 
containing points on the right of $\{-2,-1, 0\}\times[0,2r]$ (the total length of the segments this set consists of is only $6r$). Applying Lemma~\ref{lem:no-systematic-slope}, we see that this happens with probability~$0$.
Therefore, our assumption that with positive probability there exists a point satisfying~\eqref{eq:assumption-on-existence-of-fat-shock} was wrong,
so the proof of the lemma is completed.\epf

\section{Constructing stationary solutions}\label{sec:global_solutions}

Let us begin with the following auxiliary result:

\begin{lemma}\label{lem:lipschitzness} Let $K>0$. With probability $1$ there is $n_0\in\N$
such that for all $n\ge n_0$, the Lipschitz constant of  $F_\omega(n,x)$ with respect to $x\in[-Kn,Kn]$ is bounded by $\ln n$.
\end{lemma}

\bpf It is sufficient to prove that for sufficiently large $n\in\N$ and any $r\in\N$ with $|r|<Kn+1$, the Lipschitz
constant of  $F_\omega(n,x)$ w.r.t.\ $x\in[r,r+2]$ is bounded by $\ln n$. So let us take arbitrary $n$ and $r$ satisfying the conditions above.
For any $x,y\in[r,r+2]$, the definition
\eqref{eq:shot-noise-definition} and our assumptions on the functions involved in it, imply
\[
|F_\omega(n,x)-F_\omega(n,y)|\le N_\omega(n,r) L|x-y|,
\]
where 
\[
N_\omega(n,r)=\omega(\{n\}\times [r-1,r+3]\times\R\times\R)
\]
and $L$ is the Lipschitz constant of $\phi$.
So the Lipschitz constant of  $F_\omega(n,x)$ w.r.t.\ $x\in[r,r+2]$ is bounded by $N_\omega(n,r) L$.
Since for any $\lambda>0$ there is $C(\lambda)>0$ such that
\[
\Pp\{N_\omega(n,r) L>\ln n\}\le \frac{\E e^{\lambda N_\omega(n,r) L}}{e^{\lambda \ln n}}\le \frac{C(\lambda)}{n^\lambda},
\]
we can choose $\lambda=3$ to see that
\[
\sum_{n\in\N}\ \sum_{r: |r|<Kn+1}\Pp\{N_\omega(n,r) L>\ln n\}<\infty,
\]
so, the Borel--Cantelli Lemma implies that with probability 1 only finitely many events $\{N_\omega(n,r) L>\ln n\}$ happen, and the lemma follows.
\epf

Let us recall that for a fixed $v\in\R$, with probability~$1$, to each point $(n,x)\in\ZR$, we can assign the rightmost
minimizer $\gamma^+(n,x)$ and leftmost minimizer $\gamma^-(n,x)$ with slope $v$. For all but countably many $(n,x)$, 
$\gamma^+(n,x)=\gamma^-(n,x)$.
In this section, we will denote 
$\gamma(n,x)=\gamma^+(n,x)$ for brevity. We also assume without loss of generality that $v=0$ and suppress the dependence on $v$ unless specifically stated otherwise.

\begin{lemma} 
\label{lem:1-st-Busemann}
Let $n\in\Z$ and $x,y\in\R$ with $x<y$. Let $(k_j)_{j\in\N}$ satisfy $k_j\to-\infty$ and 
\begin{equation}
W_{k_j}(n,x,y)<\frac{1}{n-k_j},\quad j\in\N
\label{eq:width-at-tight-places}
\end{equation}
(the existence of such a pairing sequence follows from Theorem~\ref{thm:weak-hyperbolicity}).
Then the limit
\[
B(n,x,y)=\lim_{j\to\infty}\Delta_j
\]
is well-defined and finite, where
\[
\Delta_j=A^{k_j-1,n}(\gamma(n,y))-A^{k_j-1,n}(\gamma(n,x)),\quad j\in\N.
\]
\end{lemma}
\bpf It is sufficient to prove that $\Delta_j$ is a Cauchy sequence.
Let us estimate $|\Delta_j-\Delta_m|$ for $m>j$. To simplify the notation, we denote 
\begin{align*}
\gamma^1&=\gamma(n,x),& \gamma^2&=\gamma(n,y),\\
x_4&=\gamma_{k_m-1}^1,& y_4&=\gamma_{k_m-1}^2,\\
x_3&=\gamma_{k_m-2}^1,& y_3&=\gamma_{k_m-2}^2,\\
x_2&=\gamma_{k_j}^1,& y_2&=\gamma_{k_j}^2,\\
x_1&=\gamma_{k_j-1}^1,& y_1&=\gamma_{k_j-1}^2.
\end{align*}
Since $\gamma^2$ is an optimal path between $(k_j-1, y_1)$ and $(k_m-1, y_4)$, we have
\begin{multline*}
A^{k_j-1,k_m-1}(\gamma^2)\le F_\omega(k_j-1,y_1)+\frac{(x_2-y_1)^2}{2}+ A^{k_j,k_m-2}(\gamma^1)
\\+F_\omega(k_m-2,x_3)+ \frac{(y_4-x_3)^2}{2}.
\end{multline*}
Combining this with
\begin{multline*}
F_\omega(k_j-1,x_1)+\frac{(x_2-x_1)^2}{2}+A^{k_j,k_m-2}(\gamma^1) \\+F_\omega(k_m-2,x_3)+ \frac{(x_4-x_3)^2}{2}= A^{k_j-1,k_m-1}(\gamma^1),
\end{multline*}
we obtain
\[%\begin{equation}
%\label{eq:delta-1}
\Delta_j-\Delta_m\le \delta_{j,m},
\]%\end{equation}
where
\begin{align*}
\delta_{j,m}=& \frac{(x_2-y_1)^2}{2}-\frac{(x_2-x_1)^2}{2}+ \frac{(y_4-x_3)^2}{2} - \frac{(x_4-x_3)^2}{2}\\ 
&+ F_\omega(k_j-1,y_1)-F_\omega(k_j-1,x_1)\\
=&
(x_1-y_1)\left(x_2-\frac{x_1+y_1}{2}\right)+(x_4-y_4)\left(x_3-\frac{x_4+y_4}{2}\right)\\&+ F_\omega(k_j-1,y_1)-F_\omega(k_j-1,x_1).
\end{align*}
Therefore, using~\eqref{eq:width-at-tight-places}, the fact that $|\gamma^1_k|\vee |\gamma^2_k|=o(n-k)$ as $k\to-\infty$, and Lemma~\ref{lem:lipschitzness},
we obtain that there is a function $\beta_1(j)\downarrow 0$  and a number $J_1$ such that if $j>J_1$ and $m>j$, then
$\Delta_j-\Delta_m<\beta_1(j)$.

Interchanging the roles of $\gamma^1$ and $\gamma^2$ we also obtain that there is a function $\beta_2(j)\downarrow 0$  
and a number $J_2$ such that if $j>J_2$ and $m>j$, then
$-(\Delta_j-\Delta_m)<\beta_2(j)$.

Combining these last two statements we conclude that $(\Delta_j)$ is a Cauchy sequence.
\epf

The following lemma is an immediate extension of the previous ones. We can extend definition of $W$ to nonsimultaneous points
$(n_1,x_1), (n_2,x_2)\in\Z\times\R$:
\begin{multline*}
W_k((n_1,x_1), (n_2,x_2))=\sum_{i=0}^{2}\Bigl(\gamma^+_{k-i}(n_1,x_1)\vee\gamma^+_{k-i}(n_2,x_2)\\-\gamma^-_{k-i}(n_2,x_2)\wedge\gamma^-_{k-i}(n_1,x_1) \Bigr),\quad k\le n_1\wedge n_2.
\end{multline*}

\begin{lemma}\label{lem:existence-of-Busemann} Let $(n_1,x_1), (n_2,x_2)\in\Z\times\R$. Then 
\begin{enumerate} 
\item \label{it:minimizers-approach-at-times} There is a sequence $k_j\downarrow-\infty$ such that
\[
W_{k_j}((n_1,x_1), (n_2,x_2))\le \frac{1}{n_1\wedge n_2-k_j},\quad j\in\N.
\]
\item \label{it:exists-Busemann} For every  pairing  sequence $(k_j)$ satisfying conditions of part~\ref{it:minimizers-approach-at-times}, the following
finite limit exists:
\[
B((n_1,x_1), (n_2,x_2))=\lim_{j\to\infty} (A^{k_j-1,n_2}(\gamma(n_2,x_2))-A^{k_j-1,n_1}(\gamma(n_1,x_1))).
\]
\item The limit in part~\ref{it:exists-Busemann} does not depend on the concrete choice of $(k_j)$.
\end{enumerate}
\end{lemma}
\bpf The first part of the lemma follows since we can apply Theorem~\ref{thm:weak-hyperbolicity} to points 
$(n_1\wedge n_2,\gamma_{n_1\wedge n_2}(n_1,x_1))$ and $(n_1\wedge n_2,\gamma_{n_1\wedge n_2}(n_2,x_2))$ that share the
time coordinate. The second part holds since
\begin{multline*}
A^{k_j-1,n_1}(\gamma(n_1,x_1))=A^{k_j-1,n_1\wedge n_2}(\gamma(n_1\wedge n_2, \gamma_{n_1\wedge n_2}(n_1,x_1)))\\ +A^{n_1\wedge n_2,n_1}(\gamma(n_1,x_1))
\end{multline*}
and
\begin{multline*}
A^{k_j-1,n_2}(\gamma(n_2,x_2))=A^{k_j-1,n_1\wedge n_2}(\gamma(n_1\wedge n_2, \gamma_{n_1\wedge n_2}(n_2,x_2)))\\+A^{n_1\wedge n_2,n_1}(\gamma(n_2,x_2)),
\end{multline*}
so one can apply Lemma~\ref{lem:1-st-Busemann} to points $(n_1\wedge n_2,\gamma_{n_1\wedge n_2}(n_1,x_1))$ and $(n_1\wedge n_2,\gamma_{n_1\wedge n_2}(n_2,x_2))$.
The last part follows from the standard trick of interlacing the two sequences.
\epf

The function $B((n_1,x_1), (n_2,x_2))=B_\omega((n_1,x_1), (n_2,x_2))$ may be called the Busemann function 
in analogy to the Busemann functions used in~\cite{BCK:MR3110798} and previous
work on last-passage percolation, although we stress that in our setting we are currently 
able to prove
convergence only along appropriate subsequences $k_j$.
This function has several standard properties.
Some of them are summarized in the following lemma:
\begin{lemma} \label{lem:busemann_properties} Let $B$ be defined as above.
 \begin{enumerate}
  \item The distribution of $B$ is translation invariant: for any $\Delta \in\Z\times \R$,
\[
 B(\cdot+\Delta,\cdot+\Delta)\stackrel{distr}{=} B(\cdot,\cdot).
\]
\item\label{it:anti-symmetry-of_Busemann}
 $B$ is antisymmetric:
\[
 B((n_1,x_1),(n_2,x_2))=-B((n_2,x_2), (n_1,x_1)),\quad (n_1,x_1),(n_2,x_2)\in\Z\times\R,
\]
in particular $B(n,x),(n,x))=0$ for any $(n,x)\in\Z\times\R$.
\item $B$ is additive: for any $(n_1,x_1),(n_2,x_2),(n_3,x_3)\in\Z\times\R$,
\[
 B((n_1,x_1),(n_3,x_3))=B((n_1,x_1),(n_2,x_2))+B((n_2,x_2),(n_3,x_3)). 
\]

\item\label{it:upper-estimate-on-Busemann} For any $(n_1,x_1),(n_2,x_2)\in\Z\times\R$ satisfying $n_1<n_2$,
\begin{equation}
 B((n_1,x_1),(n_2,x_2))\le A^{n_1,n_2}(x_1,x_2).
\end{equation}

\item  For any $(n_1,x_1),(n_2,x_2)\in\Z\times\R$, $\E |B(n_1,x_1),(n_2,x_2)|<\infty$. \label{item:expectation_finite}
 \end{enumerate}
\end{lemma}
\bpf The first two parts are obvious. For the third part we need to ensure that the convergence in the definition of $B$
for all the values of arguments involved holds along the same sequence $k_j$. That is true, since all one-sided minimizers are ordered
and the gap between the two minimizers on the sides dominates the smaller gaps between the middle one and the side ones.

Let us prove part~\ref{it:upper-estimate-on-Busemann}. We denote $\gamma^1=\gamma(n_1,x_1)$ and $\gamma^2=\gamma(n_2,x_2)$.
Let us find a pairing sequence $(k_j)_{j\in\N}$ for $\gamma^1,\gamma^2$. For any $j$ let us create a path $\gamma(j)$ that starts
at $(k_j-1,\gamma^2_{k_j-1})$, makes a step to $(k_j,\gamma^1_{k_j})$, coincides with $\gamma^1$ between $k_j$ and $n_1$, and 
on $\{n_1,\ldots,n_2\}$ coincides with the optimal path between $(n_1,x_1)$ and $(n_2,x_2)$. Then
\begin{align*}
&A^{k_j-1,n_2}(\gamma^2)\le  A^{k_j-1,n_2}(\gamma(j))
\\\le& F(k_j-1,\gamma^2_{k_j-1})+\frac{(\gamma^1_{k_j}-\gamma^2_{k_j-1})^2}{2}+A^{k_j,n_1}(\gamma^1)+A^{n_1,n_2}(x_1,x_2)
\\ \le& F(k_j-1,\gamma^2_{k_j-1})-F(k_j-1,\gamma^1_{k_j-1})+\frac{(\gamma^1_{k_j}-\gamma^2_{k_j-1})^2-(\gamma^1_{k_j}-\gamma^1_{k_j-1})^2 }{2}
\\ &+A^{k_{j-1},n_1}(\gamma^1)+A^{n_1,n_2}(x_1,x_2).
\end{align*}
Moving $A^{k_{j-1},n_1}(\gamma^1)$ to the left-hand side, taking limit $j\to\infty$, using the pairing property of $k_j$, 
the sublinear growth of
$|\gamma^1|,|\gamma^2|$, and Lemma~\ref{lem:lipschitzness}, we finish the proof of part~\ref{it:upper-estimate-on-Busemann}.

Let us prove the last part.
Using additivity and translation invariance we see that it is sufficient to consider points $(0,0)$ and $(n,x)$ with $n<0$.
Parts~\ref{it:anti-symmetry-of_Busemann} and~\ref{it:upper-estimate-on-Busemann} of this lemma along with Lemma~\ref{lem:linbound} imply that
\[
\E B((0,0),(n,x))\ge - \E A^{n,0}(x,0) > -\infty.
\]
So it remains to prove an upper bound. Furthermore, 
\[
B((0,0),(n,x))=\sum_{k=1}^{|n|} B\left(\left(-(k-1),x\frac{k-1}{|n|}\right), \left(-k,x\frac{k}{|n|}\right)\right),
\]
and all the terms on the right-hand side have the same distribution. Expectation of each of them belongs to $(-\infty,\infty]$
Therefore, $\E B((0,0)(n,x))$ is finite if and only if $\E B((0,0)(-1,x/|n|))$ is finite. So it is sufficient to prove
$\E B((0,0),(-1,x))<\infty$ for all $x\in\R$. Furthermore, we have
\[
B((0,0),(-2,0))= B((0,0),(-1,x))+ B((-1,x),(-2,0)),
\]
where, due to the symmetry of the Poissonian process and the action, the distributions of two terms on the right-hand side coincide.
Therefore, $\E B((0,0),(-1,x))<\infty$ iff $\E B((0,0),(-2,0))<\infty$. Applying this once again, we see that it is sufficient
to establish $\E B((0,0),(-1,0))<\infty$.

We denote $\gamma^0=\gamma(0,0)$ and $\gamma^1=\gamma(0,-1)$ for brevity.

Let $L_1>0$ (we will later impose some conditions on $L_1$). $H=\{(n,x)\in\Z\times\R:\ |x|\le L_1(-n-1)\}$. Since $\gamma^0$ has asymptotic slope $0$, the time
$\tau$ defined by
\[
\tau=\min\{n\le 0: \gamma^0_k\notin H\}-1
\]
is finite. We define $z=\gamma^0(\tau)$
Let $(k_j)$ be a pairing sequence for $\gamma^0$ and $\gamma^1$. For $k_j<\tau$, we have
\[ %begin{align*}
A^{k_j-1,-1}(\gamma^1)-A^{k_j-1,0}(\gamma^0)\le
Q_{k_j}(\gamma^0,\gamma^1)+ A^{\tau,-1}(z,0)- A^{\tau,0}(\gamma^0),
\]%\end{align*}
where
\begin{align*}
Q_{n}(\gamma^0,\gamma^1)=&F(n-1,\gamma^1_{n-1})+\frac{(\gamma^0_n-\gamma^1_{n-1})^2}{2}\\&-F(n-1,\gamma^0_{n-1})-\frac{(\gamma^0_n-\gamma^0_{n-1})^2}{2}
\\ =&F(n-1,\gamma^1_{n-1})-F(n-1,\gamma^0_{n-1})\\&+(\gamma^1_{n-1}-\gamma^0_{n-1})\left(\frac{\gamma^0_{n-1}+\gamma^1_{n-1}}{2}-\gamma^0_n\right)
\end{align*}
Since $Q_{k_j}(\gamma^0,\gamma^1)\to 0$ as $j\to\infty$, we obtain
\[
B((0,0),(-1,0))\le A^{\tau,-1}(z,0)- A^{\tau,0}(\gamma^0)=A^{\tau,-1}(z,0)- A^{\tau,0}(z,0),
\]
and it is sufficient to prove
\begin{equation}
\label{eq:finite-expectation-1}
\E A^{\tau,-1}(z,0)<\infty
\end{equation}
and
\begin{equation}
\label{eq:finite-expectation-2}
\E A^{\tau,0}(z,0)>-\infty.
\end{equation}

Let us estimate the tail of the distribution of $\tau$. First, we choose $L_2>0$ so that $L_1-L_2>2R$,
where $R$ is chosen according to Lemma~\ref{lem:distribution-of-animal-size}.
If $\tau=n$,  $\gamma^0_{\tau+1}>0$, and $\gamma^0_{\tau}-\gamma^0_{\tau+1}< -(L_1-L_2)|n|$, then the
optimal path $\gamma'$ connecting $(n, L_1(|n|-2)-(L_1-L_2)|n|)=(n, L_2|n|-2L_1)$ to $(0,0)$ satisfies $\gamma'_{n+1}>L_1(|n|-2)$,
and thus, if $|n|$ is sufficiently large, say, greater than some $n_1$, it deviates from the straight line connecting $(n, L_2|n|-2L_1)$ to $(0,0)$ by at least 
$L_1(|n|-2)-(L_2|n|-2L_1)(|n|-1)/|n|>Rn$. Applying Lemma~\ref{lem:path_in_wide_rectangle_whp},
we obtain that
\begin{equation*}
\Pp\bigl\{\tau=n, \gamma^0_{n+1}>0, \gamma^0_{n}-\gamma^0_{n+1}< -(L_1-L_2)|n|\bigr\}\le c_1\exp(-c_2 n),\quad n\ge n_1,
\end{equation*}
for some $c_1,c_2>0,$
and, similarly,
\begin{equation*}
\Pp\bigl\{\tau=n, \gamma^0_{n+1}<0, \gamma^0_{n}-\gamma^0_{n+1}> (L_1-L_2)|n|\bigr\}\le c_1\exp(-c_2 n),\quad n\ge n_1.
\end{equation*}
Combining these two inequalities we obtain
\begin{equation}
\label{eq:time-to-enter-cone-1}
\Pp\bigl\{\tau=n, \gamma^0_{n}\notin I_n \cup (-I_n)\bigr\}\le 2c_1\exp(-c_2 n),\quad n\ge n_1,
\end{equation}
where $I_n=[L_2|n|-2L_1, L_1(|n|-1)]$.

Let us now take any $n\in-\N$, any point $y\in I_n\cup (-I_n)$ and suppose that $\tau=n$ and $|\gamma^0_n-y|\le 1$.
Since the asymptotic slope of $\gamma^0$ equals $0$, for sufficiently large values of $|m|$ we will have
$|\gamma^0_m|<L_2|m|/2$. Lemma~\ref{lem:quantitative-straightness} implies that probability of such
an event is bounded by $C_1e^{-C_2|n|^\kappa}$ for sufficiently large $|n|$, and since one can find $3(L_1-L_2)|n|$ points $y$ such that
the union of segments $[y-1,y+1]$ covers the entire set $I_n\cup(-I_n)$, we obtain that there is $n_2>0$ such that
\begin{equation}
\Pp\bigl\{\tau=n, \gamma^0_{n}\in I_n \cup (-I_n)\bigr\}\le 2(L_1-L_2)|n| C_1e^{-C_2|n|^\kappa},\quad |n|\ge n_2.
\label{eq:time-to-enter-cone-2}
\end{equation}
Combining~\eqref{eq:time-to-enter-cone-1} and~\eqref{eq:time-to-enter-cone-2}, we obtain that there are constants $\bar C_1, \bar C_2>0$
such that
\begin{equation}
\Pp\bigl\{\tau=n\}\le \bar C_1e^{-\bar C_2|n|^\kappa},\quad n\in\N.
\label{eq:time-to-enter-cone-3}
\end{equation}

\smallskip

Let us now prove~\eqref{eq:finite-expectation-1}. Notice that on $\{\tau=n\}$,
\begin{equation}
A^{\tau,-1}(z,0)\le F_\omega(n,z)+\frac{(L_1n)^2}{2}+A^{n+1,-1}(0,0),
\label{eq:upper-on-A}
\end{equation}
and
\[
A^{n-1,-1}(0,0)\le F_\omega(n-1,0)+ \frac{(L_1n)^2}{2}+A^{\tau,-1}(z,0),
\]
so
\begin{equation}
A^{\tau,-1}(z,0)\ge -F_\omega(n-1,0)-\frac{(L_1n)^2}{2}+A^{n-1,-1}(0,0).
\label{eq:lower-on-A}
\end{equation}
Combining~\eqref{eq:upper-on-A} and~\eqref{eq:lower-on-A}, we obtain that on $\{\tau=n\}$,
\begin{multline*}
|A^{\tau,-1}(z,0)|\le F_\omega^*(\{n\}\times[-L_1|n|,L_1|n|])+|F_\omega(n-1,0)|+(L_1n)^2
\\+|A^{n+1,-1}(0,0)|+|A^{n-1,-1}(0,0)|,
\end{multline*}
where $F^*_\omega$ was defined in~\eqref{eq:potential-maximum}.
The Cauchy--Schwarz inequality implies
\begin{align*}
\E |A^{\tau,-1}(z,0)|\ONE_{\{\tau=n\}}\le & \sqrt{\Pp\{\tau=n\}}
\Bigl(\E {F_\omega^*}^2(\{n\}\times[-L_1|n|,L_1|n|])
\\&+\E |F_\omega(n-1,0)|
+(L_1n)^2\\&+\E (A^{n+1,-1}(0,0))^2+\E (A^{n-1,-1}(0,0))^2\Bigr).
\end{align*}
The first term in the parentheses grows logarithmically in $|n|$, the second term does not depend on $n$, and the remaining 
terms grow quadratically
due to Lemma~\ref{lem:moments_of_At}. Combining this with~\eqref{eq:time-to-enter-cone-3}, we finish the proof 
of~\eqref{eq:finite-expectation-1}. Inequality~\eqref{eq:finite-expectation-2} can be proved in exactly the same way. \epf

Let us now define 
\[
U_\omega(n,x)=U(n,x)=B((0,0),(n,x)),\quad (n,x)\in\Z\times\R.
\]
The main claim of this section is that thus defined $U$ is skew invariant under of the
HJBHLO cocycle, and its space derivative is a global solution
of the Burgers equation.

Let us recall that the HJBHLO evolution is given by
\begin{equation}
\label{eq:Burgers_dynamics_on_potentials}
 \Phi^{m,n}W(y)=\inf_{x\in\R} \{W(x)+A^{m,n}(x,y)\}, \quad m\le n,\quad y\in\R,
\end{equation}
where $A^{m,n}(x,y)$ has been defined in~\eqref{eq:optimal_action_between_two_points}.

\begin{theorem}\label{thm:construction-of-solution} The random function $U$ is a global solution of the HJBHLO, i.e., for
almost all $\omega\in\Omega$,
\[
\Phi_{\omega}^{n_1n_2}U_\omega(n_1,\cdot)= U_\omega(n_2,\cdot),\quad n_1<n_2.
\]

\end{theorem}

\bpf Let $\gamma$ be a minimizer with endpoint $(n_2,x)$. Then
\begin{align*}
 U(n_2,x)=&U(n_1,\gamma_{n_1})+(U(n_2,x)-U(n_1,\gamma_{n_1}))\\
                      =&U(n_1,\gamma_{n_1}) + A^{n_1,n_2}(\gamma_{n_1},x).
\end{align*}
We need to show that the right-hand side is the infimum of $U(n_1,y)+A^{n_1,n_2}(y,x)$ over
all $y\in\R$. Suppose that for some $y\in\R$,
\begin{equation}
U(n_1,y)+A^{n_1,n_2}(y,x)  < U(n_1,\gamma_{n_1}) + A^{n_1,n_2}(\gamma_{n_1},x).
\label{eq:suppose_not_infimum}
\end{equation}

Let us take any minimizer $\bar \gamma$ originating at $(n_1,y)$ and 
use Lemma~\ref{lem:existence-of-Busemann} to
find a pairing
sequence $(k_j)$ for $(n_1,y)$ and $(n_1,\gamma_{n_1})$. Then
\eqref{eq:suppose_not_infimum} implies
\begin{align*}
\lim_{j\to\infty} (A^{k_j-1,n_1}(\bar \gamma)-A^{k_j-1,n_1}(\gamma))&=U(n_1,y) - U(n_1, \gamma_{n_1})\\
&< A^{n_1,n_2}(\gamma_{n_1},x) -A^{n_1,n_2}(y,x).
\end{align*}
Denoting the right-hand side by $\delta$, we conclude that for sufficiently large $j$,
\begin{align*}
A^{k_j-1,n_1}(\bar \gamma)+A^{n_1,n_2}(y,x)&\le A^{k_j-1,n_1}(\gamma)+A^{n_1,n_2}(\gamma_{n_1},x)-\delta/2\\
&\le A^{k_j-1,n_2}(\gamma) - \delta/2.
\end{align*}
Let now $\gamma'(j)$ be the path starting at $(k_j-1,\gamma_{k_j-1})$, making an immediate step to $(k_j,\bar\gamma_{k_j})$, 
coinciding with $\bar\gamma$ between times $k_j$ and $n_1$, and coinciding with the optimal path connecting $y$ to $x$ between
$n_1$ and $n_2$. We have
\begin{multline*}
A^{k_j-1,n_2}(\gamma'(j))\le A^{k_j-1,n_1}(\bar \gamma)+A^{n_1,n_2}(y,x) \\ + \frac{(\gamma_{k_j-1}-\bar \gamma_{k_j})^2-(\bar\gamma_{k_j-1}-\bar \gamma_{k_j})^2}{2}\\
+F_\omega(k_j-1,\gamma_{k_j-1})-F_\omega(k_j-1,\bar \gamma_{k_j-1}).
\end{multline*}
Combining the last two inequalities,
we obtain
\begin{align*}
A^{k_j-1,n_2}(\gamma'(j))\le A^{k_j-1,n_2}(\gamma)+r(j),
\end{align*}
where
\begin{multline*}
r(j)=-\frac{\delta}{2}+\left(\frac{\gamma_{k_j-1}+\bar \gamma_{k_j-1}}{2}-\bar \gamma_{k_j}\right) (\gamma_{k_j-1}-\bar \gamma_{k_j-1})
\\+
F_\omega(k_j-1,\gamma_{k_j-1})-F_\omega(k_j-1,\bar \gamma_{k_j-1}).
\end{multline*}
For sufficiently large $j$, $r(j)<0$. This is implied by sublinear growth of
 $\gamma_{k_j-1},\bar \gamma_{k_j-1},\bar \gamma_{k_j}$ in $k_j$, the rate of decay of 
$|\gamma_{k_j-1}-\bar \gamma_{k_j-1}|$ guaranteed by Lemma~\ref{lem:existence-of-Busemann}, and Lemma~\ref{lem:lipschitzness}.

So, $A^{k_j-1,n_2}(\gamma'(j))< A^{k_j-1,n_2}(\gamma)$ which contradicts the fact that $\gamma$ is an optimal path between
points $((k_j-1),\gamma_{k_j-1})$ and $(n_2,x)$. Therefore, our assumption
on existence of $y$ satisfying~\eqref{eq:suppose_not_infimum} was wrong, and the proof is complete.
\epf

To prove $U(n,\cdot)\in \HH$ for all $n\in\Z$, we begin with the Lipschitz property.
\begin{lemma}
For all $n\in\Z$, $U_v(n,\cdot)$ is locally Lipschitz.
\end{lemma}
\bpf Let us fix $n\in\N$ and any points $z_1,z_2\in \R$ satisfying $z_1<z_2$. Let $\gamma^-=\gamma^-(n,z_1)$ and
$\gamma^+=\gamma^+(n,z_2)$. Then $\gamma^-_k<\gamma^+_k$ for all $k\le n$ and there is a pairing sequence $k_j\downarrow -\infty$ such that
\[
|\gamma^+_{k_j}-\gamma^-_{k_j}|+|\gamma^+_{k_j-1}-\gamma^-_{k_j-1}|<\frac{1}{n-k_j},\quad j\in\N.
\] 
Let us now take any points $x,y\in (z_1,z_2)$ and denote $\gamma^1=\gamma(n,x)$, $\gamma^2=\gamma(n,y)$. These minimizers
 pass between
$\gamma^-$ and $\gamma^+$, and, in particular, the same sequence $(k_j)$ is pairing for $\gamma^1$ and $\gamma^2$, so 
\[
|\gamma^2_{k_j}-\gamma^1_{k_j}|+|\gamma^2_{k_j-1}-\gamma^1_{k_j-1}|<\frac{1}{n-k_j},\quad j\in\N.
\]
As in the proof of Lemma~\ref{lem:existence-of-Busemann}, we obtain
\begin{align*}
A^{k_j-1,n}(\gamma^2)-A^{k_j-1,n}(\gamma^1)\le &
(\gamma_{k_j-1}^1-\gamma_{k_j-1}^2)\left(\gamma_{k_j}^1-\frac{\gamma_{k_j-1}^1+\gamma_{k_j-1}^2}{2}\right)
\\&+(x-y)\left(\gamma_{n-1}^1-\frac{x+y}{2}\right)
\\&+ F_\omega(k_j-1,\gamma_{k_j-1}^2)-F_\omega(k_j-1,\gamma_{k_j-1}^1).
\end{align*}
Taking $j\to\infty$, we obtain
\[
U_n(y)-U_n(x)\le (x-y)\left(\gamma_{n-1}^1-\frac{x+y}{2}\right).
\]
Since $\gamma_{n-1}^1\in (\gamma_{n-1}^-,\gamma_{n-1}^+)$  irrespective of the choice of $x,y\in(z_1,z_2)$,
we conclude that there is $C_1(z_1,z_2)$ such that
\[
U(n,y)-U(n,x)\le C_1(z_1,z_2) |x-y|,\quad  z_1<x<y<z_2. 
\]
Similarly, for some  $C_2(z_1,z_2)$,  we  obtain 
\[
U(n,x)-U(n,y)\le C_2(z_1,z_2) |x-y|,\quad  z_1<x<y<z_2, 
\]
which completes the proof.
\epf

Although we have always assumed that $v=0$ in this section, all the definitions, constructions and results hold true for other values 
$v$ as well, due to the Galilean shear invariance. Let us denote the corresponding Busemann function and global HJBHLO solution
by $B_v$ and $U_v$.

To prove that $U_v(n,\cdot)\in \HH(v,v)$ for all $n$, we will compute the expectation of its spatial increments
(we already know that it is well defined due to part~\ref{item:expectation_finite} of
Lemma~\ref{lem:busemann_properties}), and prove that $u_v(n,\cdot)$ is ergodic with respect to the spatial variable.

\begin{lemma}\label{lem:mean_increment} For any  $(n,x)\in\Z\times\R$,
\[
 \E (U_v(n,x+1)-U_v(n,x))=\E B_v((n,x),(n,x+1))=v.
\]
\end{lemma}
\bpf
First, we consider the case $v=0$. Due to the distributional invariance of the potential process $F$ under reflections,
\[
\E B_0((n,x),(n,x+1)) = \E B_0((n,x+1),(n,x)).
\]
Combining this with the anti-symmetry of $B_0$, we obtain $\E B_0((n,x+1),(n,x))=0$, as required.

In the general case, we can apply the shear transformation $L$ of $\Z\times\R\times\R\times\R$ defined by
\[
L:(m,y,q,r)\mapsto (m,y+(n-m)v,q,r).
\]
Due to Lemma
\ref{lem:shear}, the one-sided minimizers with slope $v$ will be mapped onto
one-sided minimizers of slope $0$ for the new potential $F_{L(\omega)}$.
We already know that
\[
 \E B_{0}((n,x+1),(n,x))_{L(\omega)}=0.
\]
A direct computation based on Lemma~\ref{lem:shear} gives
\[
  B_{0,L( \omega)}((n,x),(n,x+1))=B_{v,\omega}((n,x),(n,x+1))-v,
\]
and our statement follows since $L$ preserves the driving measure $\mu$ and hence the distribution
of the Poisson process and the potential $F$.
\epf

So far we have worked with solutions of the Hamilton--Jacobi equation. One can obtain the corresponding
solutions of the Burgers equation by
\[
 u_v(n,x)=\gamma_{v,n}(n,x)-\gamma_{v,n}(n-1,x).
\]
Then $U_v(n,x)-U_v(n,0)=\int_0^x u_v(n,y)dy$.
We recall that $\Psi^{m,n}w$ denotes the solution at time $n$ of the Burgers equation with initial condition $w$
imposed at time~$m$.
\begin{lemma} The function $u_v$ defined above is a global solution of the kicked Burgers equation. If $m\le n$, then
\[
 \Psi^{m,n} u_v (m,\cdot) = u_v(n,\cdot),\quad m\le n.
\]
\end{lemma}
\bpf This statement is a direct consequence of 
Lemma~\ref{lem:properties-of-evolution-operator},
Theorem~\ref{thm:construction-of-solution}, and
the definition of the Burgers cocycle $\Psi$.\epf

\begin{lemma}\label{eq:spatial-ergodicity} For any $v$ and any $n$, the process $u_v(n,\cdot)$ is stationary and ergodic with respect to
spatial translations.
\end{lemma}
\bpf Denoting $u_v(n,0)(\omega)=\xi(\omega)$, we see that due to space-time invariance of the procedure of
constructing of one-sided minimizers, $u_v(n,x)=\xi(\tau_{x}\omega)$, where $\tau_x$ denotes
the space shift of the Poisson process by distance~$x$. Since the measure $\Pp$ is invariant and ergodic with respect to spatial translations,
we conclude that so is  $u_v(n,\cdot)$.\epf

\begin{theorem}
For any $v\in\R$ and any $n\in\N$, we have $U_{v,\omega}(n,\cdot)\in \HH(v,v)$. The sequence $u_{v,\omega}(n,\cdot), n\in\N$ is 
a stationary process with values in $\GG(v,v)$.
\end{theorem}
\bpf The first claim is a direct consequence of Lemmas~\ref{lem:mean_increment} and \ref{eq:spatial-ergodicity},
and Birkhoff's ergodic theorem. The second claim follows from the first one and the space-time invariance of
the construction of minimizers.
\epf

\section{Stationary solutions: uniqueness and basins of attraction}\label{sec:attractor}
In this section we prove Theorem~\ref{thm:pullback_attraction} and the uniqueness part in
Theorem~\ref{thm:global_solutions}.
The key step is the following observation.

\begin{lemma}\label{lem:asymptotic_slope_in_pullback_attraction}
Let $n\in\Z$ and suppose that an initial condition $W$ satisfies one of the
conditions~\eqref{eq:no_flux_from_infinity},\eqref{eq:flux_from_the_left_wins},\eqref{eq:flux_from_the_right_wins}.
 With probability one, the following holds true for every $y\in\R$. Let  $y^*(m)$ be a
solution of the optimization problem
~\eqref{eq:Burgers_dynamics_on_potentials} for each $m\le n$. Then
\[
 \lim_{m\to-\infty}\frac{y^*(m)}{m}=v.
\]
\end{lemma}
\bpf[Proof of Theorem~\ref{thm:pullback_attraction}]  Since $u_v$ is a solution
of the Burgers equation over any finite time interval, $u_v(n,\cdot)$ is continuous at $y$ iff there is a unique one-sided
minimizer $\gamma(n,y)$ with endpoint $(n,y)$. Moreover, in this case, $u_v(n,y)=\gamma_n(n,y)-\gamma_{n-1}(n,y)$.

So let us take such a point $y$. 
Lemma~\ref{lem:asymptotic_slope_in_pullback_attraction} and
Lemma~\ref{lem:limits-of-minimizers-uniquely-defined} guarantee then that solutions $y^*(m)$ for optimization problem
~\eqref{eq:Burgers_dynamics_on_potentials} and the corresponding optimal paths $\gamma^{(m,y^*(m)),(n,y)}$ realizing
$A^{m,n}(y^*(m),y)$ converge pointwise to the infinite one-sided minimizer $\gamma(n,x)$.
\marginpar{$(n,y)$?}
 In particular, 
\begin{multline*}
\Psi^{m,n}w(y)=\gamma_n^{(m,y^*(m)),(n,y)}-\gamma_{n-1}^{(m,y^*(m)),(n,y)}
\\ \to  \gamma_n(n,y)-\gamma_{n-1}(n,y)=u_v(n,y),\quad m\to-\infty,
\end{multline*}
which completes the proof.
\epf

\bpf[Proof of Lemma~\ref{lem:asymptotic_slope_in_pullback_attraction}] We will only prove the sufficiency of
condition~\eqref{eq:no_flux_from_infinity}. The proof of sufficiency of conditions
\eqref{eq:flux_from_the_left_wins} and \eqref{eq:flux_from_the_right_wins} follows the same lines and we omit it.

Let us also restrict ourselves to $n=0$ for simplicity. The proof does not change for other values of $n$.

 Since $y^*$ is increasing in $y$, it is sufficient to show that the conclusion of the lemma holds with
probability 1 for fixed $y$. The stationarity of the forcing potential implies that we can assume $y=0$.

We must show that for any $\eps>0$ it is extremely unlikely for a path $\gamma$ with $\gamma_0=0$
and $|\gamma_m|> \eps |m|$ to provide a solution to~\eqref{eq:Burgers_dynamics_on_potentials} if $|m|$ is large.
For definiteness, let us work with paths satisfying
$0^*(m)=\gamma_m>\eps |m|$.

For any $\delta>0$ and for sufficiently large $|m|$,
\[
 W(0) + A^{m,0}(0,0)< (\alpha(0)+\delta)|m|.
\]
If additionally $x=0^*(m)>\eps |m|$, then for $i=[x-\eps|m|]$, 
\[
 \inf_{z\in[\eps |m| +i,\eps |m| +i+1]} W(z)+ A^{m-1,0}(\eps |m|+i,0)< (\alpha(0)+\delta)|m|+\frac{1}{2}+F(m-1,\eps |m|+i).
\]
Condition~\eqref{eq:no_flux_from_infinity} at $+\infty$ implies that there is $r_1$ such that for $m<-r_1$ and all
$i\in\N$,
\[
 \inf_{z\in[\eps |m| +i,\eps |m| +i+1]} W(z)> -(|m|+i) \delta,
\]
so there is $r_2>0$ such that if $m<-r_2$ and $i\in\N$, then
\begin{align*}
 A^{m-1,0}(\eps |m|+i,0)-F(m-1,\eps |m|+i) &< (\alpha(0) + 2\delta)|m| +\delta i + \frac{1}{2} 
 \\&<
|m|\left(\alpha(0)+3\delta+\delta\frac{i}{|m|}\right).
\end{align*}

Let us denote by $B_{mi}$ the event defined by this inequality. Then $B_{mi}\subset C_{mi}\cup D_{mi}$,
where
\begin{align*}
C_{mi}&=\bigl\{ -F(m-1,\eps |m|+i)< - \delta(|m|+i)\bigr\}\\
D_{mi}&=\left\{ A^{m-1,0}(\eps |m|+i,0)< |m|\left(\alpha(0)+4\delta+2\delta\frac{i}{|m|}\right)\right\}
\end{align*}

Due to the Borel--Cantelli lemma, to show that with probability 1, events $B_{mi}$ can happen only for finitely many values
of $m$, it suffices to
show that for some $\beta>0$ and $c>0$,
\begin{equation}
\label{eq:double_sum-for-tails}
\sum_{m\le - c}\,\sum_{i}\Pp(C_{mi})<\infty,
\end{equation}
\begin{equation}
\label{eq:at_most_linear_growth_BC}
\sum_{m\le - c}\,\sum_{i\le \beta |m|}\Pp(D_{mi})<\infty,
\end{equation}
and
\begin{equation}
\label{eq:super_linear_growth_BC}
\sum_{m\le - c}\,\sum_{i>\beta |m|}\Pp(D_{mi})<\infty.
\end{equation}
Inequality~\eqref{eq:double_sum-for-tails} does not depend on $\beta$ and holds for any $c$ since
\[
\Pp(C_{mi})\le e^{-\delta(|m|+i)}\E e^{F(0,0)}.
\]

Denoting $\alpha_{mi}=\alpha\left(\frac{\eps |m|+i}{|m|+1}\right)$, using shear and translation invariance,
we
obtain
\begin{multline*}
\Pp(D_{mi})
 =
\Pp\Bigl\{A^{|m|+1}-\alpha(0)(|m|+1) \\ <
|m|\left(\alpha(0)+4\delta+\frac{2\delta i}{|m|}- \frac{|m|+1}{|m|}\alpha_{mi}\right)\Bigr\}.
\end{multline*}
If $\delta$ is sufficiently small, then, using Lemma~\ref{lem:shape-function},  we can find $r_2$ such
that for all $m<-r_2$ and all $i$,
\[
|m|\left(\alpha(0)+4\delta+\frac{2\delta i}{|m|}- \frac{|m|+1}{|m|}\alpha_{mi}\right)
 <
-(|m|+1)\left(\frac{\eps^2}{2}+\frac{i^2}{2|m|^2}\right),
\]
so
\begin{equation}
\label{eq:D-mi}
\Pp(D_{mi})\le \Pp\left\{A^{|m|+1}-\alpha(0)(|m|+1) < -(|m|+1)\left(\frac{\eps^2}{2}+\frac{i^2}{2|m|^2}\right)\right\}.
\end{equation}

Now~\eqref{eq:at_most_linear_growth_BC} follows (with any $c\ge r_2$ and with arbitrary choice of $\beta$) from
Theorem~\ref{thm:concentration_around_alphat}.

To prove~\eqref{eq:super_linear_growth_BC} we need an auxiliary lemma. In its statement and proof we use the notation
introduced in Section~\ref{sec:subadd}.
\begin{lemma}
\label{lem:tails_of_action}
There are constants $c_1,c_2,X_0,N_0>0$ such that for $n>N_0$, $x>X_0$,
 \[
  \Pp\{A^n\le  - xn\}\le c_1 e^{-c_2 xn}.
 \]
\end{lemma}
Now~\eqref{eq:super_linear_growth_BC} is a consequence of this lemma and~\eqref{eq:D-mi}
 if we choose $c>N_0$
and
$\beta$ satisfying
\[
\frac{\eps^2}{2}+\frac{\beta^2}{2}-\alpha(0)>X_0.
\]
It remains to prove Lemma~\ref{lem:tails_of_action}.

\bpf[Proof of Lemma~\ref{lem:tails_of_action}] 
 Let us take $c_3>0$ and write
\begin{equation}
 \Pp\{A^n\le -xn\}\le \Pp\{\Sigma \le c_3 xn,\ A^n\le -xn\}+\Pp\{\Sigma > c_3 xn\},
\label{eq:estimating_action_linear_tails}
\end{equation}
where $\Sigma=\Sigma(\gamma)$ has been defined in~\eqref{eq:Sigma} for a path $\gamma$ realizing~$A^n$.

To estimate the first term on the right-hand side,  we choose $c_3>0$ and  $x_0>0$  so that
for all $x>x_0$ and all $n\in\N$, $\lceil c_3 xn\rceil>n$ and $xn>y_0\lceil c_3 xn\rceil$,
 where $y_0$ was introduced before
inequality~\eqref{eq:large-dev-of-omega-star}. We can now
apply that inequality to conclude that for some constants $c_4>0$,
\begin{align*}
\Pp\{\Sigma \le c_3 xn,\ A^n\le -xn\}\le
 \Pp\{F^*_{\omega, n,\lceil c_3 xn\rceil}\ge xn\}
                                       \le e^{-c_4 xn}.
\end{align*}
The second term on the right-hand side of~\eqref{eq:estimating_action_linear_tails} can be estimated using
Lemma~\ref{lem:distribution-of-animal-size}. 
If $c_3x\geq R$ and $n$ is sufficiently large, then
\begin{align*}
 \Pp\{\Sigma > c_3 xn\}&\le \sum_{m\ge c_3xn} \Pp(E_{n,m}) \\
                      &\le \sum_{m\ge c_3xn}  C_1\exp(-C_2m^2/n)\\
                      &\le C'_1\exp(-C'_2x^2n),
\end{align*}
for some constants $C'_1,C'_2>0$,
which completes the proof.\epf

\bpf[Proof of uniqueness in Theorem~\ref{thm:global_solutions}] 
Let $w_\omega=W'_\omega$ be a global solution of the Burgers equation
such that for each $n\in\Z$, $w(n,\cdot)\in\HH'(v,v)$ with probability~$1$. 
Then, for any $n\le N$,  $w(0,\cdot)=\Psi_\omega^{-n,0} w(-n,\cdot)$. 
The cadlag version of $w$ belongs to $\GG(v,v)$. For any $x\in\R$, 
the trajectory solving the Euler--Lagrange equation and terminating at $x$
with velocity $w(x)$ is a minimizer on every finite interval. Therefore it is a one-sided minimizer and must have an 
asymptotic slope~$\tilde v(x)$. Notice that $\tilde v(x)$ is monotone in $x$, since on any finite interval the
minimizers cannot intersect. Due to spatial translation invariance, $\tilde v(x)$ is a stationary process in $x$, so 
$\tilde v(x)=\tilde v$ has to be a constant. Hence, $w_\omega(0,\cdot)$ almost surely coincides with $u_{\tilde v}(0,\cdot)$.
Since the latter belongs to $\GG(\tilde v,\tilde v)$ almost surely, we see that $\tilde v=v$, so
$w_\omega(0,\cdot)$ almost surely coincides with $u_{v}(0,\cdot)$, which completes the proof.\epf

\section{Metric on $\GG$}\label{sec:metric}
Let us recall that $\GG$ consists of all cadlag functions $w:\R\to\R$ such that $M_w:\R\to\R$ defined by $M_w(x)=x-w(x)$
is a strictly increasing function satisfying $\lim_{x\to\pm\infty}M_w(x)=\pm\infty$. The goal of this
section is to introduce a metric $d$ on $\GG$ such that $\lim_{n\to\infty} d(w_n,w)= 0$ is equivalent to
$\lim_{n\to\infty}w_n(x)=w(x)$ for all $x\in\Cc(w)$.
First we note that for every $w$, the inverse $M^{-1}_w$ of $M_w$ defined by
\[
M^{-1}_w(y)=\inf\{x:M_w(x)\ge y\},\quad x\in\R,
\]
is a  continuous function retaining all the information about $M_w$ and $w$.
Let us define $d$ as the metric of locally uniform convergence on continuous functions on $\R$:
\[
d(u,w)=\sum_{N=1}^{\infty} 2^{-N}d_N(u,w),\quad u,v\in\GG,
\]
where
\[
d_N(u,w)=\sup_{x\in[-N,N]} |M^{-1}_u(x)-M^{-1}_w(x)|\wedge 1,\quad u,v\in\GG.
\]

\begin{lemma} Let $(w_n)_{n\in\N}$ be a sequence in $\GG$ and $w\in\GG$. Then $d(w_n,w)\to 0$ as $n\to\infty,$ iff
$\lim_{n\to\infty}w_n(x)=w(x)$ for all $x\in\Cc(w)$.
\end{lemma}
\bpf Suppose $d(w_n,w)\to 0$ as $n\to\infty$. We need to prove $M_{w_n}(x)\to M_{w}(x)$ for each $x\in\Cc(w)=\Cc(M_w)$. 
So let us take such an $x$ and any $\eps>0$. We can find $x^-<x$ and $x^+>x$
such that \[M_w(x)-\eps < M_w(x^-)< M_w(x)<M_w(x^+)<M_w(x)+\eps.\]
Let us denote $y^-=M_w(x^-)$, $y=M_w(x)$, $y^+=M_w(x^+)$.
Since $d(w_n,w)\to 0$, we see that there is $n_0$ such that for $n>n_0$, $M_{w_n}^{-1}(y^-)<(x^-+x)/2$ and $M_{w_n}^{-1}(y^+)>(x+x^+)/2$.
Therefore, for $n>n_0$, 
$y^-<M_{w_n}(x)<y^+$, so $M_w(x)-\eps < M_{w_n}(x)<M_w(x)+\eps$, and our claim follows.

Now let us assume that $\lim_{n\to\infty}w_n(x)=w(x)$ for all $x\in\Cc(w)$. We need to show that for any $N\in\N$, $d_N(w_n,w)\to0$.
First, we can find points $z^-,z^+\in \Cc(w)$ such that $M_w(z^-)<-N$ and $M_w(z^+)>N$. There is $n_0$ such that $[-N,N]\subset 
(M_{w_n}(z^-),M_{w_n}(z^+))$ for all $n>n_0$.

For any $\eps>0$ let us find a finite collection of points $x_0,x_1,\ldots,x_m\in \Cc(w)$ such  that $z^-=x_0<x_1<\ldots<x_m=z^+$ and
 $x_{k}-x_{k-1}<\eps/2$ for all $k=1,\ldots,m$. Let us denote $y_k=M_w(x_k),$ $k=0,1,\ldots,m$. Points $y_k=M_w(x_k),$ $k=0,1,\ldots,m$
 form a strictly increasing sequence. Let us denote 
 \[\Delta=\min_{k=1,\ldots,m} (y_k-y_{k-1})\wedge (-N-y_0)\wedge (y_m-N)>0.\] Since $x_k\in\Cc(w)$, we can find $n_1\ge n_0$
 such that for all $n>n_1$ and all $k=0,\ldots,m$, $|M_{w_n}(x_k)-y_k|< \Delta/2$.
Therefore, if $y\in [y_{k-1},y_k]\cap [-N,N]$, then
$M_{w_n}^{-1}(y)\in [x_{(k-2)\vee 0},x_{(k+1)\wedge m}]$ and $M_{w}^{-1}(y)\in [x_{k-1},x_{k}]$. Hence, $|M_{w_n}^{-1}(y)-M_{w}^{-1}(y)|\le \eps$. Since $[-N,N]\subset \bigcup_{k=1}^m [y_{k-1},y_k]$,  we obtain $d_N(w_n,w)\le \eps$ for $n\ge n_1$. This is the desired uniform estimate.
 \epf

\section{Auxiliary lemmas}\label{sec:aux}
\bpf[Proof of Lemma~\ref{lem:invariant_spaces}]   The local Lipschitzness of $\Phi^{n_0,n_1}_\omega W$ is part~\ref{it:continuity_of_HJ_solution} of Lemma~\ref{lem:properties-of-evolution-operator},
so let us establish the behavior as $x\to\infty$. 
Let us begin with the second part of the Lemma. 

Due to the cocycle property it is sufficient to consider the situation where $n_1=n_0+1$.
Let us take $W\in\HH(v_-,v_+)$.
We have
\begin{equation*}
\Phi^{n_0,n_0+1}_\omega W(x)\le W(x)+F_{\omega}(n_0,x), \quad x\in\R. 
%\label{eq:lower-bound-on-cocycle}
\end{equation*}
Since $\lim_{x\to\infty}(F_{\omega}(n_0,x)/x)=0$ on $\Omega_1$, we see that
\[
\limsup_{x\to+\infty} \frac{\Phi^{n_0,n_0+1}_\omega W(x)}{x}\le  v_{+}.
\]
Let us prove that 
\[
\liminf_{x\to+\infty} \frac{\Phi^{n_0,n_0+1}_\omega W(x)}{x}\ge v_{+}.
\]
If this inequality is violated, then there is $\eps>0$ and two increasing sequences $(x_k)_{k\in\N}$, $(y_k)_{k\in\N}$ such that
$x_k\to+\infty$ and
\begin{equation}
\label{eq:proving-Hvv}
V(y_k)+\frac{1}{2}(x_k-y_k)^2< (v_+-\eps)x_k,\quad k\in\N
\end{equation}
where $V(y)=W(y)+F(y)$, $y\in\R$. Inequality~\eqref{eq:proving-Hvv} implies that $y_k$ cannot be bounded, so we obtain $y_k\to+\infty$.
Moreover, it follows from~\eqref{eq:proving-Hvv} and from $\lim_{y\to\infty}V(y)/y=v_+$ that for sufficiently large $k$,
\[
v_+ y_k < (v_+-\eps/2)x_k,
\]
so $y_k/x_k$ cannot converge to $1$. Therefore, there is a subsequence $(k')$ and a constant $c>0$ 
such that $|y_{k'}-x_{k'}|>c(y_{k'}+x_{k'})$, so
\begin{align*}
V(y_{k'})&\le -\frac{1}{2}(x_{k'}-y_{k'})^2+ (v_+-\eps)x_{k'}\\
      &\le - \frac{c^2y_{k'}^2}{2}-\frac{c^2x_{k'}^2}{2}+ (v_+-\eps)x_{k'}\\
      &\le - \frac{c^2y_{k'}^2}{2},
\end{align*}
for sufficiently large $k'$ which contradicts the linear growth of $V$. A similar analysis applies to the behavior at $-\infty$, and the second
part of the Lemma is proved completely.

To prove the first part, assume that $W\in\HH$, but $\Phi^{n_0,n_0+1}_\omega W\notin\HH$ due to the behavior, say, at $+\infty$. Then 
there are two increasing sequences $(x_k)_{k\in\N}$, $(y_k)_{k\in\N}$ such that
$x_k\to+\infty$ and
\begin{equation}
\label{eq:proving-Hvv-1}
\frac{V(y_k)+\frac{1}{2}(x_k-y_k)^2}{x_k}\to -\infty,\quad k\to\infty.
\end{equation}
This means that $V(y_k)/x_k\to-\infty$. Since $V(y_k)/y_k$ is bounded below for large $k$, we conclude that $x_k/y_k\to 0$. Now
\eqref{eq:proving-Hvv-1} implies that $V(y_k)\le - y_k^2/4$ for large $k$, which contradicts $W\in\HH$. Similar reasoning applies to the 
behavior near $-\infty$.\epf

\begin{lemma} \label{lem:limit-of-minimizers} A pointwise limit of a sequence of point-to-point minimizers is a point-to-point minimizer.
\end{lemma}
\bpf Suppose $(\gamma^k)_{k\in\N}$ is a sequence of point-to-point minimizers on a time interval $\{n_1,\ldots,n_2\}$. Suppose that
$\gamma$ is a path such that $\gamma^k_j\to \gamma_j$ as $k\to\infty$ for $j\in\{n_1,\ldots,n_2\}$. If $\gamma$ is not a 
point-to-point minimizer, then there is a path $\beta$ satisfying $\beta_{n_1}=\gamma_{n_1}$, $\beta_{n_2}=\gamma_{n_2}$,
$A^{n_1,n_2}(\beta)<A^{n_1,n_2}(\gamma)$. 
Let us introduce paths
 $\beta^k=(\gamma^k_{n_1},\beta_{n_1+1},\ldots,\beta_{n_2-1},\gamma^k_{n_2})$.
Since $A^{n_1,n_2}(\gamma^k)\to A^{n_1,n_2}(\gamma)$ and  $A^{n_1,n_2}(\beta^k)\to A^{n_1,n_2}(\beta)$, we obtain
$A^{n_1,n_2}(\beta^k)<A^{n_1,n_2}(\gamma^k)$ for sufficiently large $k$ which contradicts the minimizing property of $\gamma^k$.
\epf

\begin{lemma}
\label{lem:finite-horizon-minimizers-do-not-intersect}
 Let $\omega\in\Omega$ and $(n,x)\in\ZR$. If $\gamma^1$ and $\gamma^2$ are two distinct point-to-point
 minimizers on $\{n,\ldots,n'\}$ satisfying $\gamma^1_n=\gamma^2_n=x$, then, as curves in $\R\times\R$,
 they do not intersect on time interval $(n,n')$.  
\end{lemma}
\bpf[Proof of Lemma~\ref{lem:finite-horizon-minimizers-do-not-intersect}] If $\gamma^1$ and $\gamma^2$ have two consecutive points in common, they coincide due to the Euler--Lagrange equation. 

Suppose  $\gamma^1_m=\gamma^2_m=y$ and $\gamma^1_{m-1}<\gamma^2_{m-1}$ for some $m\in\{n+1,\ldots,n'-1\}$. By the Euler--Lagrange equation, $\gamma^1_{m+1}>\gamma^2_{m+1}$.
For $\delta>0$, we define 
\begin{align*}
\tilde\gamma^1&=(x,\gamma^1_{n+1},\ldots,\gamma^1_{m-1},y-\delta,\gamma^2_{m+1}),\\
\tilde\gamma^2&=(x,\gamma^2_{n+1},\ldots,\gamma^2_{m-1},y+\delta,\gamma^1_{m+1}).
\end{align*}
Then
\begin{align*}
&A^{n,m+1}(\tilde\gamma^1)+A^{n,m+1}(\tilde\gamma^2) -(A^{n,m+1}(\gamma^1)+A^{n,m+1}(\gamma^2)) 
\\=& F(m,y+\delta)+F(m,y-\delta)-2F(m,y)
\\&+\frac{1}{2}\Bigl((\gamma^1_{m+1}-y-\delta)^2+(\gamma^2_{m+1}-y+\delta)^2+(y+\delta-\gamma^2_{m-1})^2+(y-\delta-\gamma^1_{m-1})^2
\\&-(\gamma^1_{m+1}-y)^2-(\gamma^2_{m+1}-y)^2-(y-\gamma^2_{m-1})^2-(y-\gamma^1_{m-1})^2\Bigr)
\\=& F(m,y+\delta)+F(m,y-\delta)-2F(m,y)+2\delta^2+\delta(\gamma^2_{m+1}-\gamma^1_{m+1}+\gamma^1_{m-1}-\gamma^2_{m-1}).
\end{align*}
Since $\delta(\gamma^2_{m+1}-\gamma^1_{m+1}+\gamma^1_{m-1}-\gamma^2_{m-1})<0$ and the sum of the remaining terms is $o(\delta)$, $\delta\to0$ due to
differentiability of $F$, we obtain that for sufficiently small $\delta$,
\[
A^{n,m+1}(\tilde\gamma^1)+A^{n,m+1}(\tilde\gamma^2) < A^{n,m+1}(\gamma^1)+A^{n,m+1}(\gamma^2),
\]
so either $A^{n,m+1}(\tilde\gamma^1)<A^{n,m+1}(\gamma^2)$ or $A^{n,m+1}(\tilde\gamma^2)<A^{n,m+1}(\gamma^1)$ which contradicts the minimizing property
of $\gamma^1,\gamma^2$.

It remains to exclude the case where $\gamma^1_{m}<\gamma^2_m$ and $\gamma^1_{m+1}>\gamma^2_{m+1}$ for
some $m\in\{n+1,\ldots,n-1\}$. In this case, we define
\begin{align*}
\tilde\gamma^1&=(x,\gamma^1_{n+1},\ldots,\gamma^1_{m},\gamma^2_{m+1}),\\
\tilde\gamma^2&=(x,\gamma^2_{n+1},\ldots,\gamma^2_{m},\gamma^1_{m+1}).
\end{align*}
Then 
\begin{align*}
&A^{n,m+1}(\tilde\gamma^1)+A^{n,m+1}(\tilde\gamma^2) -(A^{n,m+1}(\gamma^1)+A^{n,m+1}(\gamma^2)) 
\\=& \frac{1}{2}\Bigl((\gamma^2_{m+1}-\gamma^1_{m})^2+(\gamma^1_{m+1}-\gamma^2_{m})^2-(\gamma^1_{m+1}-\gamma^1_{m})^2-(\gamma^2_{m+1}-\gamma^2_{m})^2\Bigr)
\\=& (\gamma^2_m-\gamma^1_m)(\gamma^2_{m+1}-\gamma^1_{m+1})<0,
\end{align*}
so, as above,  at least one of the paths $\gamma^1,\gamma^2$ is not optimal. This contradiction completes the proof of the lemma.
\epf

\bibliographystyle{alpha}
\bibliography{Burgers}
\end{document}